\documentclass[11pt]{article}

\usepackage{paper/style}

\numberwithin{equation}{section}
\numberwithin{figure}{section}

\newcommand{\titletext}{Simultaneous Frequentist Calibration of Confidence Regions for Multiple Functionals in Constrained Inverse Problems}

\title{\titletext}

\author{
    Pau Batlle\footremember{caltechcms}{Department of Computing and Mathematical Sciences, California Institute of Technology, CA 91125, USA.} 
    \\ {\small \texttt{\href{mailto:pbatllef@caltech.edu}{pbatllef@caltech.edu}}} 
    \and ~~~Pratik Patil\footremember{utstats}{Department of Statistics and Data Sciences, University of Texas, Austin, TX 78712, USA.} 
    \\ {\small ~~~\texttt{\href{mailto:pratikpatil@utexas.edu}{pratikpatil@utexas.edu}}} 
    \and Michael Stanley\footremember{ama}{Analytical Mechanics Associates, Hampton, VA 23666, USA.} 
    \\ {\small \texttt{\href{mailto:michael.c.stanley@ama-inc.com}{michael.c.stanley@ama-inc.com}}} \\ 
    \and Javier Ruiz Lupon\footrecall{caltechcms}~~~~ 
    \\ {\small \texttt{\href{mailto:jruizlup@caltech.edu}{jruizlup@caltech.edu}}~~~~} 
    \and Houman Owhadi\footrecall{caltechcms}~~~~
    \\ {\small \texttt{\href{mailto:owhadi@caltech.edu}{owhadi@caltech.edu}}~~~~} 
    \and Mikael Kuusela\footremember{cmustats}{Department of Statistics and Data Science, Carnegie Mellon University, Pittsburgh, PA 15213, USA.}~~~
    \\ {\small \texttt{\href{mailto:mkuusela@andrew.cmu.edu}{mkuusela@andrew.cmu.edu}}~~~}
}

\date{\vspace{-10pt}}

\begin{document}

\maketitle

\begin{abstract}
    Many scientific analyses require simultaneous comparison of multiple functionals of an unknown signal at once, calling for multidimensional confidence regions with guaranteed simultaneous frequentist under structural constraints (e.g., non-negativity, shape, or physics-based). This paper unifies and extends many previous optimization-based approaches to constrained confidence region construction in linear inverse problems through the lens of statistical test inversion. We begin by reviewing the historical development of optimization-based confidence intervals for the single-functional setting, from ``strict bounds'' to the Burrus conjecture and its recent refutation via the aforementioned test inversion framework. We then extend this framework to the multiple-functional setting. This framework can be used to: (i) improve the calibration constants of previous methods, yielding smaller confidence regions that still preserve frequentist coverage, (ii) obtain tractable multidimensional confidence regions that need not be hyper-rectangles to better capture functional dependence structure, and (iii) generalize beyond Gaussian error distributions to generic log-concave error distributions. We provide theory establishing nominal simultaneous coverage of our methods and show quantitative volume improvements relative to prior approaches using numerical experiments.
\end{abstract}

\newpage

\section{Introduction}
\label{sec:introduction}

We study linear inverse problems of the form
\begin{equation}
    \label{eq:linear-model}
    y = K x^* + \varepsilon,
\end{equation}
where $y \in \mathbb{R}^n$ denotes the observations, $x^* \in \mathbb{R}^p$ is a fixed unknown parameter vector, $K \in \mathbb{R}^{n \times p}$ is a known arbitrary forward operator, and the random noise vector $\varepsilon \in \mathbb{R}^n$ is drawn from a known log-concave density (with respect to the Lebesgue measure).
To connect with and extend previous work, we begin with the Gaussian setting $\varepsilon \sim \mathcal{N}(0, I)$\footnote{Assuming zero mean and identity covariance is without loss of generality after shifting and whitening.} and later generalize to generic log-concave settings. 
The true parameter $x^*$ is known to satisfy certain constraints based on prior physical knowledge. 
In this work, we focus on the historically important example of non-negative constraints $x^* \geq 0$, and then extend to two broader classes: (i) general linear constraints $Ax^*\leq b$ (including box-constraints $x_{\mathrm{lo}} \leq x^* \leq x_{\mathrm{up}}$), and (ii) cone constraints $x^* \in \mathcal{C}$ for a closed convex cone $\mathcal{C}$. 
Let $\mathcal{X}$ denote the set of parameters satisfying our constraints, and for every $x^* \in \mathcal{X}$, write $P_{x^*}$ for the distribution of $y$ under the model \eqref{eq:linear-model} with true parameter $x^*$ (e.g., $P_{x^*} = \mathcal{N}(Kx^{*}, I)$ in the Gaussian setting).

We aim not to estimate the (possibly) high-dimensional vector $x^*$ itself, but rather to make inferences about a few key linear functionals of it, such as localized averages or weighted sums over specific components. 
Formally, given a matrix $H \in \mathbb{R}^{k \times p}$, we are interested in determining a finite-sample $1-\alpha$ frequentist confidence set for the unknown vector $Hx^* \in \mathbb R^k$, i.e., a region $\mathcal{R}_\alpha(y) \subseteq \mathbb{R}^k$ such that, for any $x \in \mathcal{X}$,
\begin{equation}
    \label{eq:freq_guarantees}
    \mathbb{P}_{y \sim P_{x}}\big(Hx \in \mathcal{R}_\alpha(y)) \ge 1-\alpha.
\end{equation} 
We stress that the guarantee \eqref{eq:freq_guarantees} needs to hold over all parameters $x$ in the constrained parameter space $\cX$, so that in particular holds for the true unknown parameter $x^*$.
In the literature, this problem is sometimes formulated by identifying $Hx^* = (h_1^\top x ^*, \dots, h^\top_k x^*) \in \mathbb{R}^k$ as the image of $k$ functionals of interest and then finding $k$ intervals (or one-dimensional regions) $\mathcal{R}^i_\alpha(y), i = 1, \dots, k$. 
Under this framework, we emphasize that \eqref{eq:freq_guarantees} specifies \emph{simultaneous} coverage, also named \emph{$k$-at-a-time}: 
\begin{equation}
    \mathbb{P}_{y \sim P_{x^*}} \big(h^\top_i x^* \in \mathcal{R}^i_\alpha(y) \text{ for all } i = 1, \dots, k) \ge 1-\alpha,
\end{equation}
as opposed to \emph{marginal} or \emph{one-at-a-time} coverage:
\begin{equation}
    \mathbb{P}_{y \sim P_{x^*}} \big(h^\top_i x^* \in \mathcal{R}^i_\alpha(y)) \ge 1-\alpha \text{ for all } i = 1, \dots, k.
\end{equation}
In this work, we do not restrict the region $\mathcal R_\alpha(y)$ to be necessarily a product of one-dimensional intervals (or regions). 
We allow geometry that captures dependencies among functionals (e.g., polytopes or ellipsoids), provided the region is computationally tractable. 
We define this for this work as being able to test membership $\mu \in \mathcal{R}_\alpha(y)$ by solving a convex optimization problem whose size scales at most linearly with the problem dimensions.

\begin{figure}[!t]
    \centering
    \includegraphics[width = 0.8\columnwidth]{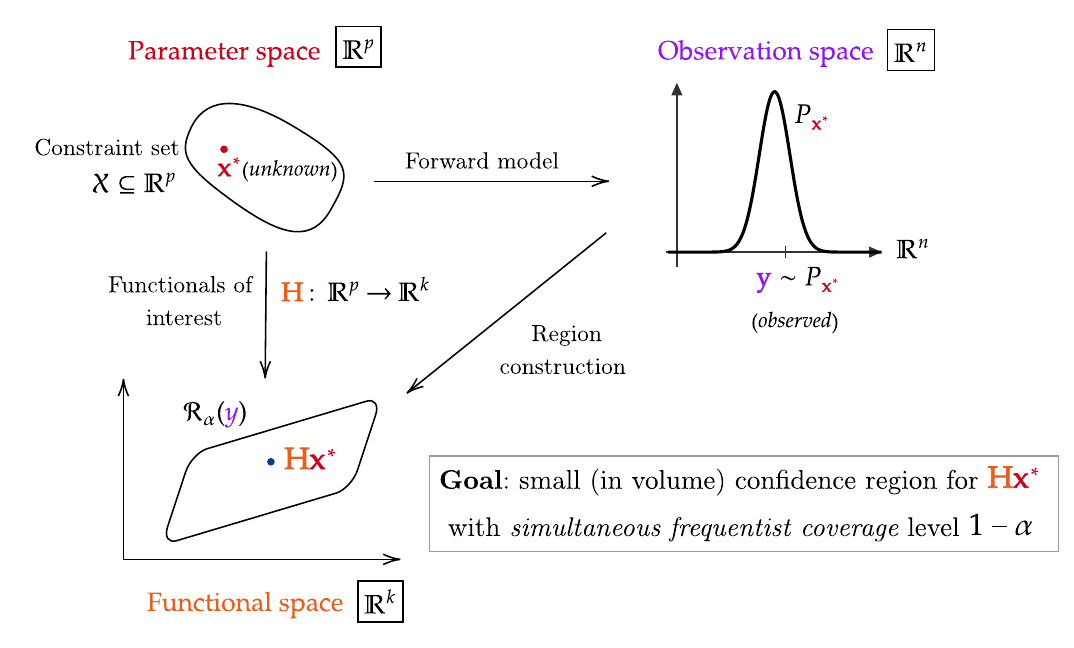}
    \caption{
    Illustration of the problem setup. 
    We seek to construct confidence region $\cR_{\alpha}(y) \subseteq \RR^k$ for $H x^* \in \RR^k$ from an observation $\by \in \RR^n$ sampled from $P_{\bx^*}$ that satisfies a frequentist coverage guarantee in finite sample while being as small (in volume) as possible.
    }
    \label{fig:setup-multiple-functionals}
\end{figure}

We generally compare and evaluate competing methods by the size of their regions $\cR_{\alpha}(y)$ (e.g., the $k$-dimensional Lebesgue measure) subject to the simultaneous coverage guarantee \eqref{eq:freq_guarantees}, since arbitrarily large regions always satisfy such a guarantee.
However, this is a challenging problem as the size of $\cR_{\alpha}(y)$ is a random quantity whose distribution depends on the unknown $x^*$. Consequently, different optimality notions (such as minimax or average expected size or tail quantiles) lead to distinct and generally challenging problems, and no uniformly best method exists across all $x^*$.
Instead, in this work, we focus on establishing provable \textit{uniform} improvements over popular existing methods: we construct procedures producing improved regions $\widetilde{\cR}_{\alpha}(y)$ such that: (i) we retain the same finite-sample simultaneous coverage as in \eqref{eq:freq_guarantees}, yet (ii) we satisfy $\widetilde{\cR}_{\alpha}(y) \subseteq \cR_{\alpha}(y)$ for every data realization $y$, where $\cR_{\alpha}(y)$ is a region given by a previous method with guarantees. 
We make no claim of global optimality among all conceivable methods; rather, we demonstrate strict dominance (by set inclusion) over commonly used prior methods. 

Our main baseline throughout this article is the classical Simultaneous Strict Bounds (SSB) method, which achieves the guarantee \eqref{eq:freq_guarantees} by first constructing a $1-\alpha$ confidence set for the full true parameter $x^*$, intersecting it with the constraint set $\mathcal{X}$, and then mapping the intersection through the functional matrix $H$ into the functional space $\mathbb{R}^k$ \cite{Stark1992}. While valid, this approach is often conservative and yields regions that are larger than necessary because the initial confidence set for $x^*$ must be valid for \emph{any} possible functional (including nonlinear ones), not just the specific linear functionals in $H$.

The key tool we use to improve upon this method is a link between two classical ways for constructing confidence regions: statistical test inversion and optimization-based methods, which we leverage to shrink regions while preserving coverage.

We briefly situate our contributions within prior work next.
For a broader historical context, including the existing connections between test-inversion and optimization-based viewpoints for inference on a single functional, see \Cref{sec:history}.

\subsection{Related work}

The problem of constructing confidence regions with frequentist guarantees has a long history (with additional context provided in \Cref{sec:history}). 

Under weaker model assumptions than the aforementioned SSB method but without constraints, classical statistics has similarly extensively developed simultaneous confidence regions \cite{miller_simultaneous, seberlee}.
This in turn has motivated a long line of research aimed at finding tighter ``one-at-a-time'' intervals, most famously captured by the Burrus conjecture \cite{Burrus1964, Rust1972-cf, rust1994confidence, tenorio2007confidence}.
This conjecture, which proposed a less conservative calibration for single-functional intervals under non-negativity constraints, remained unresolved for decades before being disproven recently through a modern framework based on the inversion of a likelihood ratio test (LRT) \cite{BatlleDisproof}.
This approach lead to a large variety of one-at-a-time methods with one computational approach developed in \cite{stanley2025confidenceintervalsfunctionalsconstrained}.

The literature on constrained inference similarly develops theory and methods to incorporate known parameter constraints into hypothesis testing and estimation and is thoroughly summarized in \cite{robertson1988, Silvapulle2001}.
These approaches heavily rely upon particular hypothesis test frameworks allowing for the use of the chi-bar-squared distribution.
The hypothesis test forms explored in this work and \cite{BatlleDisproof, stanley_unfolding} fall outside those classical frameworks, though we discuss connections in \Cref{sec:calibration_problem}.

Other statistical inference approaches that rely on splitting or resampling data can also be used for constrained inference.
Universal inference \cite{universal_inference} is one such method relying on data splitting.
Bootstrapping approaches can be used to resample constrained estimators and achieve confidence sets \cite{li_constrained_bootstrap} but the validity of these approaches is typically asymptotic.
Additionally, we emphasize that the scientific scenarios motivating this work often prohibit the use of data splitting or resampling, since the forward model defining $K$ is often only computationally available \cite{stanley_thesis}.

The proposed framework shares some parallels with conformal prediction (CP) (see, e.g., \cite{vovk2005algorithmic,lei2014distribution,lei2018distribution} for some early work and \cite{angelopoulos2023gentle} for a recent overview of some recent trends), which also typically constructs finite-sample valid regions via test inversion and quantile calibration.
The key differences are in goals, guarantees, and assumptions.
We target confidence regions for a collection of unknown parameter functionals, whereas conformal (and in general predictive) inference methods typically target prediction intervals for future responses.
CP assumes exchangeability of the data and provides marginal coverage with respect to the underlying distribution.
By contrast, we assume a (possibly nonlinear) forward model with Gaussian or log-concave noise and structural constraints on the parameter.
Under these assumptions, we obtain stronger (uniform worst-case frequentist) coverage guarantees that hold for all feasible parameters in the constraint set, here extended to simultaneous coverage across multiple functionals.

\subsection{Summary of contributions}

Broadly speaking, we construct and analyze (via test inversion) simultaneous confidence regions that retain finite-sample frequentist coverage while being uniformly less conservative than classical SSB.
More specifically:
\begin{itemize}
    \item We synthesize the extensive literature on single-functional methods, from classical strict bounds to modern optimization-based procedures, within a unified test-inversion language (\Cref{sec:unconstrained_case,sec:history}).
    \item We extend this framework from single to multiple functionals, introducing joint test statistics and their associated calibration problems; for several choices of test statistics, we derive optimal calibration constants (\Cref{sec:test_inversion,sec:test_statistics}).
    \item We propose computationally practical procedures for constructing these resulting regions, including algorithms to compute the aforementioned optimal constants and reduction and splitting techniques for high-dimensional settings (\Cref{sec:high-dimensional}).  
\end{itemize}

The paper proceeds as follows.
\Cref{sec:unconstrained_case} reviews the classical unconstrained problem as a warm-up and shows the equivalence between parameter-space and functional-space formulations and casts the standard Wald ellipsoid via test inversion.
\Cref{sec:history} gives historical context for the single-functional, non-negatively constrained case, tracing the evolution of the Burrus conjecture and its eventual resolution, which motivates the modern test-inversion framework.
In \Cref{sec:test_inversion}, we formally extend this framework to the multiple-functional case, formulating joint hypothesis tests and the construction of simultaneous confidence regions.
\Cref{sec:test_statistics} analyses the calibration problem in details and compares three key test statistics identifying their theoretical properties, which are crucial for practical implementation.
\Cref{sec:high-dimensional} presents computationally practical algorithms for high-dimensional problems, including generalized reductions for multiple functionals or and/box-constrained problems and a new row--null splitting technique for rank-deficient models.
\Cref{sec:numerical} provides numerical experiments demonstrating nominal coverage for our methods with substantially smaller (in volume) confidence regions compared to classical approaches.
We conclude with a summary and discussion of some open problems in \Cref{sec:conclusion}.

\section{Warming up with the unconstrained case}
\label{sec:unconstrained_case}

We begin by reviewing the classical problem of constructing a confidence region for $H x^*$ given $y = Kx^* + \varepsilon$ with $\varepsilon \sim \mathcal N(0,I)$ in the unconstrained setting, where the parameter space is $\cX = \RR^p$.
This well-understood case serves as an essential baseline and allows us to introduce the key recurring concepts of observability of the functionals, parameter description, and functional description of confidence sets, and the test-inversion viewpoint that we will generalize to the constrained setting.
We then preview how adding constraints can (in a sense) ``create'' observability that can result in finite regions even when unconstrained regions would be unbounded.

\subsection{Classical unconstrained framework}

In the unconstrained setting, we first consider an idealized ``oracle'' case with noiseless observations, which leads to a fundamental condition for obtaining bounded confidence regions.
We then introduce noise and review the classical Wald ellipsoid, its relationship to parameter-space representations, and its modern interpretation through test inversion.

Define the oracle compatibility region $\mathcal{R}^O$ as the set of functionals $\mu \in \mathbb{R}^k$ (i.e., values of $Hx^* = (h_1^\top x ^*, \dots, h^\top_k x^*)$) that would be compatible with our data $y$ if that data was observed without noise (i.e., $y = K x^*$):
\begin{equation} \label{eq:og_mu_description}
    \mathcal{R}^O = \{\mu \in \mathbb{R}^k: \exists \, x \in \mathbb{R}^p \text{ such that } Kx = y, Hx = \mu\}.
\end{equation}
A simple argument shows that 
\begin{equation} \label{eq:og_x_description}
    \mathcal{R}^O = \{Hx: x \in \mathbb{R}^p, Kx = y\} = \{HK^\dagger y + H(I-K^\dagger K)z : z \in \mathbb{R}^p\}.
\end{equation}
This simple equality between a set in $\mathbb{R}^k$ (functional space) and a linear image, under $H$, of a set in $\mathbb{R}^p$ (parameter space) will have strong consequences in noisy constrained settings. 
We will henceforth refer to descriptions such as \eqref{eq:og_mu_description}, directly in terms of elements of the functional space, as the $\mu$-description of the set, and to descriptions such as \eqref{eq:og_x_description}, in terms of the linear image of elements in the parameter space, as the equivalent $x$-description of the set, and their equivalence will be an important tool for the analysis of noisy and/or constrained settings. 

In this simple noiseless unconstrained setting, any functional of interest $h_k$ can be either exactly recovered by $h_k^\top K^\dagger y$ if $h_k \in \row(K)$ or cannot be estimated from the data, so that $\mathcal{R}^O$ is the coordinate-wise product of singletons and infinite intervals. 
Overall, a finite volume region for $Hx^*$ is possible if and only if all of the functionals of interest (rows of $H$) lie in $\row(K)$ (the same observability condition holds with noise as well). 
A general, not necessarily bounded, confidence region can be constructed by taking the product of the bounded confidence region of the functionals of the row space and copies of $\mathbb{R}$ for those unobservable functionals (i.e., not in the row space).
We henceforth assume that $H$ lies in $\row(K)$ to build the confidence region for observable functionals. 
We furthermore assume that the rows of $H$ are linearly independent without loss of generality, since if a functional lies in the span of the other rows of $H$, its inference follows directly by projection of the confidence region of the other functionals.

Under this assumption, there exists a matrix $B$ such that $B^\top K = H$. With Gaussian noise $\varepsilon \sim \cN(0, I)$, we have $B^\top y = B^\top Kx^* + B^\top \varepsilon \sim \mathcal{N}(Hx^*, B^\top B)$. 
The matrix satisfying this that minimizes the variances in each column is $B = K(K^\top K)^\dagger H^\top$. Since $B^\top y \sim \mathcal{N}(Hx^*, B^\top B)$, one obtains the following $1-\alpha$ ellipsoidal region:
\begin{equation}
    \label{eq:mu_space_unc}
    \mathcal{R}_y^\mu = \{\mu \in \mathbb{R}^k: \|\mu - B^\top y\|^2_{B^\top B} \leq Q_{\chi^2_k,1-\alpha}\},
\end{equation}
where we adopt the convention $\|v\|^2_{B^\top B} := v^\top (B^\top B)^{-1}v$\footnote{We will use this notation as well as $\|v\|^2_M = v^\top M^\dagger v$ for the case in which $M$ is not invertible} and $Q_{X,1-\alpha}$ denotes the $(1-\alpha)$-quantile of a random variable $X$ (satisfying $\mathbb{P}_{x\sim X}\left(x \leq Q_{X,1 - \alpha} \right) = 1 - \alpha$), with a slight abuse of notation, we will also use $Q_{P,1-\alpha}$ to denote the $(1-\alpha)$-quantile of a distribution $P$.

Since $B^\top y = H (K^\top K)^\dagger K^\top y = H K^\dagger y$, we can view \eqref{eq:mu_space_unc} as the ellipsoid around the (minimum $\ell_2$-norm) least squares estimator of $Hx^*$, which is known as the Wald form.

As in the noiseless case, we can identify \eqref{eq:mu_space_unc} as the linear image under $H$ of a set of \say{compatible parameters} in the original parameter space. 
Defining 
\begin{equation}
    \label{eq:x_space_unc}
    \mathcal{R}_{y}^x = \{Hx: x \in \mathbb{R}^p, \|Kx-y\|^2_2 \leq \min_{x'}\|Kx'-y\|^2_2 + Q_{\chi^2_k,1-\alpha}\},
\end{equation}
we have $\mathcal{R}_y^\mu = \mathcal{R}^x_y$ (equivalence between the $\mu$-description \eqref{eq:mu_space_unc} and the $x$-description \eqref{eq:x_space_unc}). 

We prove a useful generalization of this result in the following proposition, in terms of a function of the data $f(y)$ and a real-valued constant $M$, which we separate from the function $f$ to facilitate conceptual understanding of the results to follow.
\begin{proposition}[Equivalence between $x$-description and $\mu$-description]
    \label{prop: mu_x_equiv}
    Let $f: \mathbb{R}^n \to \mathbb{R}$ and $M \in \RR$, and define the set 
    \begin{equation}
        \mathcal{S}_y = \{x \in \mathbb{R}^p : \|Kx-y\|^2_2 \leq f(y) + M \}.
    \end{equation}
    Then, we have
    \begin{align}
        \{Hx : x \in \mathcal{S}_y \} 
        &= \{\mu \in \RR^k : \inf_{Hx = \mu} \|Kx-y\|^2_2 \leq f(y) + M\} \label{eq:three_descriptions_implicit} \\
        &=  \{\mu \in \RR^k : \|\mu - B^\top y\|^2_{B^\top B} \leq f(y) + M - \inf_{x'}\|Kx'-y\|^2_2  \}. \label{eq:three_derscriptions_explicit}
    \end{align}
\end{proposition}

\begin{proof}
    Let $c(y):=\min_{x'}\|Kx'-y\|_2^2=\|(I-KK^\dagger)y\|_2^2$ and $\phi(\mu):=\min_{Hx=\mu}\|Kx-y\|_2^2$. 
    Since $H\subseteq\row(K)$, there is $B=K(K^\top K)^\dagger H^\top$ with $B^\top K=H$. 
    A standard constrained LS argument gives $\phi(\mu)=c(y)+\|\mu-B^\top y\|_{B^TB}$.
    
    ($\subseteq$) If $x\in\mathcal{S}_y=\{x:\|Kx-y\|_2^2\le f(y) + M\}$ and $\mu=Hx$, then $\phi(\mu)\le\|Kx-y\|_2^2\le f(y)+M$, hence $\mu$ belongs to the RHS.
    
    ($\supseteq$) If $\phi(\mu)\le f(y) + M$, pick $x_\mu$ attaining $\phi(\mu)$. 
    Then, $\|Kx_\mu-y\|_2^2=\phi(\mu)\le f(y) + M$ and $Hx_\mu=\mu$, so $\mu\in H(\mathcal{S}_y)$.
    
    Thus, $H(\mathcal{S}_y)=\{\mu:\phi(\mu)\le f(y) + M\}$, and substituting $\phi$ gives the quadratic form description.
\end{proof}

This result lets us translate between $\mu$-descriptions and $x$-descriptions. 
This is key for our analysis of coverage: even though the set is constructed in $\mu$-space, coverage has to be enforced in $x$-space, because it has to hold for every parameter in the original space.

This coverage analysis is trivial to do in the unconstrained case, but becomes challenging in the constrained space, hence the usefulness of this translation tool between the two spaces. 

Another consequence is the ability to describe confidence regions in $\mu$-space for SSB procedures, that come from mapping $1-\alpha$ confidence regions in $x$-space through $H$. 
While this is a conservative assumption that usually leads to overcoverage, the best-known regions in the constrained case prior to our work follow this procedure. 
For example, in this setting, one can build the $1-\alpha$ confidence region $\{ x \in \mathbb{R}^p:\|Kx-y\|^2_2 \leq Q_{\chi^2_n,1-\alpha}\}$ and map it through $H$.
\Cref{prop: mu_x_equiv} shows that in this case one obtains the region
\begin{equation}
    \label{eq:unconstrained_lambda1}
    \{\mu \in \mathbb{R}^k: \|\mu - B^\top y\|^2_{B^\top B} \leq Q_{\chi^2_n,1-\alpha} - \inf_{x'}\|Kx'-y\|^2_2 \},
\end{equation} 
which is a valid $1-\alpha$ ellipsoidal confidence region (coverage is inherited from the $x$-space coverage).
Unlike the Wald ellipsoid \eqref{eq:mu_space_unc}, the length of the different axes of the ellipse in \eqref{eq:unconstrained_lambda1} is dependent on the data $y$. 
Depending on the observation $y$, the ellipse \eqref{eq:unconstrained_lambda1} might be larger or smaller than the Wald ellipsoid.

As a summary, we have two different $1-\alpha$ confidence sets, each allowing for three equivalent descriptions:
\begin{align}
    \mathcal{R}^1 &= \{\mu \in \mathbb R^k :\|\mu - B^\top y\|^2_{B^\top B} \leq Q_{\chi^2_n,1-\alpha} - \inf_{x'}\|Kx'-y\|^2_2 \} \\
    &=\{\mu \in \RR^k : \inf_{Hx = \mu} \|Kx-y\|^2_2 \leq Q_{\chi^2_n,1-\alpha}\} \\
    &= \{Hx: \|Kx-y\|^2_2 \leq Q_{\chi^2_n,1-\alpha}\},
\end{align}
corresponding to $f(y) = 0$, $M = Q_{\chi^2_n,1-\alpha}$ in \Cref{prop: mu_x_equiv}, and 
\begin{align}
    \mathcal{R}^2 &=\, \{\mu \in \mathbb R^k :\|\mu - B^\top y\|^2_{B^\top B} \leq Q_{\chi^2_k,1-\alpha} \} \\
    &=\,\{\mu \in \RR^k : \inf_{Hx = \mu} \|Kx-y\|^2_2 \leq \inf_{x'}\|Kx'-y\|^2_2 + Q_{\chi^2_k,1-\alpha}\} \\
    &=\, \{Hx: \|Kx-y\|^2_2 \leq  \inf_{x'}\|Kx'-y\|^2_2 + Q_{\chi^2_k,1-\alpha}\}, 
\end{align}
corresponding to $f(y) = \inf_{x'}\|Kx'-y\|^2_2$ and $M = Q_{\chi^2_k,1-\alpha}$ in \Cref{prop: mu_x_equiv}. 

In order to later generalize in the constrained case, it is helpful to interpret both of these sets (and any of the form of \Cref{prop: mu_x_equiv}) as the inversion set of a hypothesis test:
\begin{equation}
    \label{eq:hypothesis_test_informal}
    H_0: Hx^* = \mu \quad \text{versus} \quad H_1: Hx^* \neq \mu.
\end{equation}
with test statistic $\lambda(\mu, y) = \inf_{Hx = \mu} \|Kx-y\|^2_2 -f(y)$, resulting in confidence regions of the form $\{\mu: \lambda(\mu, y) \leq M\}$.

Therefore, when selecting $f$ in \Cref{prop: mu_x_equiv}, we are equivalently choosing a test statistic, and then choosing the constant $M$ is the \emph{calibration} problem for that test statistic.

We then observe that $\mathcal{R}^1$ is the inversion region of the test statistic with one term $\lambda^1(\mu, y) = \inf_{Hx = \mu} \|Kx-y\|^2_2 $, and that $\mathcal{R}^2$ is the inversion region of the test statistic with two terms $\lambda^2(\mu, y) = \inf_{Hx = \mu} \|Kx-y\|^2_2 - \inf_{x'} \|Kx'-y\|^2_2 $, which corresponds to a scaled negative log-likelihood ratio (LLR) test-statistic.

We re-prove coverage of the regions obtained by taking $M = Q_{\chi^2_n,1-\alpha}$ for $\lambda^1$ and $M = Q_{\chi^2_k,1-\alpha}$ for $\lambda^2$ below
through the test inversion lens by analyzing the distributions of the test statistics under the null. 
This test-inversion argument for coverage will extend cleanly to the constrained setting, in which other direct techniques for coverage validity become more challenging. 
We do so with full generality over $K$ and $H$, we then recover the case in which $H \subseteq \row(K)$ and its rows are linearly independent.

\begin{proposition}
    [Exact null laws for both test statistics and orthogonal decomposition]
    \label{prop:both-lambdas}
    Assume $H\subseteq\row(K)$, and let $r:=\rank(H)$ and $R:=\rank(K)$. 
    Define $B:=K(K^\top K)^\dagger H^\top$ so that $B^\top K=H$, and set $c(y):=\inf_{x}\|Kx-y\|_2^2 =\|(I-P_{K})y\|_2^2$, where $P_{K}$ projects into the column space of $K$. 
    For $\mu\in\RR^k$, consider the statistics
    \[
        \lambda^1(\mu,y)=\inf_{Hx=\mu}\|Kx-y\|_2^2,
        \qquad
        \lambda^2(\mu,y)=\inf_{Hx=\mu}\|Kx-y\|_2^2-\inf_x\|Kx-y\|_2^2.
    \]
    Then, under $H_0:Hx^*=\mu$ with $y=Kx^*+\varepsilon$, $\varepsilon\sim\cN(0,I)$, the following hold:
    \begin{itemize}
    \item (Exact null distributions)
    Under $H_0$,
    \[
        c(y)\sim\chi^2_{\,n-R}, 
        \qquad 
        \lambda^1(\mu,y)\sim\chi^2_{\,n-R+r},
        \qquad
        \lambda^2(\mu,y)\sim\chi^2_{\,r}.
    \]
    Moreover, $c(y)$ and $\lambda^2(\mu,y)$ are independent. 
    If $H\subseteq \row(K)$ and $H$ has full row rank, then $r=\rank(H)=k \le R=\rank(K)\le n$. 
    
    \item (Inversion sets with exact calibration)
    For any $0<\alpha<1$, the regions
    \begin{align*}
        \mathcal{R}^1 &= \{\mu:\lambda^1(\mu,y)\le Q_{\chi^2_{n-R+r},1-\alpha}\}
        =\Big\{\mu:\|\mu-B^\top y\|_{B^\top B}^2\le Q_{\chi^2_{n-R+r},1-\alpha}-c(y)\Big\}\\
        \mathcal{R}^2 &= \{\mu:\lambda^2(\mu,y)\le Q_{\chi^2_r,1-\alpha}\}
        =\{\mu:\|\mu-B^\top y\|_{B^\top B}^2\le Q_{\chi^2_r,1-\alpha}\}
    \end{align*}
    have $1-\alpha$ frequentist coverage.
    \end{itemize}
\end{proposition}
\begin{proof}
    Recall, from the proof of \Cref{prop: mu_x_equiv}:
    \begin{equation}\label{eq:ls-identity}
    \inf_{Hx=\mu}\|Kx-y\|_2^2 = c(y)+ \|\mu-B^\top y\|_{B^TB}^2
    \end{equation}
    This immediately shows $\lambda^2(\mu,y)=\|\mu-B^\top y\|_{B^TB}^2$. 
    For the first claim, write $y=Kx^*+\varepsilon$ with $\varepsilon\sim\cN(0,I)$ under $H_0:Hx^*=\mu$. 
    Then $B^\top y=B^\top(Kx^*+\varepsilon)=Hx^*+B^\top\varepsilon\sim\cN(\mu,B^\top B)$, so
    \[
        \lambda^2(\mu,y)= \|B^\top y-\mu\|^2_{B^\top B}\sim\chi^2_{r}.
    \]
    On the other hand,
    \[
        c(y)=\|(I-P_K)y\|_2^2=\|(I-P_K)\varepsilon\|_2^2\sim\chi^2_{n-R}.
    \]
    Because $(I-P_K)\varepsilon$ is orthogonal to $P_K\varepsilon$ and $B^\top\varepsilon$ depends only on $P_K\varepsilon$ (as $B^\top$ maps into $\row(K)$), the random variables $c(y)$ and $\lambda^2(\mu,y)$ are independent. 
    Hence
    \[
        \lambda^1(\mu,y)=c(y)+\lambda^2(\mu,y)\sim\chi^2_{\,n-R}+\chi^2_{\,r}\stackrel{\textup{d}}{=}\chi^2_{\,n-R+r}.
    \]
    The second claim follows by test inversion using the exact null laws proved in the first claim.
\end{proof}
\begin{remark}
    \Cref{prop:both-lambdas} implies coverage of the SSB region calibrated by $\chi^2_n$ \eqref{eq:unconstrained_lambda1}, though the region can be shrunk when $r < R$ while maintaining coverage by using the smaller quantile of $\chi^2_{n-R+r}$.
\end{remark}

When $k = 1$, the projection through $H$ corresponds to taking the minimum and maximum of the quantity of interest over a particular data-dependent set. 
Writing $h^\top x$ for $Hx$ when $k = 1$, \eqref{eq:mu_space_unc} becomes the classic least squares interval $h^\top \underbrace{(K^\top K)^{-1}Ky}_{\hat{x}} \pm z_{1-\alpha/2}$. 
Noting that $z^2_{1-\alpha/2} = Q_{\chi^2_1,1-\alpha}$ and applying \Cref{prop: mu_x_equiv} shows that this is equivalent to the interval:
\begin{equation}
    \bigg[\min_{x: \|Kx-y\|^2_2 \leq \min_{x'}\|Kx'-y\|^2_2 + Q_{\chi^2_1,1-\alpha}} h^\top x, \; \max_{x: \|Kx-y\|^2_2 \leq \min_{x'}\|Kx'-y\|^2_2 + Q_{\chi^2_1,1-\alpha}} h^\top x \bigg].
\end{equation}
We will henceforth write intervals coming from minimizing and maximizing a quantity over the same constraint set compactly as:
\begin{equation}
    \begin{aligned}
         \min_{x}/\max_{x} \quad & h^\top x \\
         \st  \quad  & \Vert Kx-y \Vert_2^2 \leq \min_{x'} \|Kx'-y\|^2_2 + Q_{\chi^2_1,1-\alpha}.
    \end{aligned}
\end{equation}

\subsection{Virtues of constraints}

In the unconstrained setting, functionals outside the row space of $K$ are unobservable, leading to infinitely wide confidence intervals.
By contrast, incorporating physical constraints can make these same functionals identifiable by restricting the parameter space. 
This is one of the key benefits of constraints: for specific quantities of interest, they can transform an ill-posed problem into a well-posed one. 
In what follows, we begin with linear constraints $Ax^* \leq b$ and then note the minor changes needed for cone constraints $x^* \in \mathcal{C}$.

Analogously to the unconstrained case, define the oracle compatibility region as:
\begin{equation}
    \label{eq:Oracle_bounded}
    \mathcal{R}^O = \{\mu \in \mathbb{R}^k: \exists \, x \in \mathbb{R}^p \text{ s.t. } Ax \leq b, Kx = y, Hx = \mu\} = \{ Hx: x \in \mathbb{R}^p, Ax\leq b, Kx = y\}.
\end{equation}

We assume that \eqref{eq:Oracle_bounded} is non-empty in this section, which will happen if the true parameter $x^*$ satisfies the imposed constraints. We note that a set similar to \eqref{eq:Oracle_bounded} can also be constructed before observing the data by eliminating the constraint $Kx = y$, which then leads to worst-case methods used in the uncertainty quantification literature. 

We next show that the presence of noise does not affect the question of boundedness. 
Specifically, boundedness of $\mathcal{R}^O$ is equivalent to the existence of a procedure that gives a.s.-bounded $1-\alpha$ confidence regions, as formalized below.
\begin{lemma}
    [Noiseless versus noisy boundedness equivalence]
    \label{lemma:noise_less}
    For any $0 \leq \alpha \leq 1$, there exists a procedure to obtain an almost surely bounded confidence region for $Hx^*$ given the observation of $y = Kx + \varepsilon$, $\varepsilon \sim \mathcal{N}(0,I)$, and the constraint $Ax^*\leq b$ if and only if the oracle compatibility region \eqref{eq:Oracle_bounded} is bounded.
\end{lemma}
\begin{proof}
    We use throughout the proof that a non-empty closed convex set is bounded if and only if its recession cone is the origin. 
    We first note that the boundedness of $\mathcal{R}^O$ is independent of the observation $y$, since its recession cone is $\mathcal{R} = \{Hd: Ad \leq 0, Kd = 0\}$. 
    To prove that if $\mathcal{R}^O$ is not bounded, a $1-\alpha$ almost surely bounded confidence region cannot exist, consider a feasible $x_0$ and points of the form $x_t = x_0 + td$, where $t \in \mathbb{R}$ and $0 \neq Hd \in \mathcal R$. 
    Since $Kx_t = Kx_0$ and $Hx_t-Hx_0 = tHd$, there exist points that produce the same data distribution that are arbitrarily far away in functional space.
    This disallows finite confidence sets almost surely. 
    Formally, fix $\varepsilon\in(0,1)$; by a.s. finiteness there exists $M<\infty$ with $\mathbb P(D(\mathcal C(y))\le M)\ge 1-\varepsilon$, where $D$ is the diameter of a set. 
    Choose $t$ so large that $\|Hx_t-Hx_0\|>M$.
    On the event $\{D(\mathcal C(y))\le M\}$, the set $\mathcal C(Y)$ cannot contain both $Hx_0$ and $Hx_t$.
    Yet the data law under $x_t$ equals that under $x_0$, and uniform $1-\alpha$-coverage gives
    $\mathbb P(Hx_0\in\mathcal C(Y))\ge 1-\alpha$ and $\mathbb P(Hx_t\in\mathcal C(Y))\ge 1-\alpha$.
    Therefore
    \begin{equation}
        \mathbb P\big(\{Hx_0\in\mathcal C(Y)\}\cap\{Hx_t\in\mathcal C(Y)\}\cap\{D(\mathcal C(Y))\le M\}\big)
        \ \ge\ 1-2\alpha-\varepsilon.
    \end{equation}
    For $\alpha<1/2$ (or for any $\alpha<1/l$ if we pick $l$ points), this probability is positive, contradicting $\{D(Y)\le M\}$. 
    Hence no a.s.-bounded $1-\alpha$ confidence set exists when $\mathcal C^O$ is unbounded.
    
    In the other direction, for any $c \geq 0, \{\mu \in \mathbb{R}^k: \exists \, x \in \mathbb{R}^p \text{ s.t. } Ax \leq b, Kx = y, Hx = \mu\}$ is bounded if and only if $\{\mu \in \mathbb{R}^k: \exists \, x \in \mathbb{R}^p \text{ s.t. } Ax \leq b, \|Kx-y\|^2_2 \leq c, Hx = \mu\}$ is bounded. 
    This is because both sets share the same recession cone. 
    The quantity $c$ can then always be picked to obtain a $1-\alpha$ confidence set for any required $\alpha$ based on the distribution of the norm of the noise, so a bounded almost surely method exists when $\mathcal{R}^O$ is finite.
\end{proof}

\begin{theorem}
    [Boundedness characterization]
    \label{thm:finite_or_infinite_regions}
    The compatibility region \eqref{eq:Oracle_bounded} is bounded, and therefore a method to obtain an almost surely bounded confidence region exists, if and only if
    \begin{equation}
        \label{eq:bdd_chatacter}
        \forall \, d \in \mathbb{R}^p:\quad Kd=0,\ Ad\le 0 \ \Rightarrow\ Hd=0.    
    \end{equation}
\end{theorem}
\begin{proof}
    As we have used in the proof of \Cref{lemma:noise_less}, $\mathcal R^O$ is bounded if and only if its recession cone $\mathcal{R} = \{Hd: d\in \mathbb{R}^p, Kd = 0, Ad\leq 0\}$ equals the origin. 
    This happens if and only if \eqref{eq:bdd_chatacter} holds.
\end{proof}

We can furthermore derive a more refined test that can be performed to check whether upper/lower bounded regions can be constructed for each individual functional, by looking at the projection of the recession cone in each coordinate. 
This can be understood as a compatibility test of the observation model $K$, constraint set $\mathcal{X}$, and functionals of interest $H$, answering the question \say{Is there enough information in the constraint set $\mathcal{X}$ and the observation $K$ to be able to obtain bounded confidence regions for $H$?}

\begin{theorem}
    [Coordinate-wise boundedness via the recession cone]
    \label{thm:coordinate_bounded}
    Let $\mathcal{R} = \{Hd: d\in \mathbb{R}^p, Kd = 0, Ad\leq 0\}$ be the recession cone of $\mathcal{R}^O$ and let $\mathcal{D} := \{d \in \mathbb{R}^p: Ad \leq 0, Kd = 0\}$ be its preimage of $\mathcal{R}$ under $H$. 
    The projection of $\mathcal{R}^O$ into the $i$-th coordinate, $P_i$, is a (possibly unbounded) interval that falls into exactly one of:
    \[
    \begin{cases}
        \forall \, d \in \mathcal{D} \;\;  h_i^\top d=0 & \Longrightarrow\ P_i\ \text{is a finite interval},\\[4pt]
        \exists \, d \in \mathcal{D}:\ h_i^\top d<0,\ \nexists \, d\in \mathcal{D}:\ h_i^\top d>0
        & \Longrightarrow\ P_i=(-\infty,U_i],\\[4pt]
        \exists \, d\in \mathcal{D}:\ h_i^\top d>0,\ \nexists \, d\in \mathcal{D}:\ h_i^\top d<0
        & \Longrightarrow\ P_i=[L_i,\infty),\\[4pt]
        \exists \, d_+,d_-\in \mathcal{D}:\ h_i^\top d_+>0,\ h_i^\top d_-<0
        & \Longrightarrow\ P_i=(-\infty,\infty).
    \end{cases}
    \]
    The boundedness (though not the numerical values) of the endpoints $L_i$, $U_i$ depends only on $h_i$ and $\mathcal D$, not on $y$ or on the other rows $\{h_j\}_{j\ne i}$.
    Moreover, an upper/lower a.s.-bounded confidence set for $h_i^\top x$ exists if and only if $P_i$ is upper/lower bounded. 
\end{theorem}
Note that each functional can be checked independently, and the existence or lack thereof of bounded confidence sets for a given functional is independent of other observed functionals.
\begin{proof}
    Since the compatibility set $\mathcal{R}^O$ is non-empty, closed, and convex, its linear one-dimensional projection $P_i := \{\mu_i : \mu \in \mathcal{R}^O\}$ is a possibly unbounded interval of the form:
    \begin{align*}
        \bigg[\inf_{\mu \in \mathcal{R}^O} \mu_i, \sup_{\mu \in \mathcal{R}^O} \mu_i\bigg]
        &= \bigg[\inf_{x \in \mathbb{R^p}: Ax\leq b, Kx = y} (Hx)_i, \sup_{x \in \mathbb{R^p}: Ax\leq b, Kx = y} (Hx)_i\bigg] \\
        &= \bigg[\inf_{x \in \mathbb{R^p}: Ax\leq b, Kx = y} h_i^\top x, \sup_{x \in \mathbb{R^p}: Ax\leq b, Kx = y} h_i^\top x \bigg].
    \end{align*}
    It generally holds for any polyhedron $\mathcal{P}$ and vector $h$ that $\sup_{x \in \mathcal{P}} h^\top x < \infty \iff h^\top d \leq 0$ for all $d \in \text{rec}(\mathcal{P})$ and $\inf_{x \in \mathcal{P}} h^\top x > -\infty \iff h^\top d \geq 0$ for all $d \in \text{rec}(\mathcal{P})$. 
    Since $\text{rec}(\{ x \in \mathbb{R^p}: Ax\leq b, Kx = y\}) = \mathcal{D}$, this concludes the proof of the criterion. 
    One can then use the same argument in \Cref{lemma:noise_less} to prove the equivalence of this criterion to the ability of finding a.s.-bounded confidence sets for each individual functional $h_i$.
\end{proof}

\begin{remark}
    In the unconstrained case (in which $A$ is absent), $\mathcal{D}=\ker K$ is a subspace, so $P_i$ is either a finite interval or $(-\infty,\infty)$.
    Moreover, $P_i$ is finite if and only if $h_i\in\row(K)$.
\end{remark}

 \paragraph{Cone constraints.}
All of the results above extend almost verbatim to cone constraint $x^* \in \mathcal{C}$, where $\mathcal{C}$ is a closed, convex cone. 
Since $\text{rec}(\mathcal{C}) = \mathcal{C}$, the recession cone of the updated compatibility region
\begin{equation}
    \label{eq:Oracle_bounded_Cone}
    \mathcal{R}^O = \{\mu \in \mathbb{R}^k: \exists \, x \in \mathbb{R}^p \text{ s.t. } x \in \mathcal{C}, Kx = y, Hx = \mu\} = \{ Hx: x \in \mathbb{R}^p, x \in \mathcal{C}, Kx = y\}
\end{equation}
is $\text{rec}(\mathcal{R}^O) = \{Hd : d \in \mathcal{C}, Kd = 0\}$. 
Thus, \Cref{thm:finite_or_infinite_regions,thm:coordinate_bounded} follow directly by replacing $Ad \leq 0$ with $d \in \mathcal{C}$ throughout.

\section{Classical and modern perspectives on the non-negative single-functional case}
\label{sec:history}

In this section, we summarize key developments for the case of a single functional under non-negativity constraints.
This setting has received the most attention since at least 1964.
Following the notation of the previous section, when $k=1$ we identify $Hx$ with the scalar $h^\top x$.

\subsection{Classical strict bounds perspectives and the Burrus conjecture}

In his Ph.D.\ thesis \cite[p.~64]{Burrus1964}, Walter Burrus observed the $k = 1$ analogue of \Cref{prop: mu_x_equiv} in the unconstrained, full column rank case: the usual Gaussian (Wald) interval around the least squares estimator:
\begin{equation}
    \big[
        h^\top \hat{x} - z_{\alpha/2} \| h \|_{K^\top K},
        \,
        h^\top \hat{x} + z_{\alpha/2} \| h \|_{K^\top K}
    \big]
    \text{ with }
    \hat{x} = (K^\top K)^{-1}K^\top y,
\end{equation}
can be obtained as the solution of the pair of optimization problems:
\begin{equation}
    \begin{aligned}
         \min_{x}/\max_{x} \quad & h^\top x \\
         \st \quad  & \Vert Kx-y \Vert_2^2 \leq \min_{x'} \|Kx'-y\|^2_2 + z^2_{\alpha/2}.
    \end{aligned}
\end{equation}
On p.~67, he further claimed \say{if it is known a priori that $x^* \geq 0$, then the confidence level resulting from the solution will not be reduced if only $x \geq 0$ is considered}, and claims
\begin{equation}
    \label{eq:burrus1}
    \begin{aligned}
         \min_{x}/\max_{x} \quad & h^\top x \\
         \st \quad  & \Vert Kx-y \Vert_2^2 \leq \min_{x'} \|Kx'-y\|^2_2 + z^2_{\alpha/2} \\
         & x \geq 0,
    \end{aligned}
\end{equation}
has exact $1-\alpha$ coverage even when $K$ is not full column rank. 
He considered that the fact that the constraint set might be empty for some observations is \say{a serious conceptual problem} and proposes the more relaxed version:
\begin{equation}
    \label{eq:burrus2}
    \begin{aligned}
         \min_{x}/\max_{x} \quad & h^\top x \\
         \st \quad  & \Vert Kx-y \Vert_2^2 \leq \min_{x'\geq 0} \|Kx'-y\|^2_2 + z^2_{\alpha/2} \\
         & x \geq 0,
    \end{aligned}
\end{equation}
claiming that this interval will always be non-empty, but will have coverage at least $1-\alpha$. 
He believes that \say{in many cases, a weak interval that always exists is preferable to a stronger interval that may occasionally fail to exist}. 
It is now known that his claim that the confidence level will not be reduced whenever we add the extra constraint to the optimization problem is incorrect, and furthermore, his relaxed version may not even achieve $1-\alpha$: the claim that \eqref{eq:burrus2} has $1-\alpha$ coverage was later known as the Burrus conjecture, and was open until the disproof first appeared in \cite{BatlleDisproof} in 2023.
Note that since \eqref{eq:burrus2} is necessarily a larger interval than \eqref{eq:burrus1}, this also disproves coverage of \eqref{eq:burrus1}.

In the later book by Burrus with Rust \cite[p.~200]{Rust1972-cf}, they provide geometric evidence that \eqref{eq:burrus1} may not yield a $1-\alpha$ interval, while also claiming \eqref{eq:burrus2} has coverage. 
They correctly noted, however, that coverage can be guaranteed by first constructing any $1-\alpha$ confidence set $C(y)$ for $x^*$ and then optimizing:
\begin{equation}
    \label{eq:burrusgeneral}
    \begin{aligned}
         \min_{x}/\max_{x} \quad & h^\top x \\
         \st  \quad & x \in C(y) \\
         & x \geq 0.
    \end{aligned}
\end{equation}
This is the template later formalized and generalized as the strict bounds method \cite{stark1987velocity, kuusela2017strict}.
In their instance, $C(y)$ is constructed using the distribution of $(x-\hat{x})^\top K^\top K(x-\hat{x}) = \|K(x-\hat{x})\|^2_2$, where $\hat{x} \in \argmin_x \|Kx-y\|^2_2$ is any least squares estimator.
This can be seen as building a confidence set for $x^*$ using the log-likelihood ratio test. 
We point out here the important equality: 
\begin{equation}
    \label{eq:least_squares_eq}
    \|Kx-y\|^2_2 - \min_{x}\|Kx-y\|^2_2 = \|Kx-y\|^2_2 - \|K\hat{x}-y\|^2_2 = \|K(x-\hat{x})\|^2_2, 
\end{equation} 
a consequence of the normal equations $K^\top K\hat{x} = K^\top y$.
From this equality, a condition of the form $\|Kx-y\|^2_2 - \min_{x}\|Kx-y\|^2_2 \leq a$ can be turned into an equivalent condition of the form $\|K(x-\hat{x})\|^2_2 \leq a$. 
The book identifies the distribution of $\|K(x-\hat{x})\|^2_2$ as $\chi^2_r$ where $r = \rank(K)$, proposing therefore the two equivalent intervals: 
\begin{equation}
    \label{eq:burrus3a}
    \begin{aligned}
        \min_{x}/\max_{x} \quad & h^\top x\\
        \st  \quad & \|K(x-\hat{x})\|^2_2 \leq Q_{\chi^2_r,1-\alpha}\\ 
        & x \geq 0,
    \end{aligned}
\end{equation}
and
\begin{equation}
    \label{eq:burrus3b}
    \begin{aligned}
        \min_{x}/\max_{x} \quad & h^\top x\\
        \st  \quad & \|Kx-y\|^2_2 \leq \min_{x'}\|Kx'-y\|^2_2 + Q_{\chi^2_r,1-\alpha}\\ 
        & x \geq 0.
    \end{aligned}
\end{equation}

Rust and O'Leary \cite{OLeary1986} built upon this result to propose a practical algorithm based on \eqref{eq:burrus3b} in 1986. 
They identify the fact that the solutions of the optimization problems in \eqref{eq:burrus3b} can be found by defining the function $L(\phi) := \min_{h^\top x = \phi, x\geq 0} \|Kx-y\|^2_2$ and solving for the two roots of the equation $L(\phi) = \min_{x'}\|Kx-y\|^2_2 + Q_{\chi^2_r,1-\alpha}$, a result that can be seen as a corollary of $\Cref{prop: mu_x_equiv}$. 
They proved that $L$ is a convex, differentiable, and piecewise quadratic.

In 1994, Rust and O'Leary published an attempted proof of the Burrus conjecture \cite{rust1994confidence}, based on the dualization of \eqref{eq:burrus2}. 
The dual of any interval of the family:
\begin{equation}
    \label{eq:generic}
    \begin{aligned}
        \min_{x}/\max_{x} \quad & h^\top x\\
        \st  \quad & \|Kx-y\|^2_2 \leq g(y)\\ 
        & x \geq 0,
    \end{aligned}
\end{equation}
where $g(y)$ is a non-negative function, given by: 
\begin{align}
    \label{eq:duals}
    & \max_{u^\top K \leq h^\top}u^\top y - \sqrt{g(y)} \| u \|_2,  \\
    & \min_{u^\top K \geq h^\top}u^\top y + \sqrt{g(y)} \| u \|_2. \label{eq:duals2}
\end{align}
When the constraint $x\geq 0$ is absent, the inequalities of the constraints in \eqref{eq:duals} and \eqref{eq:duals2} become equality constraints. 
Note that this constraint corresponds to the bias of the estimator $u^\top y$, making it, in expectation, either equal or an upper/lower bound of $h$. 
We observe here that the constraint allows us to solve the optimization problem in the case in which $h$ is not in the row space of $K$, for which the set of $u$ satisfying $K^\top u = h$ would be empty. 

The proof attempt of the Burrus conjecture in that paper relies on studying this dual problem, under the claim that, since for every fixed vector $u$, $\frac{u^\top(Kx-y)}{\| u \|_2} \sim \mathcal N(0,1)$, this is also true when $u$ is the solution of the dual problems. 
However, \cite{tenorio2007confidence} correctly pointed out in 2007 that the dependence of $u$ on the data $y$ invalidates this ``pivotal'' argument (no pun intended) and thereby the proof of the conjecture.
In modern terms, the proof attempt can be seen as an early attempt to characterize the distribution of a given random function of the data $y$, which would later be formalized with test inversion theory, and prove that pivotal or almost-pivotal quantities such as the one they intended to define are not easy to find in this problem, making calibration quite challenging.

Tenorio et al. in \cite{tenorio2007confidence} also provided an (alleged and later disproved in \cite{BatlleDisproof}) counterexample to the conjecture and provided one-dimensional and two-dimensional reductions for any non-negativity constrained problem in which $K$ has full column rank. 
We defer the exposition of such reductions to \Cref{sec:reductions}.

In parallel to the analysis and developments of the particular methods tailored for the Gaussian linear model, the strict bounds method was developed in the seminal work of \cite{Stark1992} as a way to build confidence regions in arbitrary and even infinite-dimensional constrained settings.
The basic idea is to construct a $1-\alpha$ confidence region for $x^*$ in the original space, and then to minimize and maximize an arbitrary number of quantities of interest over this set.
When instantiated in our linear setting, this recovers the same structure as \eqref{eq:burrusgeneral}, and the confidence set is usually obtained by bounding the norm of the noise, so that the final interval reads:
\begin{equation}
    \label{eq:starkSSB}
    \begin{aligned}
        \min_{x}/\max_{x} \quad & h^\top x\\
        \st  \quad & \|Kx-y\|^2_2 \leq Q_{\chi^2_n,1-\alpha} \\ 
        & x \geq 0.
    \end{aligned}
\end{equation}
This interval was used in \cite{stanley_unfolding} in a particle unfolding application. 

As pointed out in \cite{Stark1992}, the coverage in \eqref{eq:starkSSB}, and more generally in any interval of the form \eqref{eq:burrusgeneral} whenever $C(y)$ is a $1-\alpha$ confidence set for $x^*$, holds for any number of functionals, in particular for $Hx^*$ for $k > 1$ by building the region defined as the product of the following $k$ intervals:
\begin{equation}
    \label{eq:starkSSBk}
        \prod_{i = 1}^k
        \quad
    \begin{aligned}
        \min_{x}/\max_{x} \quad & h_i^\top x \\
        \st  \quad & \|Kx-y\|^2_2 \leq Q_{\chi^2_n,1-\alpha} \\ 
        & x \geq 0.
    \end{aligned}
\end{equation}
This provides a first benchmark method for the constrained $k > 1$ setting that we will further explore in \Cref{sec:test_inversion}. 
Another alternative extension to the $k > 1$ is given by applying Bonferroni-adjusted one-at-a-time intervals: if $\mathbb{P}(h_i^\top x \in I_i) \geq 1-\alpha_i$ for all $i = 1,\dots,k$, and $\sum_{i = 1}^K \alpha_i \leq \alpha$, then $\mathbb{P}(Hx \in \prod_{i =1}^kI_i) \geq 1-\alpha$. 
However, note that this method is strictly worse than $\eqref{eq:starkSSBk}$ whenever the individual intervals are built using the strict bounds procedure $\eqref{eq:burrusgeneral}$, as optimizing a $1-\alpha$ confidence set over $x$ already guarantees joint coverage without the need for the Bonferroni correction.

Finally, \cite{Stark1992} explicitly suggested seeking confidence sets $C(y)$ tailored to the specific functional in the $k = 1$ case, rather than using the generic choice $\{x: \|Kx-y\|^2_2 \leq Q_{\chi^2_r,1-\alpha}\}$.
This idea anticipates the modern test-inversion viewpoint and motivates our developments in the sections that follow.

\subsection{Modern test-inversion perspectives and resolution of the Burrus conjecture}

Recently, \cite{BatlleDisproof} presented a theoretical framework based on test inversion that explains the coverage properties of many different optimization-based intervals, including all of those reviewed in this section. 
Prior to this work, rigorous coverage guarantees were typically available only when the optimization was performed over a preconstructed $1-\alpha$ confidence set for $x^*$ (i.e., intervals of the form $\eqref{eq:burrusgeneral}$, including $\eqref{eq:burrus3a}$, $\eqref{eq:burrus3b}$, and $\eqref{eq:starkSSB}$). 

Informally, for intervals obtained by solving
\begin{equation}
    \label{eq:test_inversion_informal}
    \begin{aligned}
        \min_{x}/\max_{x} \quad & h^\top x\\
        \st  \quad & \|Kx-y\|^2_2 \leq f(y) + M \\ 
        & x \in \mathcal{X},
    \end{aligned}
\end{equation}
the test inversion framework provides, for each choice of $f$, a quantity $\mathcal{Q}_{f, 1-\alpha}(x)$ (dependent on $K$, $h$, and $\mathcal{X}$) such that the interval has at least $1-\alpha$ coverage for a given $x$ (for the functional of interest $h^\top x$ under $y \sim P_{x}$ ) if and only if $M \ge \mathcal{Q}_{f, 1-\alpha}(x)$. Thus, uniform calibration for all $x \in \mathcal{X}$ is achieved by taking $M \geq \sup_{x \in \mathcal{X}}\mathcal{Q}_{f, 1-\alpha}(x)$.

Though computing this quantity can itself be computationally challenging.
We defer to \Cref{sec:test_inversion} the formal definition of $\mathcal{Q}_{f, 1-\alpha}(x)$ and the rest of the test inversion setup. We discuss below important consequences obtained from defining and computing this function $\mathcal{Q}_{f, 1-\alpha}(x)$:
\begin{itemize}
    \item 
    \emph{Tightening known valid methods.}
    Procedures that already have provable coverage (in particular, $\eqref{eq:burrus3a}$, $\eqref{eq:burrus3b}$, $\eqref{eq:starkSSB}$) correspond to specific choices of $f$ in \eqref{eq:test_inversion_informal}, but with $M > \sup_{x \in \cX} \cQ_{f,1-\alpha}(x)$.
    The test inversion machinery can then provide calibrated constants $M$ that retain the nominal coverage while being smaller than the original constants.
    This is because we obtain coverage in functional space, avoiding the conservativeness of first building a $1-\alpha$ set in parameter space, which in turn can provide smaller confidence regions in practice.
    \item 
    \emph{Diagnosing and fixing invalid methods.}
    Procedures that do not yet have provable coverage (in particular, \eqref{eq:burrus1}, \eqref{eq:burrus2}) also correspond to specific choices of $f$ in \eqref{eq:test_inversion_informal} but with $M$ not satisfying the condition $M \geq \sup_{x \in \mathcal{X}}\mathcal{Q}_{f, 1-\alpha}(x)$. Identifying points $x$ such that $\mathcal{Q}_{f, 1-\alpha}(x) > M$ shows that the procedure can undercover (as an instance, this was used to disprove the Burrus conjecture).
    However, the test inversion machinery can provide an alternative constant (in particular, any $M \ge \sup_{x \in \cX} \cQ_{f,1-\alpha}(x)$) for intervals of this form such that the uniform coverage is attained.
\end{itemize}

On the computational side, \cite{stanley2025confidenceintervalsfunctionalsconstrained} presented a way to calibrate the interval \eqref{eq:burrus2} using the theoretical setup in \cite{BatlleDisproof} by solving a quantile optimization problem restricted over a suitable confidence set.

In the next section, we present the test inversion setup formally, together with how it can be used to recover all the intervals presented in this section and how it generalizes to $k > 1$ functionals and thresholds $M$ that might depend on $x$ in \eqref{eq:test_inversion_informal}. 
In this work, we focus on the linear Gaussian case, but we note that the exposition in \cite{BatlleDisproof} treats a broader class of settings. 

\section{Test-inversion framework for the constrained multiple-functional case}
\label{sec:test_inversion}

In this section, we present the test inversion framework which will be used to calibrate confidence regions. 
We will only focus on the linear Gaussian setting. A more general description of test inversion for constrained single-functional inference can be found in \cite{BatlleDisproof}, which can be extended to the multiple-functional case as we do below.

\subsection{Summary of the single-functional case}
\label{subsec:k_equals_1}

This subsection follows the exposition from \cite{BatlleDisproof, stanley2025confidenceintervalsfunctionalsconstrained}.
Following the duality between hypothesis tests and confidence sets, we define, for every $\mu \in \mathbb R$, $\Phi_\mu := \{x \in \mathbb R^p: h^\top x =\mu\}$ and consider the family of hypothesis tests that test whether $h^\top x^* = \mu$ while preserving the known constraints of $x^* \in \mathcal{X}$:
\begin{equation}
    \label{eq:hypothesis_test}
    H_0: x^* \in \Phi_\mu \cap \mathcal{X} \quad \text{versus} \quad H_1: x^* \in \mathcal{X} \setminus \Phi_\mu.
\end{equation}
For every parameter $\mu$ and data $y$, we define a test statistic $\lambda(\mu, y)$ that rejects the test when the value is large. 
We use $\mathcal{X} = \{x \in \mathbb{R}^p : Ax\leq b\}$ throughout the section to connect with previous methods, but the same results hold for cone constraints by replacing $Ax \leq b$ for $x \in \mathcal{C}$ unless explicitly noted. 
All of the intervals previously discussed in the literature, including those in \Cref{sec:history}, arise from considering one of three different test statistics to test \eqref{eq:hypothesis_test}, and the analysis of each of them will become the main focus of this paper: 
\begin{align}
    \lambda^2_c(\mu, y) &= \min_{h^\top x = \mu, Ax\leq b} \|Kx-y\|^2_2 - \min_{Ax\leq b} \|Kx-y\|^2_2, \label{eq:lambda2c_k_equals_1} \\
    \lambda^2_u(\mu, y) &= \min_{h^\top x = \mu, Ax\leq b} \|Kx-y\|^2_2 - \min_{x \in \mathbb R^p} \|Kx-y\|^2_2, \label{eq:lambda2u_k_equals_1}\\
    \lambda^1(\mu, y) &= \min_{h^\top x = \mu, Ax\leq b} \|Kx-y\|^2_2. \label{eq:lambda1_k_equals_1}
\end{align}

We follow the naming convention $\lambda^1, \lambda^2_c, \lambda^2_u$ to indicate whether one or two terms are used, and when two terms, whether the second minimization is constrained ($c$) or unconstrained ($u$). 
We note that $\lambda^2_c$ corresponds to a re-scaled log likelihood ratio test statistic. 
It is understood in the above definitions that if there is no $x$ satisfying $h^\top x = \mu$ and the constraint $Ax \leq b$, we take $+\infty$ as the test statistic value so that the hypothesis is rejected even before seeing the data (alternatively, we can consider the test \eqref{eq:hypothesis_test} only for $\mu$ in the image of the constraint set under the quantity of interest  to perform test inversion).

For any general test statistic $\lambda(\mu, y)$, we build one of two types of confidence regions, which we define as sliced and global, by obtaining a suitable decision threshold function $d(\mu)$ or constant threshold $D$ and then defining:
\begin{align} 
    \mathcal{R}_{\text{s}} &= \{\mu \in \mathbb{R} : \lambda(\mu, y) \leq d(\mu) \}, \label{eq:int_def_s} \\ 
    \mathcal{R}_{\text{g}} &= \{\mu \in \mathbb{R} : \lambda(\mu, y) \leq D\}. \label{eq:int_def_g}
\end{align}
Under this setup, \cite[Lemma 2.2]{BatlleDisproof} characterizes the optimal (smallest) function $d(\mu)$ and constant $D$ such that $1-\alpha$ coverage is preserved for all $x$ satisfying $Ax\leq b$: they are obtained as optimum values of quantile optimization problems.

To define such problems, we first define, for each $x \in \mathcal{X}$, the random variable $Z_x$ corresponding to $\lambda(\mu, y)$ when $\mu$ is fixed to be $h^\top x$ and $y \sim P_x$, i.e., $y = Kx + \varepsilon$ with $\varepsilon \sim \mathcal N (0, I)$. More formally, letting $\mathcal{T}_x$ be the measurable map $\mathcal{T}_x: \varepsilon \to \lambda(h^\top x, Kx+\varepsilon)$, the random variable $Z_x$ has pushforward law $ (\mathcal{T}_x)_{*}(\mathcal{N}(0,I))$.

We note that the null hypothesis is composite, and while $x_1 \neq x_2$ with $h^\top x_1 = h^\top x_2$ will share the same null, we generally have the laws of $Z_{x_1}$ and $Z_{x_2}$ be different. 
We let $Q_{x,1-\alpha}$ be the quantile of $Z_x$ at level $1-\alpha$, so that $\mathbb{P}_{\varepsilon\sim\mathcal{N}(0,I)}\left(\lambda(h^\top x, Kx+\varepsilon) \leq Q_{x,1-\alpha} \right) = 1 - \alpha$. 
The next result characterizes the optimal thresholds as functions of these quantiles:

\begin{lemma}
    [Lemma 2.2 of \cite{BatlleDisproof}] 
    \label{lemma:coverage}
    The optimal (smallest) thresholds ensuring uniform $1-\alpha$ coverage over $\cX$ are:
    \begin{align} 
        d^*(\mu) &= \mu \mapsto \sup_{h^\top x = \mu, Ax \leq b} Q_{x,1-\alpha} \label{eq:max_q_mu} \\
        D^* &= \sup_{Ax \leq b} Q_{x,1-\alpha}. \label{eq:max_q}
    \end{align}
    Thus, the desired uniform coverage of \eqref{eq:int_def_s} and \eqref{eq:int_def_g} can be guaranteed at level $1-\alpha$ if and only if $d(\mu) \geq d^*(\mu)$ and $D \geq D^*$, respectively.
    Furthermore, coverage at a particular $x$ is achieved at level $1-\alpha$ if and only if $d(h^\top x)$ and $D$ is at least $Q_{x,1-\alpha}$.
\end{lemma}

We therefore identify the quantity $\mathcal{Q}_{f, 1-\alpha}(x)$ from the previous section as $Q_{x, 1-\alpha}$, the $1-\alpha$ quantile of $Z_x$. 
We postpone to \Cref{sec:test_statistics} the discussion of optimization problems and null statistic distributions appearing in \eqref{eq:max_q_mu} and \eqref{eq:max_q} for the three test statistics discussed, and assume for now that suitable $d(\mu)$ and $D$ ensuring the desired coverage can be computed. 

Finally, as a consequence of \cite[Theorem 2.4]{BatlleDisproof} (which uses the same argument used in \Cref{prop: mu_x_equiv}, but with the added constraint $x \in \mathcal{X}$), when $\lambda(\mu, y)$ is of the form $\min_{h^\top x = \mu, Ax\leq b} \|Kx-y\|^2_2 - f(y)$, we have the following equivalences: 
\begin{align}
    \mathcal{R}_{\text{s}} &= \{\mu \in \mathbb{R} : \min_{h^\top x = \mu, Ax\leq b} \|Kx-y\|^2_2 - f(y) \leq d(\mu) \} \label{eq:Rsmu} \\ &= \{h^\top x: x\in \mathbb{R}^p, Ax\leq b, \|Kx-y\|^2_2 \leq f(y) + d(h^\top x) \} \label{eq:Rsx} 
\end{align}
and
\begin{align}
    \mathcal{R}_{\text{g}} &= \{\mu \in \mathbb{R} : \min_{h^\top x = \mu, Ax\leq b} \|Kx-y\|^2_2 + f(y) \leq D\} \label{eq:Rgmu} \\  &= \{h^\top x: x\in \mathbb{R}^p, Ax\leq b, \|Kx-y\|^2_2 \leq f(y) + D \}, \label{eq:Rgx}
\end{align}
which correspond to the $\mu$-descriptions and $x-$descriptions of the $\mathcal{R}_s$ and $\mathcal R_g$ regions. 
Furthermore, since $\eqref{eq:Rgmu}$ is a convex set in the real line, it holds that:
\begin{align}
    \label{eq:Rg_x_opt}
        \mathcal{R}_g ~ = ~ 
    \begin{aligned}
        \min_{x}/\max_{x} \quad & h^\top x\\
        \st  \quad & \|Kx-y\|^2_2 \leq f(y) +  D\\ 
        & Ax \leq b.
    \end{aligned}
\end{align}
Since $\eqref{eq:Rsmu}$ is not necessarily convex unless $d$ is concave, $\mathcal{R}_s$ is not necessarily an interval, but it is nevertheless included in the following interval:
\begin{align}
    \label{eq:Rs_x_opt}
        \mathcal{R}_s ~ \subseteq ~
    \begin{aligned}
        \min_{x}/\max_{x} \quad & h^\top x\\
        \st  \quad & \|Kx-y\|^2_2 \leq f(y) + d(h^\top x)\\ 
        & Ax \leq b.
    \end{aligned}
\end{align}

Combining \eqref{eq:Rg_x_opt} with \Cref{lemma:coverage} is the key insight that lets us analyze coverage of a broad class of optimization-based methods and to tighten classical strict bounds methods that first construct a $1-\alpha$ confidence set for $x$ in parameter space. 
In \Cref{sec:test_statistics}, we connect the test statistics in \eqref{eq:lambda2c_k_equals_1}--\eqref{eq:lambda1_k_equals_1} to previously proposed optimization-based methods, together with the consequences of such analysis.

\subsection{Generalization to the multiple-functional case}

When the quantity of interest is the vector $Hx^* \in \RR^{k}$, the setup of \Cref{subsec:k_equals_1} extends directly. 
We test whether $Hx^* = \mu$ for each $\mu \in \mathbb{R}^k$ while preserving the known constraints. 
Defining $\Phi_\mu := \{x \in \mathbb R^p: H x =\mu\}$ analogously as above, the test remains:
\begin{equation}
    \label{eq:hypothesis_test_md}
    H_0: x^* \in \Phi_\mu \cap \mathcal{X} \quad \text{versus} \quad H_1: x^* \in \mathcal{X} \setminus \Phi_\mu.
\end{equation}
Given a test statistic $\lambda : \RR^{k} \times \RR^{n} \to \RR$ mapping a pair $(\mu, y)$ to the real line, the acceptance region is now the region in $\mathbb{R}^k$ defined by either of:
\begin{align} 
    \mathcal{R}_{\text{s}} &= \{\mu \in \mathbb{R}^k : \lambda(\mu, y) \leq d(\mu) \}, \label{eq:int_def_s_md} \\ 
    \mathcal{R}_{\text{g}} &= \{\mu \in \mathbb{R}^k : \lambda(\mu, y) \leq D\}, \label{eq:int_def_g_md}
\end{align}
for a slice-dependent threshold $d(\cdot)$ or a global constant $D$.
Let $\varepsilon \sim \cN(0,I)$ and define $Z_x := \lambda(Hx, Kx + \varepsilon)$ as before.
Write $Q_{x,1-\alpha}$ for the $(1-\alpha)$-quantile of $Z_x$.
Then, analogously to \Cref{lemma:coverage}, the optimal (smallest) thresholds ensuring uniform $1-\alpha$ coverage over $\cX$ are given by:
\begin{align}
    \label{eq:max_q_mu_md}
    d^*(\mu) &= \mu \mapsto \sup_{Hx = \mu, Ax \leq b} Q_{x,1-\alpha}, \\
    D^* &= \sup_{Ax \leq b} Q_{x,1-\alpha}. \label{eq:max_q_md}
\end{align}

The three single-functional test statistics \eqref{eq:lambda2c_k_equals_1}, \eqref{eq:lambda2u_k_equals_1}, and \eqref{eq:lambda1_k_equals_1} easily generalize to this case: 
\begin{align}
    \lambda^2_c(\mu, y) &= \min_{Hx = \mu, Ax\leq b} \|Kx-y\|^2_2 - \min_{Ax \leq b} \|Kx-y\|^2_2, \label{eq:lambda2c_k_bigger_1_md} \\
    \lambda^2_u(\mu, y) &= \min_{Hx = \mu, Ax\leq b} \|Kx-y\|^2_2 - \min_{x \in \mathbb R^p} \|Kx-y\|^2_2, \label{eq:lambda2u_k_bigger_1_md}\\
    \lambda^1(\mu, y) &= \min_{Hx = \mu, Ax\leq b} \|Kx-y\|^2_2. \label{eq:lambda1_k_bigger_1_md}
\end{align}

When $k > 1$ and for test statistics of the family $\lambda(\mu, y) = \min_{Hx = \mu, Ax\leq b} \|Kx-y\|^2_2 - f(y)$, the global acceptance region $\eqref{eq:int_def_g_md}$ is a convex set in $\mathbb{R}^k$ that is not necessarily a product of intervals, and the sliced acceptance region \eqref{eq:int_def_s_md} is a general region of $\mathbb{R}^k$. 
Similarly to the $k = 1$ case, we can express these regions as linear images under the quantity of interest matrix $H$ of certain regions in the original parameter space ($x$-description), which will be helpful to connect with optimization-based methods. 
\begin{lemma}[Equivalence between $x$-description and $\mu$-description in multiple dimensions]
    \label{lemma:x_mu_md}
    For any $f : \RR^{n} \to \RR$ and $d: \RR^K \to \RR$, we have
    \begin{align}
        \mathcal{R}_{\text{s}} &= \{\mu \in \mathbb{R}^k : \min_{H x = \mu, Ax\leq b} \|Kx-y\|^2_2 -f(y) \leq d(\mu) \} \label{eq:Rsmu_md} \\ &= \{Hx: x\in \mathbb{R}^p, Ax\leq b, \|Kx-y\|^2_2 \leq f(y) + d(Hx) \}, \label{eq:Rsx_md} 
    \end{align}
    and
    \begin{align}
        \mathcal{R}_{\text{g}} &= \{\mu \in \mathbb{R}^k : \min_{H x = \mu, Ax\leq b} \|Kx-y\|^2_2 - f(y) \leq D\} \label{eq:Rgmu_md} \\  &= \{H x: x\in \mathbb{R}^p, Ax\leq b, \|Kx-y\|^2_2 \leq f(y) + D \}. \label{eq:Rgx_md}
    \end{align}
\end{lemma}
\begin{proof}
    The result follows directly from the proof of \Cref{prop: mu_x_equiv}, by replacing the single-row $h^\top$ with the multi-row matrix $H$. 
    The same least squares decomposition and orthogonality argument apply without modification, although, because of the extra constraint $Ax \leq b$, the optimization problem $\min_{H x = \mu, Ax\leq b} \|Kx-y\|^2_2$ cannot be explicitly solved.
\end{proof}

Observe from the statement of the previous lemma that both $\mathcal{R}_g$ and $\mathcal{R}_s$ admit a $\mu-$description and an $x-$description. 
Since $\mathcal{R}_g$ is not necessarily a product of intervals, there is no direct equality such as $\eqref{eq:Rg_x_opt}$, which we had in the $k = 1$ case. 
Nevertheless, we can consider, for both the sliced and the global regions, their following $x$-description bounding boxes in $\mathbb{R}^k$:
\begin{align}
    \label{eq:prod_opt_rs_md}
    \mathcal{R}_{\text{s}}  \subseteq \prod_{i = 1}^k
    \quad
    \begin{aligned}
     \min_{x}/\max_{x} \quad & h_i^\top x\\
      \st  \quad & \|Kx-y\|^2_2 \leq f(y) + d(h^\top x)\\ 
      & Ax \leq b,
    \end{aligned}
\end{align}
and 
\begin{align}
    \label{eq:prod_opt_rg_md}
    \mathcal{R}_{\text{g}} \subseteq  \prod_{i = 1}^k
    \quad 
    \begin{aligned}
     \min_{x}/\max_{x} \quad & h_i^\top x\\
      \st  \quad & \|Kx-y\|^2_2 \leq f(y) + D\\ 
      & Ax \leq b.
    \end{aligned}
\end{align}
We call this the $x-$description bounding box. 
In particular, this shows that the previously discussed simultaneous method for multiple functionals \eqref{eq:starkSSBk} is a bounding box encapsulating a convex region that, while not being a product intervals, can already provide $1-\alpha$ coverage while being a set of smaller volume.

\begin{remark}
    If $\rank(H) = p$, then testing $Hx^* =\mu$ is equivalent to testing $x^* = x$, because there is enough information in the functionals to recover the full true parameter. 
    Therefore, inverting the test yields a $1-\alpha$ confidence test for $x^*$, recovering strict bounds-type methods. 
    Hence, the gains of the methods presented here (with respect to strict bounds-type methods) are largest when $\rank(H) \ll p$.
\end{remark}

\section{Unifying different test statistics and calibration strategies}
\label{sec:test_statistics}

As mentioned in the previous section, most optimization-based approaches to constrained confidence regions can be expressed as test inversions of one of three test statistics: 
\begin{align}
    \lambda^2_c(\mu, y) &= \min_{Hx = \mu, Ax\leq b} \|Kx-y\|^2_2 - \min_{Ax \leq b} \|Kx-y\|^2_2, \label{eq:lambda2c_k_bigger_1_md_bis} \\
    \lambda^2_u(\mu, y) &= \min_{Hx = \mu, Ax\leq b} \|Kx-y\|^2_2 - \min_{x \in \mathbb R^p} \|Kx-y\|^2_2, \label{eq:lambda2u_k_bigger_1_md_bis}\\
    \lambda^1(\mu, y) &= \min_{Hx = \mu, Ax\leq b} \|Kx-y\|^2_2. \label{eq:lambda1_k_bigger_1_md_bis}
\end{align}

In this section, we make the mapping to prior methods explicit and then analyze each statistic from a calibration point of view: how to choose thresholds to guarantee $1 - \alpha$ coverage, when optimal constants exist, and where the quantile maximizers lie, highlighting convexity and extremal structure that make some cases tractable.

\subsection{Mapping methods to test statistics}

\begin{table}[!t]
    \centering
    \caption{Acceptance regions for the three test statistics \eqref{eq:lambda2c_k_bigger_1_md}--\eqref{eq:lambda1_k_bigger_1_md}.Here, $\delta$ denotes either the global constant $D$ or the sliced threshold $d(\mu)$ (in the $\mu$-description) or $d(Hx)$ (in the $x$-description).}
    \setlength{\tabcolsep}{2pt}
    \small
    \begin{tabularx}{\textwidth}{l l l} 
    \toprule
    Statistic
    & \makecell[c]{%
       \(\mu\)-description 
       } 
    & \makecell[c]{%
       \(x\)-description (bounding box)
       } 
       \\
    \midrule
    \textbf{\(\lambda_c^2\)} &
    \(\{\mu:\min_{\substack{Hx=\mu,\\ Ax\le b}}\|Kx-y\|_2^2
    -\min_{\substack{Ax\le b}}\|Kx-y\|_2^2
    \le \delta^2_c\}\) &$
    \displaystyle \prod_{i = 1}^k
    ~ \begin{aligned}
    & \min_{x}/\max_{x} \quad h_i^\top x\\
    & \text{s.t. } \|Kx-y\|^2_2 \leq \min_{Ax' \leq b}\|Kx'-y\|_2^2 + \delta^2_c\\ 
    & \qquad Ax \leq b.
    \end{aligned}$\\
    \arrayrulecolor{black!25}\midrule
    \textbf{\(\lambda_u^2\)} &
    \(\{\mu:\min_{\substack{Hx=\mu,\\ Ax\le b}}\|Kx-y\|_2^2
    -\min_{x \in \mathbb R^p}\|Kx-y\|_2^2
    \le \delta^2_u\}\) &$
    \displaystyle \prod_{i = 1}^k
    ~ \begin{aligned}
      &\min_{x}/\max_{x}\quad h_i^\top x\\
     & \text{s.t. } \|Kx-y\|^2_2 \leq \min_{x '\in \mathbb R^p}\|Kx'-y\|_2^2 + \delta^2_u\\ 
      & \qquad Ax \leq b .
    \end{aligned}$\\
    \arrayrulecolor{black!25}\midrule
    \textbf{\(\lambda^1\)} &
    \(\{\mu:\min_{\substack{Hx=\mu,\\ Ax\le b}}\|Kx-y\|_2^2\le \delta^1\}\) &
    $
    \displaystyle \prod_{i = 1}^k 
    ~ \begin{aligned}
     & \min_{x}/\max_{x} \quad h_i^\top x\\
      & \text{s.t. }  \|Kx-y\|^2_2 \leq \delta^1\\ 
       & \qquad Ax \leq b.
    \end{aligned}$\\
    \arrayrulecolor{black}\bottomrule
    \end{tabularx}
    \label{table:rosetta_stone}
\end{table}

For each statistic in \eqref{eq:lambda2c_k_bigger_1_md_bis}--\eqref{eq:lambda1_k_bigger_1_md_bis}, there are four natural acceptance sets: the sliced set \eqref{eq:Rsmu_md} or the global set \eqref{eq:Rgmu_md}, each optionally embedded into a product of intervals via \eqref{eq:prod_opt_rs_md} or \eqref{eq:prod_opt_rg_md}. 
Note they all live in the same functional space $\mathbb{R}^k$, and by construction, we always have that the $\mu$-description region is a subset of the $x$-description bounding box region. 
\Cref{table:rosetta_stone} shows the $\mu$-descriptions and the $x$-description bounding boxes of the acceptance regions of the three test statistics. 
We emphasize through the notation $\delta^2_u, \delta^2_c, \delta^1$, that valid thresholds achieving the nominal coverage are not the same for different test statistics. 
We formulate the corresponding optimization problems that yield the optimum $\delta$ in each case in the next subsection. 
Since there is a pointwise inequality between the three test statistics of the form $\lambda^2_c(\mu, y) \leq \lambda^2_u(\mu, y) \leq \lambda^1(\mu, y)$, it immediately follows that any threshold that ensures valid coverage for $\lambda^1$ also ensures valid coverage for $\lambda^2_u$ and $\lambda^2_c$; similarly, a valid threshold for $\lambda^2_u$ is valid for $\lambda^2_c$.
 
However, if we have a valid threshold for one test statistic, using it as a threshold for a smaller test statistic will provide a larger region than using the threshold for the original test statistic, so it is generally suboptimal in practice. 
Nonetheless, the inequality of valid thresholds also implies that the optimum thresholds obey $(\delta^{2}_c)^* \leq (\delta^{2}_u)^* \leq (\delta^{1})^*$, for both global constants $D$ and, pointwise, for sliced thresholds $d(\mu)$.

\begin{table}[!t]
    \centering
    \caption{Reinterpretation of historical methods for the non-negatively constrained case via the unified test inversion framework. All methods use the global version of test inversion, so $\delta = D$ is a constant.}
    \begin{tabular}{l c c c c c}
    \toprule
    \textbf{Method} & \textbf{Equation} & \textbf{Statistic} & $\delta$ & \textbf{ $1-\alpha$}?& \textbf{$\delta$ optimal?}  \\
    \midrule
    Burrus unconstrained (1964)                & \eqref{eq:burrus1} & $\lambda^2_u$  & $Q(\chi^2_1)$       & No  & N/A \\
    Burrus constrained (1964)                & \eqref{eq:burrus2} & $\lambda^2_c$  & $Q(\chi^2_1)$       & No  & N/A \\
    Rust and Burrus (1972)        & \eqref{eq:burrus3a}, \eqref{eq:burrus3b} & $\lambda^2_u$      & $Q(\chi^2_r)$                  & Yes & No \\
    Strict bounds $(k = 1)$               & \eqref{eq:starkSSB} & $\lambda^1$      & $Q(\chi^2_n)$            & Yes & No \\
    Strict bounds $(k > 1)$    & \eqref{eq:starkSSBk} & $\lambda^1$   & $Q(\chi^2_n)$       & Yes  & No \\
    \bottomrule
    \end{tabular}
    \label{table:old_methods}
\end{table}

Furthermore, using \Cref{table:rosetta_stone}, the intervals in \Cref{sec:history} can be interpreted as test inversions with specific choices of statistic and threshold.
\Cref{table:old_methods} summarizes the mapping, whether the stated $\delta$ attains at least $1 - \alpha$ coverage, and whether it is optimal.
Two takeaways are: (i) for methods with coverage, the constants can be improved; (ii) the Burrus intervals require calibration.
In the next subsection, we address calibration for each statistic.
For the non-negative case, we improve on the methods in \Cref{table:old_methods}: we provide theory and an algorithm to obtain the optimal constant $\delta^*$ for $\lambda^2_u$-based regions such as \eqref{eq:burrus3a} or \eqref{eq:burrus3b} (\Cref{cor:improved_constants_2u}) and for $\lambda^1$-based regions such as \eqref{eq:starkSSB} or \eqref{eq:starkSSBk}  (\Cref{cor:improved_constants_1}), including the extensions for $k \geq 1$. 
By contrast, an optimal calibration algorithm of $\lambda^2_c$-based regions such as $\eqref{eq:burrus2}$ remains open in general.

\subsection{The calibration problem}
\label{sec:calibration_problem}

This section compiles the known results on choosing thresholds $\delta$ so that the regions in \Cref{table:rosetta_stone} achieve $1-\alpha$ coverage. 
Test inversion yields optimal constants for the $\mu$-descriptions, which automatically provide coverage for their relaxed $x$-description bounding boxes. While this already significantly improves prior methods, quantifying how much the bounding box constant can be reduced below the $\mu$-description optimum while preserving coverage remains an open question. In what follows, we focus on optimal constants for the $\mu$-descriptions.

Recall that for any test statistic $\lambda$,
\begin{align}
    \label{eq:max_q_mu_md_2}
    d^*(\mu) &= \mu \mapsto \sup_{Hx = \mu, Ax \leq b} Q_{x,1-\alpha} \\
    D^* &= \sup_{Ax \leq b} Q_{x,1-\alpha} \label{eq:max_q_md_2} \\
    \mathbb{P}_{\varepsilon\sim\mathcal{N}(0,I)}&\left(\lambda(Hx, Kx+\varepsilon) \leq Q_{x,1-\alpha}\right) = 1-\alpha \label{eq:quantile_def},
\end{align}
and we define $Z_x$ to be the random variable $\mathcal{T}_x(\varepsilon) = \lambda(Hx, Kx+\varepsilon)$ for $\varepsilon \sim \mathcal N(0,I)$.
We therefore will start the analysis of each test statistic by describing $\mathcal{T}_x(\varepsilon)$ and analyzing properties of the quantile function $Q_x$ as a function of $x \in \mathcal{X}$, to understand the maximization problems in \eqref{eq:max_q_mu_md_2} and \eqref{eq:max_q_md_2}. 
All of the previously proposed methods we analyze and improve upon here correspond to the global case \eqref{eq:max_q_md_2}; in general, computing $\eqref{eq:max_q_mu_md_2}$ is typically challenging, since it requires the whole function $d^*:\mathbb{R}^k \to \mathbb{R}$, so our results primarily target the global case.

Beside properties of the quantile function $Q_x$, stochastic upper bounds for the test statistics $\lambda(Hx, Kx + \varepsilon)$ are also of interest. 
This is because if a random variable $X$ stochastically dominates $\lambda(Hx, Kx+\varepsilon)$ for all $x \in \mathcal{X}$, then by stochastic dominance $Q_{x,1-\alpha} \leq Q_{X,1-\alpha}$ for all $\alpha$ and for all $x$. 
Therefore, $Q_{X,1-\alpha} \geq D^*$ and it can be used to obtain $1-\alpha$ coverage (possibly conservative).
See \cite[Section 3.2]{BatlleDisproof} for an in-depth treatment of stochastic dominance in this context. From the pointwise inequalities, there is a chain of stochastic dominance relations between our three test statistics and $\chi^2$ distribution of varying degrees.
\begin{proposition}
    [Stochastic dominance chain between test statistics]  
    Let $\lambda^2_c$, $\lambda^2_u$, $\lambda^1$ be as in \eqref{eq:lambda1_k_bigger_1_md_bis}--\eqref{eq:lambda2c_k_bigger_1_md_bis}. 
    For each $x \in \mathcal{X}$, let $Z^{2}_{c,x}, Z^{2}_{u,x}, Z^{1}_{x}$ denote $\lambda(Hx, Kx+\varepsilon)$ for $\lambda \in \{\lambda^2_c, \lambda^2_u, \lambda^1\}$ with $\varepsilon \sim \mathcal N(0,I)$. 
    Then, we have
    \begin{equation}
        \label{eq:chain_stochastic_dominance}
        Z^2_{c,x} \preceq Z^2_{u,x} \preceq Z^1_x \preceq\chi^2_n,
    \end{equation}
    and 
    \begin{equation}
        \label{eq:mini_chain_stochastic_dominance}
        Z^2_{u,x} \preceq \chi^2_{\rank(K)} \preceq\chi^2_n.
    \end{equation}
\end{proposition}
\begin{proof}
    It follows immediately from the pointwise inequalities:
    \begin{align}
         \lambda_c^2(Hx, Kx + \varepsilon) &= \min_{H\xi = Hx,\, A\xi \leq b} \|K\xi - Kx - \varepsilon\|_2^2
        - \min_{A\xi \leq b} \|K\xi - Kx - \varepsilon\|_2^2 
        \notag \\
        &\leq \min_{H\xi = Hx,\, A\xi \leq b} \|K\xi - Kx - \varepsilon\|_2^2
        - \min_{\xi} \|K\xi - Kx - \varepsilon\|_2^2 
        &= \lambda_u^2(Hx, Kx + \varepsilon) \notag \\
        &\leq \min_{H\xi = Hx,\, A\xi \leq b} \|K\xi - Kx - \varepsilon\|_2^2 
        &= \lambda^1(Hx, Kx + \varepsilon) \notag \\
        &\leq \|\varepsilon\|_2^2, \notag
    \end{align}
    and
    \begin{align*}
         \lambda_u^2(Hx, Kx + \varepsilon) &= \min_{H\xi = Hx,\, A\xi \leq b} \|K\xi - Kx - \varepsilon\|_2^2
        - \min_{\xi} \|K\xi - Kx - \varepsilon\|_2^2 
        \notag \\
        &\leq \|\varepsilon\|_2^2 -\min_{\xi} \|K\xi - Kx - \varepsilon\|_2^2 \\     &= \|\varepsilon\|_2^2 -\min_{\xi'} \|K\xi' - \varepsilon\|_2^2  \\
        & \leq  \|\varepsilon\|_2^2.
    \end{align*}
    Because the inequalities hold almost surely, this implies stochastic dominance when $\varepsilon \sim \mathcal{N}(0,I)$, in which case $\|\varepsilon\|_2^2 \sim \chi^2_n$ and $\|\varepsilon\|_2^2 -\min_{\xi'} \|K\xi' - \varepsilon\|_2^2 \sim \chi^2_{\rank(K)}$.
    This concludes the proof.
\end{proof}

\subsubsection{Test statistic with constrained second term \texorpdfstring{$\lambda^2_c$}{lambda1}}

We first analyze $\lambda^2_c(\mu, y) = \min_{Hx = \mu, Ax\leq b} \|Kx-y\|^2_2 - \min_{Ax \leq b} \|Kx-y\|^2_2$, which is a natural choice as it corresponds to the log-likelihood ratio statistic. 
However, among the three statistics, this has the least structure in terms of the quantile function and developing a general theory about where the maximizer of the quantile optimization resides remains open. 

We begin describing the function $\mathcal{T}^2_{c,x}(\varepsilon) := \lambda^2_c(Hx, Kx+\varepsilon)$. 
Relabeling the optimization variable to $\xi$ and changing coordinates, we have
\begin{equation*}
    \mathcal{T}^2_{c,x}(\varepsilon) = \min_{H\xi = Hx, A\xi\leq b} \|K\xi- Kx -\varepsilon\|^2_2 - \min_{A\xi \leq b} \|K\xi-Kx -\varepsilon\|^2_2.
\end{equation*}
Now, we define $z := \xi - x$, so that 
\begin{equation*}
    \mathcal{T}^2_{c,x}(\varepsilon)= \min_{Hz = 0, Az\leq b-Ax} \|Kz -\varepsilon\|^2_2 -\min_{ Az\leq b-Ax} \|Kz -\varepsilon\|^2_2.
\end{equation*}
Let $Q^2_{c}(x)$ be the $1-\alpha$ quantile of $\mathcal{T}^2_{c,x}(\varepsilon)$ for $\varepsilon \sim \mathcal{N}(0, I)$.
We fix $\alpha$ throughout so we omit its dependence to simplify the notation.
We summarize below the known properties and open questions regarding this function:

\begin{itemize}[leftmargin=0em]
    \item 
    \emph{Convexity.}
    $Q^2_c(x)$ is generally neither a convex or a concave function.
    As a simple counterexample, consider $H = (1,1)$, $y = x^*+ \varepsilon \in \mathbb R^2$, $x^*\geq 0$, $\varepsilon \sim \mathcal N(0, I)$. 
    The quantile function can be shown to be neither convex nor concave; see, e.g., \cite[Figure 1.1]{stanley2025confidenceintervalsfunctionalsconstrained} and the slice $x_2 =0$.

    \item
    \emph{Maximization.}
    No general result is known for the maximizer of $Q^2_c(x)$ over $Ax\leq b$ (or over $Ax \leq b$, $Hx = \mu$ for sliced versions) for a given choice of problem parameters $(K, H, A, b)$.
    Furthermore, the quantile maximizer can (i) escape to $\infty$ for unbounded regions, and (ii) depend on the quantile level $\alpha$. 
    (These comments also apply to cone constraints.)
    The same two-dimensional setup as above exemplifies this behavior, as it can be shown that the maximizer over $x\geq 0$ of the quantile is at $x = 0$ for $1-\alpha$ sufficiently close to 1 and the quantile function grows towards $x = \infty$ when $1-\alpha$ is sufficiently close to $0$ (note that the optimal value is still finite in both cases). 
    This is a counterexample for both the linear constraints and the cone constraints case.
    
    \item
    \emph{Bounds.}
    If $K$ has full column rank and there are no constraints, \Cref{prop:both-lambdas} shows $Z^2_{c, x} \stackrel{\textup{d}}{=} \chi^2_{\rank(H)}$. 
    Under constraints, it is however not true in general that $Q^2_{c}(x) \leq Q_{\chi^2_{\rank(H)}}$ for all $\alpha$ (or, in terms of stochastic dominance, that $Z^2_{c,x} \preceq \chi^2_{\rank(H)}$). 
    The $k=1$ version of this statement is the Burrus conjecture, which is recently shown to be false \cite{BatlleDisproof}. 
    From $\eqref{eq:chain_stochastic_dominance}$ and $\eqref{eq:mini_chain_stochastic_dominance}$, we have $Q^2_{c}(x) \leq Q_{\chi^2_{\rank(K)}}$.
    One might ask if this bound is tight.
    While the bound is not necessarily tight for a given $\alpha$ and problem instance (specific choice of $K$, $H$, $\mathcal{X}$), we show below that if we aim for a bound that holds for all $\alpha$ and problem instances, the stochastic bound $Z^2_{c,x} \preceq \chi^2_{\rank(K)}$ is tight in the family of gamma distributions $\{\Gamma(t, 2)\}_{t \in \mathbb{R}^+}$ (which includes $\chi^2_k$, as $\Gamma(k/2, 2)$), a family over which stochastic order is a total order.
    \begin{lemma}
        \label{lemma:javier}
        For every $\epsilon >0$, there exists $K$, $H$, and $x^* \geq 0$ such that the non-negatively constrained problem induces a random variable $Z^2_{c,x^*} \npreceq \Gamma(\rank(K)/2 -\epsilon, 2)$. 
        Therefore, $\rank(K)$ is the smallest $t > 0$ such that $Z^2_{c, x} \preceq \Gamma(t/2, 2)$ holds uniformly over $x$, constraint, and $K$, $H$. 
    \end{lemma}
    \begin{proof}
        We will prove that $\Gamma(\rank(K)/2 -\epsilon, 2)$ does not stochastically dominate $Z^2_{c,x^*}$ by comparing their moment generating functions $\mathbb{E}[e^{tX}]$ at a particular $t > 0$ in a given family of non-negatively constrained examples of increasing dimension. 
        Let $K = I_n$ ($\rank(K) = n$), $H = (1,\dots,1,-1)$, and $x^* = (0,\dots,0,M)$. 
        Consider the following region of $\mathbb{R}^n$: 
        \begin{equation*}
            \mathcal{A} = \{ \varepsilon \in \mathbb{R}^n : \varepsilon_i \geq 0, \ \varepsilon_i + \varepsilon_n \leq 0 \ (i = 1, \dots, n-1) \}.
        \end{equation*}
        On $\mathcal A$, feasibility is controlled entirely by the last coordinate: if $M$ is sufficiently large, the constraint induced by $H$ is inactive except on $\mathcal A$. 
        Thus, for every fixed $\varepsilon\in\mathcal A$, we eventually have $\mathcal{T}^2_{c, x^*} (\varepsilon)=\|\varepsilon\|_2^2$ as $M\to\infty$. 
        More precisely, for each $\delta>0$ there exists $M$ large enough such that there is a subset $\mathcal A_M\subseteq \mathcal A$ with $\mathbb P(\mathcal A_M)\ge (1-\delta)\mathbb P(\mathcal A)$ on which $\mathcal{T}^2_{c, x^*}=\|\varepsilon\|_2^2$. 
        Together with the integrand being non-negative, it follows that:
        \begin{equation*}
            \mathbb{E}_{\varepsilon\sim \mathcal{N}(0,I_n)}\!\big[e^{tZ^2_{c,x^*}}\big]
            = \int_{\mathbb{R}^n} e^{t\,\mathcal{T}^2_{c,x^*}(\varepsilon)} f(\varepsilon)
            \ge \int_{\mathcal{A}_M} e^{t \|\varepsilon\|_2^2} f(\varepsilon)
            \ge (1-\delta)\!\int_{\mathcal{A}} e^{ t \|\varepsilon\|_2^2} f(\varepsilon),
        \end{equation*}
        where $t \in (0,\tfrac{1}{2})$, $f$ is the density of $\mathcal N(0,I_n)$, $\mathcal A_M \subseteq \mathcal A$ with $1_{\mathcal A_M}\uparrow 1_{\mathcal A}$ as $M\to\infty$, and for any fixed $\delta>0$, there exists $M$ large enough so that the last inequality holds. 
        Then, we have
        \begin{align*}
              \int_{\mathcal{A}}e^{ t \|\varepsilon\|_2^2} f(\varepsilon) 
              &= \int_{\mathcal{A}} e^{t \|\varepsilon\|_2^2} \frac{1}{(2\pi)^{n/2}}e^{-\frac{1}{2}\| \varepsilon\|_2^2} \\ 
              &= \int_{\mathcal{A}} \frac{1}{(2\pi)^{n/2}}e^{-\frac{1}{2} (1-2t)\| \varepsilon\|_2^2} \\
              &= \frac{1}{(1-2t)^{n/2}}\int_{\mathcal{A}} \frac{(1-2t)^{n/2}}{(2\pi)^{n/2}}e^{-\frac{1}{2} (1-2t)\| \varepsilon\|_2^2} \\
              &= \frac{1}{(1-2t)^{n/2}} \, \mathbb{P}(\mathcal{Z} \in \mathcal{A}),
        \end{align*}
        where $\mathcal{Z} \sim N(0, (1-2t)I_n)$. 
        To calculate this value, set $X_i=\varepsilon_i$ for $i<n$, and $Y=-\varepsilon_n$.
        Since the distribution is spherical, $X_1,\dots,X_{n-1},Y$ are i.i.d.\ $\mathcal N(0,\sigma^2)$ with $\sigma^2=1-2t$.  
        The region can be written as
        \begin{equation*}
            \mathcal A=\{Y\ge0,\ 0\le X_i\le Y \ (i=1,\dots,n-1)\}.
        \end{equation*}
        Conditional on $Y=y\ge0$, we have:
        \begin{equation*}
            \mathbb P(\mathcal{Z} \in \mathcal A \mid Y=y) 
            = \big(\mathbb P(0\le X \le y)\big)^{n-1}
            = \big(\Phi(y/\sigma)-\tfrac12\big)^{n-1},
        \end{equation*}
        where $\Phi$ is the standard normal CDF.
        Hence
        \begin{align*}
           \mathbb P( \mathcal{Z} \in \mathcal A) 
            &= \int_0^\infty \big(\Phi(y/\sigma)-\tfrac12\big)^{n-1}\, f_Y(y)\,dy \\
            &=  \int_0^\infty \big(\Phi(u)-\tfrac12\big)^{n-1}\phi(u)\,du \\
            &= \int_0^{1/2} v^{n-1}\,dv 
            = \frac{1}{n}\left(\frac12\right)^n,
        \end{align*} 
        where $f_Y(y)=\tfrac1\sigma\phi(y/\sigma)$ and $\phi$ the standard normal PDF, and we used the variable changes $u=y/\sigma$ and $v=\Phi(u)-\tfrac12$.
        Writing all the chain of inequalities, we obtain:
        \begin{align*}
            \mathbb{E}_{\varepsilon\sim \mathcal{N}(0,I_n)}[e^{tZ^2_{c,x^*}}] 
            &\geq (1-\delta)\int_{\mathcal{A}}e^{ t \|\varepsilon\|_2^2} f(\varepsilon) \\
            &= (1-\delta)\left(\frac12\right)^n\frac{1}{n}\frac{1}{(1-2t)^{n/2}} \\
            & > M_{\Gamma\left(\tfrac{n-\epsilon}{2},\,2\right)}(t) =\frac{1}{(1-2t)^{\frac{n-\epsilon}{2}}}.
        \end{align*}
        
        The last inequality holds when
        \begin{equation*}
            t \in \left( \frac{1}{2} - \dfrac12 \Big(\dfrac{n 2^n}{1-\delta}\Big)^{-2/\epsilon}, \frac{1}{2} \right).
        \end{equation*}
        Since $\delta$ can be taken arbitrarily small by increasing $M$, this proves that for every $\varepsilon$ there exists a $t$ such that the Gamma distribution of $\rank(K) - \epsilon$ fails to dominate $Z^2_{c,x^*}$, and therefore $\rank(K)$ is the minimal parameter for which stochastic domination can hold in general.
    \end{proof}
    
    \begin{corollary}
        The non-negatively constrained problem $K = I_2$, $h = (1,-1)$, $x^* = (0, M)$ is a two-dimensional counterexample to the Burrus conjecture $Z^2_{c,x} \preceq \chi^2_1$ for $M > 0$ large enough. 
    \end{corollary}
    \begin{proof}
        The proof follows by taking $n = p = 2$ in the collection of counterexamples in \Cref{lemma:javier}.
    \end{proof}
    To the best of our knowledge, this is the first known counterexample to the Burrus conjecture in dimension smaller than $3$. 
\end{itemize}

\subsubsection{Test statistic with unconstrained second term \texorpdfstring{$\lambda^2_u$}{lambda2u}}

Defining $\mathcal{T}^2_{u,x}(\varepsilon) := \lambda^2_u(Hx, Kx+\varepsilon)$ with $\lambda^2_u(\mu, y) = \min_{Hx = \mu, Ax\leq b} \|Kx-y\|^2_2 - \min_{x} \|Kx-y\|^2_2$, we obtain after a similar computation as the previous subsection:
\begin{equation*}
    \mathcal{T}^2_{u,x}(\varepsilon)= \min_{Hz = 0, Az\leq b-Ax} \|Kz -\varepsilon\|^2_2 -\min_{z} \|Kz -\varepsilon\|^2_2.
\end{equation*}
Let $Q^2_u(x)$ be the $1-\alpha$ quantile of $\mathcal{T}^2_{u,x}(\varepsilon)$ for $\varepsilon \sim \mathcal{N}(0,I)$.
We summarize some relevant properties of this function below.

\begin{itemize}[leftmargin=0em]
    \item 
    \emph{Convexity.}
    We have the following important convexity result:
    \begin{theorem}
        [Convexity of the quantile function $Q^2_u(x)$]
        \label{thm:convexityl2u}
        For any fixed $0 < \alpha < 1$, the quantile function $Q^2_u(x)$ at fixed level $1-\alpha$ is a convex function of $x$, both for linear and cone constraints. 
    \end{theorem}
    \begin{proof}
        It is enough to show that $(x, \varepsilon) \mapsto \mathcal{T}^2_{u, x}(\varepsilon)$ is jointly convex in $(x, \varepsilon)$. 
        If that is the case, \cite[Lemma 2.28]{Kibzun1995-oz} directly shows that the quantile function is convex. 
        An alternative reasoning comes from Pr\'{e}kopa's theorem \cite{Prekopa1971,Prekopa1973} which implies that log-concavity is preserved by marginalization: joint convexity implies $\mathcal{S}_t:=\{(x,\varepsilon):\ \mathcal{T}^2_{u,x}(\varepsilon)\le t\}$ is a convex set (therefore its indicator function is log-concave), and if $f(\varepsilon)$ is any log-concave density, $h(x,t,\varepsilon):=\mathbf{1}_{\mathcal{S}_t}(x,\varepsilon)\,f(\varepsilon)$ is product of log-concave functions and hence log-concave.
        Then, by the Pr\'{e}kopa's theorem $\varphi(x,t):=\int h(x,t,\varepsilon)\,d\varepsilon =\mathbb{P}\!\left(\mathcal{T}^2_{u,x}(\varepsilon)\le t\right)$
        is log-concave in $(x,t)$. 
        Therefore, the superlevel sets $\{(x,t):\ \varphi(x,t)\ge 1-\alpha\}$ are convex, which correspond exactly to the epigraph of the quantile function, making $Q^2_u$ convex.
        
        We now show that $(x, \varepsilon) \mapsto \mathcal{T}^2_{u, x}(\varepsilon)$ is jointly convex in $(x, \varepsilon)$, finishing the proof. 
        Let $\hat{z}(\varepsilon)  = K^\dagger \varepsilon\in \argmin_{z} \|Kz-\varepsilon\|^2_2$. 
        Then, for any $z$, by the Pythagorean theorem and orthogonality of $\varepsilon - K \hat{z}$ to the range of $K$, we have: 
        \begin{equation*}
            \|Kz-\varepsilon\|^2_2 = \|Kz-K\hat{z}\|^2_2 + \|K\hat{z}-\varepsilon\|^2_2,
        \end{equation*}
        which we used earlier in \eqref{eq:least_squares_eq} to relate the confidence regions around the least square estimator with inverted regions of $\lambda^2_u$. 
        Therefore, the difference cancels the constant $\|K\hat{z}-\varepsilon\|^2_2$, and we have:
        \begin{equation*}
            \mathcal{T}^2_{u,x}(\varepsilon)
            = \min_{Hz=0,Az\le b-Ax} \|K(z-\hat{z}(\varepsilon))\|^2_2.
        \end{equation*} 
        Write $\mathcal{Z}(x) := \{z \ | Hz = 0, Az \leq b-Ax\}$ and $g(u, x) := \min_{z \in \mathcal{Z}(x)} \|K(z-u)\|^2_2$. 
        The map $(z,u)\mapsto \|K(z-u)\|_2^2$ is jointly convex, and the feasible set $\mathcal{Z}(x)$ depends affinely on $(z,x)$. 
        By standard results on convex analysis, $g(u, x)$ is jointly convex in $(u,x)$. 
        Since $\hat z(\varepsilon)$ is affine in $\varepsilon$, the composition $\mathcal{T}^2_{u,x}(\varepsilon) = g(\hat z(\varepsilon),x)$ is jointly convex in $(x, \varepsilon)$. 
        The same argument applies for $\mathcal{Z}(x) := \{z \ | Hz = 0, z+x \in \mathcal{X}\}$ where $\mathcal{X}$ is any convex set, so in particular, the result applies for convex constraints as well.  
    \end{proof}
    
    \begin{lemma}
        \label{lemma:Extreme_ray_decrease_2u}
        Let $y$ be in the recession cone of the constraint set ($Ay \leq 0$ for linear constraints $Ax^* \leq b$ and all of the convex cone $\mathcal{C}$ for $x^* \in \mathcal{X}$), and let $x$ be any point. 
        Then, $Q^2_u(x+y) \leq Q^2_u(x)$.
    \end{lemma}
    \begin{proof}
        For the linearly constrained case, we have, for every fixed $\varepsilon$:
        \begin{equation}
            \label{eq:recession_inequality_one_term}
            \min_{Hz = 0, Az \leq b-Ax-Ay} \|Kz-\varepsilon\|^2_2 \leq \min_{Hz = 0, Az \leq b-Ax} \|Kz-\varepsilon\|^2_2,
        \end{equation}
        since the feasible set enlarges when $-Ay \ge 0$ (i.e., $Ay \le 0$). It therefore directly holds from the definition of $\mathcal{T}^2_{u, x}(\varepsilon)$ that 
        $\mathcal{T}^2_{u, x+y}(\varepsilon) \leq \mathcal{T}^2_{u, x}(\varepsilon)$ for all $\varepsilon$.
        The quantile inequality then follows. 
        
        Similarly, for the cone constraint case, we have:
        \begin{equation}
            \label{eq:recession_inequality_one_term_cone}
            \min_{Hz = 0, (z+x+y)\in \mathcal{C}} \|Kz-\varepsilon\|^2_2 \leq  \min_{Hz = 0, (z+x) \in \mathcal{C}} \|Kz-\varepsilon\|^2_2.
        \end{equation}
        Because $\mathcal{C}-x \subseteq \mathcal{C}-x-y$ for any $x,y \in \mathcal{C}$, the same conclusion follows.
    \end{proof}

    \item
    \emph{Maximization.}
    \Cref{thm:convexityl2u}, together with \Cref{lemma:Extreme_ray_decrease_2u}, have important consequences for locating where the quantile optimization problems \eqref{eq:max_q_mu_md_2} and \eqref{eq:max_q_md_2} are maximized. 
    We separate into the cases of linear and cone constraints.
    
    \begin{theorem}
        [Maximization of $Q^2_u(x)$ with linear constraints]
        \label{thm:maximizing_2u}
        In the linearly constrained case $Ax^* \leq b$, consider 
        \begin{equation*}
            \lambda^2_u(\mu, y) = \min_{Hx = \mu, Ax\leq b} \|Kx-y\|^2_2 - \min_{x} \|Kx-y\|^2_2
        \end{equation*}
        with its corresponding quantile $Q^2_u(x)$.
        Let $\mathcal{P}$ be a polyhedron of the form $\{x: Ax \leq b, Gx \leq d\}$, where $A$ and $b$ are the original constraints on $x$, and $G$, $d$ are arbitrary. 
        Let $\{p_i\}_{i = 1}^m$ be the set of extreme points of $\mathcal{P}$. 
        Then, $\sup_{x \in \mathcal{C}} Q^2_u(x)$ is achieved in $\mathcal{P}$, i.e.,
        \begin{equation}
            \sup_{x \in \mathcal P} Q^2_u(x) = \max_{i = 1:m} Q^2_u(p_i).
        \end{equation}
    \end{theorem}
    \begin{proof}
        Let $x \in \mathcal{P}$ and use the Weyl–Minkowski theorem to express
        \begin{equation*}
        x = \sum_{i=1}^{m} \lambda_i p_i + \sum_{j=1}^{\ell} \mu_j r_j,
        \quad
        \lambda_i \ge 0,\ \sum_{i=1}^{m} \lambda_i = 1,\ \ \mu_j \ge 0,\ \ Ar_j \le 0,\ Gr_j \le 0.
        \end{equation*}
        Since $A\sum_{j=1}^\ell \mu_jr_j \leq 0$, using \Cref{lemma:Extreme_ray_decrease_2u} and then \Cref{thm:convexityl2u}, we have that
        \begin{align*}
            Q^2_u(x) 
            &= Q^2_u\Big(\sum_{i=1}^m \lambda_ip_i + \sum_{j=1}^\ell \mu_jr_j\Big) \\
            &\leq Q^2_u\Big(\sum_{i=1}^m \lambda_ip_i\Big) \\
            &\leq \sum_{i = 1}^m \lambda_iQ^2_u(p_i) \\
            &\leq \sum_{i = 1 }^m \lambda_i \max_{j =1:m} Q^2_u(p_j) = \max_{j =1:m} Q^2_u(p_j).
        \end{align*}
        Since $x \in \mathcal P$ was arbitrary, this concludes the proof.
    \end{proof}
    Thanks to this result, for linear constraints, we can obtain the optimal thresholds as:
    \begin{align*}
        d^*(\mu) &= \sup_{Hx = \mu, Ax \leq b} Q^2_u(x) = \max_{p_i \in E(Ax \leq b, Hx =\mu)} Q^2_u(p_i), \\
        D^* &= \sup_{Ax \leq b} Q^2_u(x) = \max_{q_i \in E(Ax \leq b)} Q^2_u(q_i).
    \end{align*}

    \begin{theorem}
        [Global maximization of $Q^2_u(x)$ with cone constraints]
        \label{thm:maximizing_2u_cone}
        In the cone constrained case $x^* \in \mathcal{C}$, consider 
        \begin{equation*}
            \lambda^2_u(\mu, y) = \min_{Hx = \mu, x \in \mathcal{C}} \|Kx-y\|^2_2 - \min_{x} \|Kx-y\|^2_2
        \end{equation*}
        with its corresponding quantile $Q^2_u(x)$.
        Then, $\sup_{x \in \mathcal{C}} Q^2_u(x)$ is achieved at $x = 0$, i.e., 
        \begin{equation}
            \label{eq:global_max_l2u}
            \sup_{x \in \mathcal{C}} Q^2_u(x) = Q^2_u(0).
        \end{equation}
        Furthermore, letting $\Phi_\mu = \{ x \in \mathcal{C} | Hx = \mu \}$, and $\operatorname{rb}(\Phi_\mu)$ denote its relative boundary, we have:
        \begin{equation}
            \sup_{x\in\Phi_\mu} Q(x)=\sup_{x\in \operatorname{rb}(\Phi_\mu)} Q(x).
        \end{equation}
    \end{theorem}
    \begin{proof}
        The global maximization result \eqref{eq:global_max_l2u} comes directly from $\eqref{eq:recession_inequality_one_term_cone}$, since the recession cone of $\mathcal{C}$ is itself and $Q^2_u$ satisfies $Q^2_u(x+y) \leq Q^2_u(x)$ for all $x, y \in \mathcal{C}$.
        Taking $x = 0$ and $y \in \mathcal{C}$ arbitrary yields the desired result.
        For the sliced result, if $\Phi_\mu=\operatorname{rb}(\Phi_\mu)$, there is nothing to prove. 
        Otherwise, take any $x\in \Phi_\mu\setminus \operatorname{rb}(\Phi_\mu)$, and let us show that $Q^2_u(x) \leq \sup_{z \in \operatorname{rb}(\Phi_\mu)}Q^2_u(z)$. 
        
        Consider a nonzero direction $u\in Ker(K)$ and the maximal interval $I:=\{t\in\mathbb R:\ x+tu\in \Phi_\mu\}$.
        There are two cases:
        \begin{itemize}
            \item 
            If $I=(t_-,t_+)$ with finite endpoints, then the endpoints $a:=x+t_-u$ and $b:=x+t_+u$ belong to $\operatorname{rb}(\Phi_\mu)$, and by convexity of $Q$ along the line $x+tu$, $Q(x)\le\max\{Q(a),Q(b)\}\le\sup_{z\in \operatorname{rb}(\Phi_\mu)}Q(z)$. 
            \item
            If one endpoint is $+\infty$ or $-\infty$ (suppose $t_+=+\infty$ without loss of generality), then $u\in\operatorname{rec}(\Phi_\mu)=\mathcal C\cap Ker(K)$.
            By the cone monotonicity $Q(x+tu)\le Q(x)$ for all $t\ge 0$, hence again $Q(x)\le \sup_{z\in \operatorname{rb}(\Phi_\mu)}Q(z)$ (the other endpoint of $I$ is finite and lies on $\operatorname{rb}(\Phi_\mu)$, or the inequality is trivial if both directions are recession.
            Since $x\in \Phi_\mu\setminus \operatorname{rb}(\Phi_\mu)$ was arbitrary, $\sup_{z \in \Phi_\mu}Q^2_u(z)\le \sup_{z \in\operatorname{rb}(\Phi_\mu)}Q^2_u(z)$. 
        \end{itemize}
        The reverse inequality is immediate because $\operatorname{rb}(\Phi_\mu)\subseteq \Phi_\mu$.
        This concludes the proof.
    \end{proof}
    
    \begin{corollary}
        [Non-negative constraints]
        \label{cor:non-negative-constraints}
        Under non-negativity constraints $x^*\geq 0$, we have
        \begin{equation}
            \sup_{x\geq 0} Q^2_u(x) = Q^2_u(0),
        \end{equation}
        and if $k = 1$ and the functional of interest is $h^\top x^*$,  
        \begin{equation}
            \label{eq:mu_argmaxes}
            \argmax_{x \geq 0, h^\top x = \mu} Q^2_u(x) \in \{(\mu/h_i)  e_i, i: \mu/h_i > 0\} 
        \end{equation}
        where $e_i$ are the standard basis vectors.
    \end{corollary}
    \begin{proof}
        Directly from either \Cref{thm:maximizing_2u} or \Cref{thm:maximizing_2u_cone}, as the non-negative constraints are both linear constraints and cone constraints.
    \end{proof}
    A similar result holds with more functionals, but the number of extreme points to check grows as at most $\binom{p}{k}$ with $k$ functionals, which can become unfeasible to manually search even with $k = 2$ if $p$ is large.
    
    \begin{corollary}
        [Optimal constant for non-negative $\lambda^2_u$ constructions in \Cref{table:rosetta_stone}]
        \label{cor:improved_constants_2u}
        In the non-negatively constrained case, consider the $\mu$-description confidence region in \Cref{table:rosetta_stone}:
        \begin{equation}
            \label{eq:mu_desc_nn}
            \{\mu\in\mathbb R^k:\min_{\substack{Hx=\mu,\\ x \geq 0}}\|Kx-y\|_2^2
            -\min_{x \in \mathbb R^p}\|Kx-y\|_2^2
            \le \delta^2_u\},
        \end{equation}
        and its related $x$-description bounding box:
        \begin{equation*}
        \displaystyle \prod_{i = 1}^k
        \quad \begin{aligned}
         \min_{x}/\max_{x} \quad & h_i^\top x\\
          \st  \quad &  \|Kx-y\|^2_2 &\leq \min_{x '\in \mathbb R^p}\|Kx'-y\|_2^2 + \delta^2_u\\ 
           & x \geq 0.
        \end{aligned}
        \end{equation*}
        Then, the $1-\alpha$ quantile of 
        \begin{equation}
           \min_{Hz = 0, z \geq 0 } \|Kz-\varepsilon\|^2_2 - \min_{z} \|Kz-\varepsilon\|^2_2 = \min_{H z = 0, z \geq 0 } \|K(z-\hat{z}(\varepsilon))\|^2_2, \; \hat{z}(\varepsilon) \in \argmin_z \|Kz-\varepsilon\|^2_2 
        \end{equation}
        with $\varepsilon \sim \mathcal{N}(0, I)$, is the optimal (smallest) $\delta^2_u$ satisfying frequentist coverage for all $x \geq 0$ in \eqref{eq:mu_desc_nn}.
    \end{corollary}
    \begin{proof}
        By using \Cref{cor:non-negative-constraints} and the description of optimal constants given by \eqref{eq:max_q_mu_md} and \eqref{eq:max_q_md}.
    \end{proof}
    
    In particular, for $k = 1$, consider the construction \eqref{eq:burrus3a} or \eqref{eq:burrus3b}, which is of the form
    \begin{equation}
        \label{eq:burrus3b_delta}
        \begin{aligned}
            \min_{x}/\max_{x} \quad & h^\top x\\
            \st  \quad & \|Kx-y\|^2_2 \leq \min_{x'}\|Kx'-y\|^2_2 + \delta\\ 
            & x \geq 0
        \end{aligned}
    \end{equation}
    with $\delta = Q_{\chi^2_r,1-\alpha}$. 
    Then, the optimal constant $\delta$ such that \eqref{eq:burrus3b_delta} has $1-\alpha$ frequentist coverage for every $x \geq 0$ is the $1-\alpha$ quantile of $Z^2_u(x =0)$, i.e., the quantile of 
    \begin{equation}
       \min_{h^\top z = 0, z \geq 0 } \|Kz-\varepsilon\|^2_2 - \min_{z} \|Kz-\varepsilon\|^2_2 = \min_{h^\top z = 0, z \geq 0 } \|K(z-\hat{z}(\varepsilon))\|^2_2, \; \hat{z}(\varepsilon) \in \argmin_z \|Kz-\varepsilon\|^2_2 
    \end{equation}
    with $\varepsilon \sim \mathcal{N}(0,I)$, which is always equal or smaller than $Q_{\chi^2_r,1-\alpha}$.
    In \Cref{sec:chi_bar_sq}, we identify this random variable to follow a chi-bar-squared distribution. 
    An improvement is possible by allowing $\delta$ to depend on $x$, at the expense of destroying the convexity of the optimization problem, by taking $\delta(x) := \max_{z \in  \{(h^\top x/h_i) e_i,  i: h^\top x/h_i > 0\}} Q^2_u(z)$.

    \item
    \emph{Bounds.}
    We note that \Cref{lemma:javier}, together with \eqref{eq:chain_stochastic_dominance}, shows that there is no better uniform (in the sense of holding for all $K,H,x$) stochastic dominance bound than $Z^2_{u,x} \preceq \chi^2_{\rank(K)}$.
\end{itemize}

\subsubsection{One term test statistic \texorpdfstring{$\lambda^1$}{lambda1}}

Our final test statistic under consideration is $\lambda^1(\mu, y) = \min_{Hx = \mu, Ax\leq b} \|Kx-y\|^2_2$, which we will show has many similar properties to $\lambda^2_u$. 
It easily follows from the same computations in the previous subsections that 
\begin{equation*}
    \mathcal{T}^1_{x}(\varepsilon)= \min_{Hz = 0, Az\leq b-Ax} \|Kz -\varepsilon\|^2_2.
\end{equation*}
Let $Q^1(x)$ as the $1-\alpha$ quantile of $\mathcal{T}^1_{x}(\varepsilon)$ for $\varepsilon \sim \mathcal{N}(0,I)$.
We summarize below that similar properties of this function as with $Q^2_u$.

\begin{itemize}[leftmargin=0em]
    \item 
    \emph{Convexity.}
    We have the following convexity result, similar to \Cref{thm:convexityl2u}.
    \begin{theorem}
        \label{thm:convexity_Q1}
        For any fixed $0 < \alpha < 1$, the quantile function $Q^1(x)$ is a convex function of $x$.
    \end{theorem}
    \begin{proof}
        Using the same proof argument as in \Cref{thm:convexityl2u}, it suffices to show that $(x, \varepsilon) \mapsto \mathcal{T}^1_{x}(\varepsilon)$ is jointly convex in $(x, \varepsilon)$. 
        This is proven similarly to $\mathcal{T}^2_{u}$.
        Define $\mathcal{Z}(x) := \{z \ | Hz = 0, Az \leq b-Ax\}$ and $\mathcal{T}^1_{x}(\varepsilon) := \min_{z \in \mathcal{Z}(x)} \|Kz-\varepsilon\|^2_2$. 
        The map $(z,\varepsilon)\mapsto \|Kz-\varepsilon\|_2^2$ is jointly convex, and the feasible set $\mathcal{Z}(x)$ depends affinely on $(z,x)$, so the result follows.
    \end{proof}

    \begin{lemma}
        \label{lemma:Extreme_ray_decrease_1} 
        Let $y$ be in the recession cone of the constraint set ($Ay \leq 0$ for linear constraints $Ax^* \leq b$ and all of the convex cone $\mathcal{C}$ for $x^* \in \mathcal{X}$), and let $x$ be any point. 
        Then, $Q^1(x+y) \leq Q^1(x)$.
    \end{lemma}
    \begin{proof}
        Follows directly from the inequalities \eqref{eq:recession_inequality_one_term} and \eqref{eq:recession_inequality_one_term_cone} established in the proof of \Cref{lemma:Extreme_ray_decrease_2u}.
    \end{proof}
    
    \item
    \emph{Maximization.}
    \Cref{thm:convexity_Q1} along with \Cref{lemma:Extreme_ray_decrease_1} lead to the following consequence regarding the locations where the underlying quantile optimization problems are maximized.
    
    \begin{theorem}
        \label{thm:maximizing_one_term}
        \Cref{thm:maximizing_2u,thm:maximizing_2u_cone} apply verbatim for the maximizers of $Q^1$ in the linearly constrained and cone constrained problems, respectively.
    \end{theorem}
    \begin{proof}
        The proof of \Cref{thm:maximizing_2u,thm:maximizing_2u_cone} only use convexity of $Q^2_u(x)$ (\Cref{thm:convexityl2u}) and the recession property in \Cref{lemma:Extreme_ray_decrease_2u}, properties which also apply to $Q^1$ (\Cref{thm:convexity_Q1} and \Cref{lemma:Extreme_ray_decrease_1}, respectively).
    \end{proof}
    Thus, we can directly invoke the equivalent of \Cref{cor:improved_constants_2u} for constructions based on $\lambda^1$.
        \begin{corollary}
        [Optimal constant for non-negative $\lambda^1$ constructions in \Cref{table:rosetta_stone}]
        \label{cor:improved_constants_1}
        In the non-negatively constrained case, consider the $\mu$-description confidence region in \Cref{table:rosetta_stone}:
        \begin{equation}
            \label{eq:mu_desc_nn2}
            \{\mu\in\mathbb R^k:\min_{\substack{Hx=\mu,\\ x \geq 0}}\|Kx-y\|_2^2
           \le \delta^1\},
        \end{equation}
        and its related $x$-description bounding box:
        \begin{equation*}
        \displaystyle \prod_{i = 1}^k
        \quad \begin{aligned}
         \min_{x}/\max_{x} \quad & h_i^\top x\\
          \st  \quad &  \|Kx-y\|^2_2 &\leq \delta^1\\ 
           & x \geq 0.
        \end{aligned}
        \end{equation*}
        Then, the $1-\alpha$ quantile of 
           $\min_{Hz = 0, z \geq 0 } \|Kz-\varepsilon\|^2_2$
        with $\varepsilon \sim \mathcal{N}(0, I)$, is the optimal (smallest) $\delta^1$ satisfying frequentist coverage for all $x \geq 0$ in \eqref{eq:mu_desc_nn2}.
    \end{corollary}
    \begin{proof}
        By using \Cref{thm:maximizing_one_term}, we can show that \Cref{cor:non-negative-constraints} applies for $Q^1(x)$, then the description of optimal constants given by \eqref{eq:max_q_mu_md} and \eqref{eq:max_q_md} gives the desired result.
    \end{proof}
    
    In particular, for $k = 1$, consider the construction \eqref{eq:starkSSB} which is of the form
    \begin{equation}
        \label{eq:ssb_delta}
        \begin{aligned}
            \min_{x}/\max_{x} \quad & h^\top x\\
            \st  \quad & \|Kx-y\|^2_2 \leq \delta\\ 
            & x \geq 0
        \end{aligned}
    \end{equation}
    with $\delta = Q_{\chi^2_n,1-\alpha}$. 
    Then, the optimal constant $\delta$ such that \eqref{eq:ssb_delta} has $1-\alpha$ frequentist coverage for every $x \geq 0$ is the $1-\alpha$ quantile of $Z^1(x =0)$, i.e., the quantile of $\min_{h^\top z = 0, z \geq 0 } \|Kz-\varepsilon\|^2_2$ with $\varepsilon \sim \mathcal{N}(0,I)$, which is always equal or smaller than $Q_{\chi^2_n,1-\alpha}$.
\end{itemize}

\subsubsection[Comparison between two statistics, and their distributions at the origin]{Comparison between $\lambda^2_u$ and $\lambda^1$, and their distributions at the origin} 
\label{sec:chi_bar_sq}

As derived in the previous subsections, both $\lambda^2_u$ and $\lambda^1$ have identical behavior with respect to the locations of maximizers of their induced quantile functions.
We now give a preliminary comparison to guide the choice between them in practice.

When $K$ is surjective (full row rank), we have $\min_{x} \|Kx-y\|^2_2 = 0$ for all $y$, so $\lambda^2_u = \lambda^1$. 
We therefore focus on the non-surjective case, under cone constraints, and for global quantile optimization. A comparison under general linear constraints $Ax^*\le b$ and/or in the sliced setting remains open.

For a fixed $y$, the resulting confidence regions are:
\begin{align*}
    \mathcal{R}^1 &= \{ \mu: \min_{Hx = \mu, x \in \mathcal{C}} \|Kx-y\|^2_2 \leq D^1\}, \\
    \mathcal{R}^2_u &= \{ \mu: \min_{Hx = \mu, x \in \mathcal{C}} \|Kx-y\|^2_2 \leq \min_{x} \|Kx-y\|^2_2  + D^2_u \}, 
\end{align*}
for given constants $D^1$ and $D^2_u$. 
Coverage holds whenever $D^1 \geq \max_{x \in \mathcal{C}} Q^1(x) = Q^1(0)$ and $D^2_u \geq \max_{x \in \mathcal{C}} Q^2_u(x) = Q^2_u(0)$. 
For simplicity, we begin analyzing not the optimal constants given by their quantiles at $x = 0$, but their upper bounds independent of $K$ and $H$ given by the stochastic dominance chain $\eqref{eq:chain_stochastic_dominance}$, i.e., $D^1 = Q_{\chi^2_n,1-\alpha}$, and $D^2_u = Q_{\chi^2_r,1-\alpha}$

Observe that the distribution of $\min_{x} \|Kx-y\|^2_2$ when $y \sim \mathcal{N}(Kx^*, I)$ is the same for all $x^* \in \mathcal{C}$. Indeed, writing $y = Kx^* + \varepsilon$, we can translate the optimization variables without changing the optimal value.
We therefore may take $x^* = 0$, and study the distribution of $\min_{x} \|Kx-\varepsilon\|^2_2$ with $\varepsilon \sim \mathcal{N}(0,I)$.
A standard argument shows that this distribution is $\chi^2$ distribution with $n-r$ degrees of freedom.   
Thus, with $\xi_{n,r} \sim \chi^2_{n-r}$, we have $\mathcal{R}^1 \subseteq \mathcal{R}^2_u$ if and only if $Q_{\chi^2_n,1-\alpha} \leq \xi_{n, r} + Q_{\chi^2_r,1-\alpha}$. 

\begin{figure}[!t]
    \centering
    \includegraphics[width=0.75\linewidth]{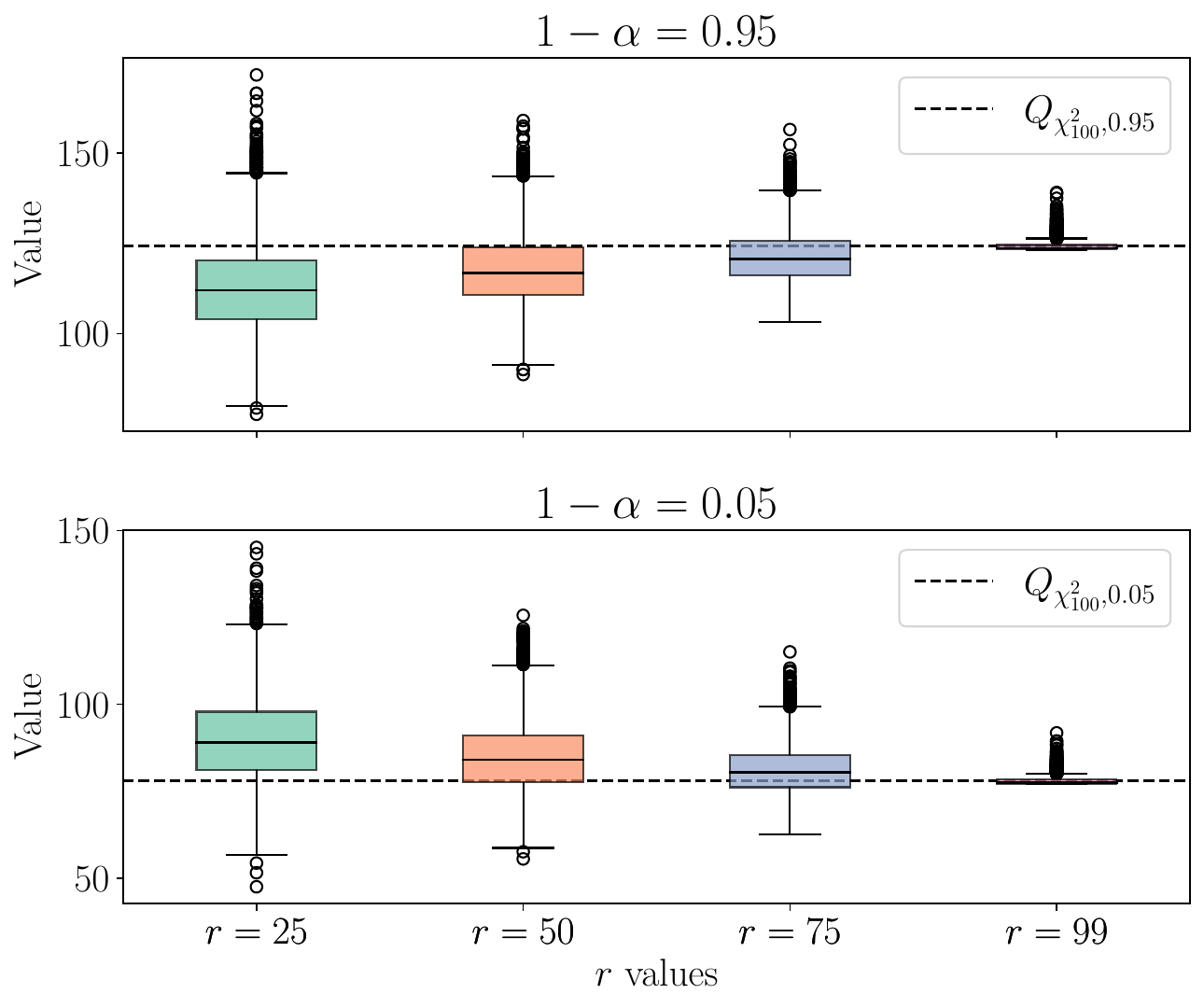}
    \caption{Comparison between the fixed chi-square quantile 
    $Q_{\chi^2_n,\,1-\alpha}$ (dashed line) and the distribution 
    $Q_{\chi^2_r,\,1-\alpha}+\xi_{n,r}$ (boxplots), for $n=100$ and 
    $r\in\{25,50,75,99\}$. The top panel corresponds to $1-\alpha=0.05$, 
    and the bottom panel to $1-\alpha=0.25$.}
    \label{fig:chicomparison}
\end{figure}

Define the quantile gap $\Delta(n,r,\alpha) = Q_{\chi^2_n,1-\alpha}-Q_{\chi^2_r,1-\alpha}$. 
With $\xi_{n,r}\sim\chi^2_{\,n-r}$, we have $\mathcal{R}^1 \subseteq \mathcal{R}^2_u$.
This is equivalent to $\xi_{n,r} \geq \Delta(n,r,\alpha)$, which happens with probability $p(n,r,\alpha)$. Using the Cornish-Fisher expansion \cite{cornish1938moments,fisher1960percentile}, we can approximate this quantity as $Q_{\chi^2_k,1-\alpha} = k + z_{1-\alpha}\sqrt{2k}+O(1)$ (with $z_{1-\alpha}$ the $1-\alpha$ normal quantile).
This yields the approximation:
\begin{equation}
    \label{eq:p-asymptotic}
    p(n,r,\alpha)\approx 1-\Phi\!\bigl(z_{1-\alpha}\,c(\rho)\bigr),
\end{equation}
valid up to $O(n^{-1/2})$ corrections, where $\rho := r/n$ and
$
    c(\rho) := \frac{1-\sqrt{\rho}}{\sqrt{1-\rho}}.
$
From \eqref{eq:p-asymptotic}, we read off three asymptotic regimes for $p(n,r,\alpha)$:
(i) if $r$ is fixed and $n \to \infty$, then
$
    p(n,r,\alpha) \rightarrow \alpha,
$
(ii) if $r/n \to \rho \in (0,1)$, then
$
    p(n,r,\alpha) \rightarrow 1-\Phi\!\bigl(z_{1-\alpha}c(\rho)\bigr),
$
and (iii) if the difference $s=n-r$ remains fixed while $n,r \to \infty$, then
$
    p(n,r,\alpha) \rightarrow \tfrac{1}{2}.
$
In all cases, for small $\alpha$, we typically have $p \leq 1/2$, so the random threshold $\mathcal{R}^2_u$, $Q_{\chi^2_{r, 1-\alpha}} + \xi_{n,r}$ is often smaller than the fixed threshold $Q_{\chi^2_{n, 1-\alpha}}$ in $\mathcal{R}^1$, at the cost of a worse worst-case over the data. 
This behavior is illustrated in \Cref{fig:chicomparison}: for small $\alpha$, most of the mass of the distribution lies below the fixed constant, and for larger $\alpha$, this relation flips. 

We next show that, at the origin $x = 0$ (where the quantiles are maximized over closed convex cones), the distributions of $Z^2_{u, x}$ and $Z^{1}_{x}$ to chi-bar-squared distributions, mixtures of $\chi^2$ distributions with different degrees of freedom. 
Importantly, their weights are related, so some conclusions of the previous subsection carry over to this case, as the distribution at the origin for $\lambda^1$ is a mixture of $\chi^2$ distributions whose degrees of freedom (but not weights) get shifted. 
We begin with a simple two-dimensional example to illustrate:
\begin{example}[Two-dimensional example]
    For $n=2$, $p=1$, $K=\begin{bmatrix}1\\0\end{bmatrix}$, and $F=\mathbb R_+$, the cone 
    $S=\{(t,0):t\ge0\}$ gives
    \[
        Z^2_{u,0}\sim\tfrac12\chi^2_0+\tfrac12\chi^2_1,\qquad 
        Z^1_0\sim\tfrac12\chi^2_1+\tfrac12\chi^2_2.
    \]
\end{example}
Therefore, when using the bounds, the choice would be between the quantile of a $\chi^2_2$ and the quantile of $\chi^2_1$ plus a draw of a $\chi^2_1$, and with the exact distributions, the fixed threshold of $\lambda^1$ is the quantile of $\tfrac12\chi^2_1+\tfrac12\chi^2_2$, and the random threshold of $\lambda^2_u$ is the quantile of $\tfrac12\chi^2_0+\tfrac12\chi^2_1$ plus a draw of a $\chi^2_1$. 
We prove below a general statement about the laws of 
$Z^1_0$ and $Z^2_{u, 0}$. 
While we expect similar conclusions to the geometry-independent bounds to apply, we defer their theoretical analysis to future work.
\begin{lemma}[Chi-bar-squared laws for $Z^1_0$ and $Z^2_{u,0}$]
    Let $L: =\{x\in\mathbb R^p:Hx=0\}$ and $F:=\mathcal C\cap L$ be a closed convex cone. 
    Let $R:=\mathrm{range}(K)\subseteq\mathbb R^n$ with $r:=\dim(R)$, and define the image cone $S:=K(F)\subseteq R$.
    Consider
    \begin{equation*}
        \mathcal{T}^1_0(\varepsilon):=\min_{x\in F}\,\|Kx-\varepsilon\|_2^2,\qquad
        \mathcal{T}^2_{u,0}(\varepsilon):=\min_{x\in F}\,\|Kx-\varepsilon\|_2^2-\min_{x\in\mathbb R^p}\,\|Kx-\varepsilon\|_2^2,
    \end{equation*}
    and let $Z^1_0, Z^2_{u,0}$ be the corresponding random variables for $\varepsilon \sim \mathcal N(0,I_n)$.
    Then, we have:
    \begin{enumerate}
        \item[(i)] 
        \begin{equation*}
            \mathcal T^1_0 = \mathrm{dist}^2(\varepsilon,S), \quad  Z^1_0\stackrel{\textup{d}}{=}\sum_{j=0}^{n} v_{\,n-j}^{(\mathbb R^n)}(S)\,\chi^2_{j},
        \end{equation*}
        where $\{v_k^{(\mathbb R^n)}(S)\}_{k=0}^n$ are the conic intrinsic volumes of $S$ as a subset of $\mathbb R^n$, which are constants depending on $\mathcal{C}, H$, and $K$.
        \item[(ii)] Writing $\varepsilon_R:=P_R\varepsilon\sim\mathcal N(0,I_r)$, we have:
        \begin{equation*}
            \mathcal T^2_{u,0} = \mathrm{dist}^2(\varepsilon_R,S), \quad Z^2_{u, 0}\stackrel{\textup{d}}{=}\sum_{j=0}^{r} v_{\,r-j}^{(R)}(S)\,\chi^2_{j},
        \end{equation*}
        where $\{v_k^{(R)}(S)\}_{k=0}^r$ are the intrinsic volumes of $S$ computed in the ambient space $R$.
        Moreover, the two sets of intrinsic volumes satisfy:
        \[
        v_k^{(\mathbb R^n)}(S)=v_k^{(R)}(S)\quad\text{for }k=0,\dots,r,
        \qquad
        v_k^{(\mathbb R^n)}(S)=0\quad\text{for }k>r.
        \]
    \end{enumerate}
\end{lemma}
\begin{proof}
    Because $\{Kx:x\in F\}=S$ and $\mathrm{range}(K)=R$, we have:
    \begin{equation}
        \label{eq:ids}
        \mathcal T^1_0(\varepsilon)=\min_{x\in F}\|Kx-\varepsilon\|_2^2=\mathrm{dist}^2(\varepsilon,S),\qquad
        \min_{x\in\mathbb R^p}\|Kx-\varepsilon\|_2^2=\mathrm{dist}^2(\varepsilon,R)=\|(I-P_R)\varepsilon\|_2^2.
    \end{equation}
    
    (i) Since $S$ is a closed convex cone in $\mathbb R^n$, Moreau’s decomposition gives
    $\varepsilon=\Pi_S\varepsilon+\Pi_{S^\circ}\varepsilon$ with orthogonal parts, hence
    $\mathrm{dist}^2(\varepsilon,S)=\|\Pi_{S^\circ}\varepsilon\|_2^2$. 
    For Gaussian $\varepsilon$, the
    distribution of $\|\Pi_{S^\circ}\varepsilon\|_2^2$ is a chi-bar-squared mixture:
    \begin{equation}
        \label{eq:chisqbarS_rn}
        \mathrm{dist}^2(\varepsilon,S)\sim \sum_{j=0}^n v_{n-j}^{(\mathbb R^n)}(S)\chi^2_j,
    \end{equation}
    where the weights are the conic intrinsic volumes of $S$ in $\mathbb R^n$. 
    (For more details about the chi-bar-squared distribution, see \cite{Silvapulle2001}.)
    
    (ii) Since $S\subseteq R$, subtracting the identities in \eqref{eq:ids} gives
    \begin{equation}
        \label{eq:orth_noise_decomposition}
        \mathcal T^2_{u,0}(\varepsilon)=\mathrm{dist}^2(\varepsilon,S)-\mathrm{dist}^2(\varepsilon,R)
        =\min_{s\in S}\|\varepsilon-s\|_2^2-\|\varepsilon_{R^\perp}\|_2^2
        =\min_{s\in S}\|\varepsilon_R-s\|_2^2
        =\mathrm{dist}^2(\varepsilon_R,S),
    \end{equation}
    where we have written the orthogonal decomposition $\mathbb R^n=R\oplus R^\perp$ and $\varepsilon_R\sim\mathcal N(0,I_r)$ and $\varepsilon_{R^\perp} \sim \mathcal N(0, I_{n-r})$. 
    Viewing $S$ as a cone in $R\simeq\mathbb R^r$, the same reasoning as in (i) yields
    \begin{equation*} 
        \mathrm{dist}^2(\varepsilon_R,S)\sim\sum_{j=0}^r v_{r-j}^{(R)}(S)\chi^2_j.
    \end{equation*}
    It remains to relate $v^{(\mathbb R^n)}$ and $v^{(R)}$. 
    From \eqref{eq:orth_noise_decomposition}, 
    \begin{equation*}
        \mathrm{dist}^2(\varepsilon,S) =\min_{s\in S}\|\varepsilon_R-s\|_2^2+\|\varepsilon_{R^\perp}\|_2^2.
    \end{equation*}
    When $\varepsilon$ is randomized, $\|\varepsilon_{R^\perp}\|_2^2 \sim \chi^2_{n-r}$, $\mathrm{dist}^2(\varepsilon,S) \sim \sum_{j=0}^r v^{(R)}_{\,r-j}(S)\,\chi^2_j$. 
    Using the convolution identity $\chi^2_j+\chi^2_{\,n-r}\stackrel{\textup{d}}{=}\chi^2_{j+(n-r)}$, we obtain:
    \begin{equation*}
        \mathrm{dist}^2(\varepsilon,S) \sim \sum_{j=0}^r v^{(R)}_{\,r-j}(S)\,\chi^2_{j+(n-r)}.
    \end{equation*}
    Re-indexing with $\ell=j+(n-r)$ (so $\ell=n-r,\dots,n$) and comparing with \eqref{eq:chisqbarS_rn} gives precisely the ambient-$\mathbb R^n$ mixture:
    \[
        \mathrm{dist}^2(\varepsilon,S)\ \stackrel{\textup{d}}{=}\ \sum_{\ell=0}^{n} v^{(\mathbb R^n)}_{\,n-\ell}(S)\,\chi^2_{\ell},
        \quad\text{with}\quad
        v^{(\mathbb R^n)}_{\,n-\ell}(S)=
        \begin{cases}
        v^{(R)}_{\,r-(\ell-(n-r))}(S)=v^{(R)}_{\,n-\ell}(S), & \ell\ge n-r,\\[2pt]
        0, & \ell<n-r.
        \end{cases}
    \]
    Equivalently,
    \begin{equation*}
         v_k^{(\mathbb R^n)}(S)=v_k^{(R)}(S)\ \text{for }k=0,\dots,r,\qquad v_k^{(\mathbb R^n)}(S)=0\ \text{for }k>r.
    \end{equation*}
    This establishes the claimed relationship.
\end{proof}

\begin{table}[!t]
    \centering
    \caption{Location of the maximizers of quantiles for different test statistic and constraint sets}
    \begin{tabular}{l c c}
    \toprule
    \textbf{Statistic} & \textbf{Global optimization} $(D^*)$ & \textbf{Sliced optimization} $(d^*(\mu))$\\
    \midrule
    $\lambda^2_c$ (all cases)      & General results not available & General results not available \\
    $\lambda^1,\lambda^2_u$ ($x\geq 0$)& $x = 0$ & Extreme point of $x\geq 0, Hx = \mu$ (at most $\binom{p}{k}$)\\
    $\lambda^1,\lambda^2_u$ ($Ax\leq b$)      & Extreme point of $Ax \leq b$ & Extreme point of $Ax \leq b, Hx = \mu$\\
    $\lambda^1,\lambda^2_u$ ($x \in \mathcal{C}$)       & $x = 0$ & Relative boundary of $x\in\mathcal{C}, Hx = \mu$\\
    \bottomrule
    \end{tabular}
    \label{tab:where_max_quantiles_at}
\end{table}

\begin{table}[!t]
    \centering
    \caption{Comparative summary of the test statistics.}
    \setlength{\tabcolsep}{2pt}
    \small
    \begin{tabular}{@{}L{4cm} L{4.1cm} L{4.1cm} L{4.1cm}@{}}
    \toprule
    \textbf{Property} & $\lambda^2_c$ & $\lambda^2_u$ & $\lambda^1$ \\
    \midrule
    Quantile convexity $Q_x$
    & No (in general)
    & Yes (Gaussian noise) [\Cref{thm:convexityl2u}]
    & Yes (log-concave noise) [\Cref{thm:convexity_Q1}] \\
    \arrayrulecolor{black!25}\midrule
    Acceptance region emptiness
    & Always non-empty
    & May be empty for some $y$
    & May be empty for some $y$ \\
    \arrayrulecolor{black!25}\midrule
    Calibration ease
    & Typically challenging
    & Typically easier
    & Typically easier \\
    \arrayrulecolor{black!25}\midrule
    Calibration status (non-negative constraints)
    & Optimal $D^*$ unknown
    & Optimal $D^*$ known [\Cref{cor:improved_constants_2u}]
    & Optimal $D^*$ known [\Cref{cor:improved_constants_1}] \\
    \arrayrulecolor{black!25}\midrule
    Average vs.\ worst-case length (small $\alpha$)
    & \cellcolor{lightgray!25}
    & Often shorter on average; worse worst-case than $\lambda^1$ (random radius)
    & Often longer on average; better worst-case than $\lambda^2_u$ (fixed radius) \\
    \arrayrulecolor{black}\bottomrule
    \end{tabular}
    \label{tab:adv_disadv_compact_matrix}
\end{table}

We summarize the results of this section regarding the test statistics $\lambda^1$, $\lambda^2_u$, and $\lambda^2_c$ in \Cref{tab:where_max_quantiles_at,tab:adv_disadv_compact_matrix}. \Cref{tab:where_max_quantiles_at} summarizes the results of quantile maximizations, including \Cref{thm:maximizing_2u}, \Cref{thm:maximizing_2u_cone}, and \Cref{thm:maximizing_one_term}. \Cref{tab:adv_disadv_compact_matrix} provides a summary of the advantages and disadvantages of using each test statistic.

\section{Practical methods for high-dimensional problems}
\label{sec:high-dimensional}

This section uses the results developed in \Cref{sec:test_statistics} to propose two approaches that scale to high-dimensional problems: 
(i) reductions for cases in which $K$ has full column rank\footnote{Depending on the specific alignment of $H$ and $K$, reductions might also be used for some non-full column rank $K$, as we explain below}, 
and 
(ii) a split technique that separates a full column rank part and a null space part for cases in which $K$ is rank-deficient. 

\subsection{Generalizing and improving TFM reductions}
\label{sec:reductions}

A TFM reduction (which we name after Tenorio, Fleck, and Moses, who first introduced them in \cite{tenorio2007confidence}) maps a broad class of problems to a single canonical form.
In their original formulation, the family
\begin{equation*}
   y = Kx^* + \varepsilon, \quad \varepsilon \sim \mathcal{N}(0, I), \quad x \geq 0, \quad K \in \mathbb{R}^{n\times p}, \quad h \in \mathbb{R}^p,
\end{equation*}
with arbitrary $K$ and a scalar functionals $h^\top x^*$, is reduced to a two-dimensional problem:
\begin{equation*}
  y = x^* + \varepsilon \in \mathbb{R}^2, 
   \quad \varepsilon \sim \mathcal{N}(0,\Sigma), \quad x^* \geq 0,
   \quad h = (1,-1),
\end{equation*}
which is independent of $(K,h)$ except through the covariance $\Sigma$ of the transformed noise. In this example, we allow $\Sigma$ to depend on $(K,h)$, even though the original proposal never requires computing it explicitly. The problem might be further reduced to one dimension if $h \geq 0$ or $h \leq 0$, as we formalize below.
Such reductions lead to computationally inexpensive algorithms that remain valid in high dimensions, at the cost of discarding some constraint information during the reduction. In what follows, we first recall the original reduction of \cite{tenorio2007confidence}, then use our results in \Cref{sec:test_statistics} to improve it and to generalize it to multiple functionals and box constraints.

\subsubsection[Review and improvements to TFM reduction for the single-functional case]{Review and improvements to TFM reduction for the single-functional case}

Let $K$ be a full column rank matrix, and consider $k = 1$. 
We consider two cases depending on $h$.

\emph{When $h \geq 0$ or $h \leq 0$.} Without loss of generality, suppose $h \geq 0$. 
Given $y = Kx^* + \varepsilon$, the reduction in \cite{tenorio2007confidence} constructs the one-dimensional vector:
\begin{equation}
    \label{eq:tenorio_og_1d}
    \tilde{y} =
         h^\top(K^\top K)^{-1} K^\top y
    = 
         h^\top x^*
    + \tilde{\varepsilon} =: \tilde{x}^* + \tilde{\varepsilon},
\end{equation}
with $\tilde{x}^* = h^\top x^*$ and $\tilde \varepsilon \sim \mathcal{N}(0,\sigma^2)$ with $\sigma^2 = h^{\top}(K^\top K)^{-1}h$.
Observe that $\tilde{x}^* \geq 0$.
Thus, this yields a one-dimensional constrained problem.
The original paper then builds the interval:
\begin{equation}
    \label{eq:og_tenorio_reduction_1d}
    \begin{aligned}
        \min_{\tilde{x}}/\max_{\tilde{x}} \quad & \tilde{x}\\
        \st \quad  &  \frac{1}{\sigma^2}(\tilde{x}-y)^2 \leq \frac{1}{\sigma^2}\min_{\tilde{x}'\geq 0} (\tilde{x}'-y)^2 +  Q_{\chi^2_1,1-\alpha}\\ 
        & \tilde{x} \geq 0,
    \end{aligned}
\end{equation}
which is a valid $1-\alpha$ confidence interval for the original $h^\top x^*$.

\emph{When $h$ is general.}
Decompose $h = h^+ - h^-$ with $h^+, h^- \geq 0$. 

Given $y = Kx^* + \varepsilon$, the reduction in \cite{tenorio2007confidence} constructs the two-dimensional vector:
\begin{equation}
    \label{eq:tenorio_og_2d}
    \tilde{y} 
    =
    \begin{bmatrix}
         h^{+,\top}(K^\top K)^{-1}K^\top y  \\
         h^{-,\top}(K^\top K)^{-1}K^\top y  \\
     \end{bmatrix}
     = 
    \begin{bmatrix}
         h^{+,\top}x^*  \\
         h^{-,\top}x^*  \\
    \end{bmatrix} 
    + \tilde{\varepsilon} =: \tilde{x}^* + \tilde{\varepsilon}
\end{equation}
with $(1,-1)^\top\tilde{x}^* = h^\top x^*$ and $\tilde \varepsilon \sim \mathcal{N}(0,\Sigma)$ with 
\begin{equation*}
    \Sigma
    = 
    \begin{pmatrix}
    h^{+,\top}(K^\top K)^{-1}h^+ & h^{+,\top}(K^\top K)^{-1}h^- \\[6pt]
    h^{-,\top}(K^\top K)^{-1}h^+ & h^{-,\top}(K^\top K)^{-1}h^-
    \end{pmatrix}.
\end{equation*}
Observe that $\tilde{x}^* \geq 0$.
The original paper then builds the interval:
\begin{equation}
    \label{eq:og_tenorio_reduction_2d}
    \begin{aligned}
        \min_{\tilde{x}}/\max_{\tilde{x}} \quad & (1, -1)^\top \tilde{x}\\
        \st \quad  & \Vert \tilde{x}-y \Vert_\Sigma^2 \leq \min_{\tilde{x}'\geq 0} \|\tilde{x}'-y\|^2_\Sigma +  Q_{\chi^2_2,1-\alpha}\\ 
        & \tilde{x} \geq 0,
    \end{aligned}
\end{equation}
where $\|v\|_\Sigma = v^\top \Sigma^{-1}v$.
This is a valid $1-\alpha$ confidence interval for the original $h^\top x^*$.

It is important to note that some information is lost in these reductions.
While it is true that $x^\star \geq 0$ implies $\tilde{x}^* \geq 0$, the reverse is not necessarily true.
Therefore, we have weakened the constraint from $\{x : x\geq 0\}$ to $\{x: \min(h^{+, \top}x, h^{-, \top}x) \geq 0\}$. 
Consequently, we expect to pay a price in interval length and overcoverage in exchange for the reduction that simplified the computations.

Both \eqref{eq:og_tenorio_reduction_1d} and \eqref{eq:og_tenorio_reduction_2d} can be identified as $\lambda^2_c$ test inversion regions with thresholds taken from the $\chi^2$ bounds in \eqref{eq:chain_stochastic_dominance}, which are not tight in general.
Moreover, in \eqref{eq:og_tenorio_reduction_2d}, the term $\min_{x'\geq 0} \|x'-y\|^2_\Sigma$  is unnecessary for coverage, since $\{x : \|x-y\|^2_\Sigma \leq  Q_{\chi^2_2,1-\alpha} \}$ is already a $1-\alpha$ confidence set for $x$. 
Two types of refinements can be done, depending on whether we wish to keep the term $\min_{\tilde{x}'\geq 0} \|\tilde{x}'-y\|^2_\Sigma$, which we briefly discuss below.

If one keeps the extra term, we work with region coming from $\lambda^2_c$ inversion for analysis.
This is possible for the one-dimensional $h \geq 0$ or $h \leq 0$ case, but becomes challenging in two dimensions due to the arbitrary covariance matrix $\Sigma$. 
Use of $\lambda^2_c$ for the one-dimensional problem $y = x^* + \varepsilon$ with $x^* \geq 0$ was studied in \cite[Section 2.4]{BatlleDisproof} for unit variance noise, but the results apply easily after rescaling to the non-unit variance case. 
Define
\begin{equation*}
    q_{\alpha}(\mu) = Q_{\mu, 1-\alpha} =
    \begin{cases}
    Q_{\chi^2_1,1-\alpha}, & 1 - \alpha < \chi^2_1(\mu^2), \\[6pt]
    r_{\mu,\alpha}, & 1 - \alpha \geq \chi^2_1(\mu^2),
    \end{cases}
\end{equation*}
where $r_{\mu,\alpha}$ is the unique non-negative root of the function $x \mapsto \Phi(\sqrt{x}) - \Phi\!\left(\frac{-\mu^2 - x}{2\mu}\right) - (1-\alpha)$, which can be found using numerical methods. 
We can then build the improved region $\{x : \frac{1}{\sigma^2}(y - x)^2 \leq q_{\alpha}(x) + \frac{1}{\sigma^2}\min_{x' \geq 0}(y - x')^2\}$, which is strictly included in \eqref{eq:og_tenorio_reduction_1d}.

If one drops the extra term, then one obtains regions coming from $\lambda^1$ or $\lambda^2_u$ (here $\lambda^1 = \lambda^2_u$ as $K = I$) inversion for analysis.
The constant $Q_{\chi^2_2,1-\alpha}$ can then be improved according to \Cref{cor:improved_constants_1}, with the optimal global constant being the quantile at zero, yielding the following intervals for the one-dimensional and two-dimensional reductions, respectively:
\begin{equation*}
    \begin{aligned}
        \min_{x}/\max_{x} \quad &  x \\
        \st \quad  & \frac{1}{\sigma^2}(x-y)^2  \leq q_1 \\ 
        & x \geq 0,
    \end{aligned}
\end{equation*}
and
\begin{equation*}
    \begin{aligned}
        \min_{x}/\max_{x} \quad & (1, -1)^\top x \\
        \st \quad  & \Vert x-y \Vert_\Sigma^2 \leq q_2 \\ 
        & x \geq 0,
    \end{aligned}
\end{equation*}
with 
\[
    q_1 := \frac{1}{\sigma^2}Q_{1-\alpha}\Big(\min_{z= 0} (z-\varepsilon)^2\Big)=Q_{\chi^2_1, 1-\alpha},
\]
and 
\[
    q_2 := Q_{1-\alpha}\Big(\min_{z_1 = z_2, z_1 \geq 0, z_2 \geq 0} \|z-\varepsilon\|^2_\Sigma\Big),
\]
with $\varepsilon \sim \mathcal{N}(0, \Sigma)$. 
In both of these settings, one can also show that the optimal sliced function equals the global constant $d^*(\mu) = D^*$ for all $\mu$, so no further improvement is possible by slicing.

\emph{Extension to some rank-deficient cases.} 
While we have assumed so far that $K$ has full column rank so that the method is valid for \emph{every} quantity of interest $h$, the only requirement on a given $K$ and $h$ for this construction to work is to be able to write $h = h_1 - h_2$, with $h_i \geq 0, h_i \in \row(K), i = 1,2$. 
Note that $h_1$ and $h_2$ can, but do not necessarily need to be, the positive and negative parts of the vector $h$ and that $h \in \row(K)$ is not sufficient for this condition to hold when $K$ is not full column rank.
The existence of the vectors $h_1$ and $h_2$ can be checked by solving the following linear feasibility program:
\begin{align*}
    \begin{aligned}
        \text{find } & v_1, v_2 \in \mathbb{R}^n \\
        \text{s.t. } & (v_1 - v_2)^\top K = h \\
                     & v_1^\top K \geq 0 \\
                     & v_2^\top K \geq 0.
    \end{aligned}
\end{align*}
Whenever this program has a solution, one sets $h_i = v_i^\top K$ and replaces $(K^\top K)^{-1}$ by $(K^\top K)^{\dagger}$ throughout (so $h_i^\top (K^\top K)^\dagger K^\top y = h_i^\top K^\dagger y$).

\subsubsection[Generalizing TFM reduction to the multiple-functional case]{Generalizing TFM reduction to the multiple-functional case}

Consider now $H \in \mathbb{R}^{k \times p}$ and $K$ having full column rank.
Let us assume that out of these $k$ functionals, $k_1 \leq k$ have positive and negative parts, and $k_2 = k-k_1$ are non-negative (or non-positive, but we consider them non-negative without loss of generality), and sort $H$ so that the first $k_1$ are the ones with positive and negative parts. After observing $y = Kx+\varepsilon$, form the $2k_1+k_2$ dimensional reduced observation:
\begin{equation}
    \tilde{y} 
    = 
    \begin{bmatrix}
         h_1^{+,\top}(K^\top K)^{-1}K^\top y\\
        h_1^{-,\top}  (K^\top K)^{-1}K^\top y\\
        \cdots\\
         h_{k_1}^{+,\top}(K^\top K)^{-1}K^\top y\\
        h_{k_1}^{-,\top}(K^\top K)^{-1}K^\top y  \\
        \cdots\\
         h_{k_1 + 1}^{\top}(K^\top K)^{-1}K^\top y \\
        \cdots\\
         h_{k}^{\top}(K^\top K)^{-1}K^\top y
     \end{bmatrix} 
     = 
     \begin{bmatrix}
         h_1^{+,\top}x^*\\
        h_1^{-,\top}x^*  \\
        \cdots\\
         h_{k_1}^{+,\top}x^*\\
        h_{k_1}^{-,\top}x^*  \\
        \cdots\\
         h_{k_1 + 1}^{\top}x^* \\
        \cdots\\
         h_{k}^{\top}x^*
     \end{bmatrix} 
     + \tilde{\varepsilon} =: \tilde{x}^* + \tilde{\varepsilon}.
\end{equation}
And now we have $\tilde y = \tilde x + \tilde \varepsilon$, with $\tilde x \geq 0 $ and $\tilde \varepsilon \sim \mathcal N(0,\Sigma)$. 
Letting $h^{(1)},\dots,h^{(2k_1+k_2)}$ denote the expanded collection of vectors obtained by replacing each functional $h_\ell$ by $(h_\ell^+,h_\ell^-)$ if it has mixed signs, and by $h_\ell$ itself if it is non-negative, the covariance of the reduced noise $\Sigma$ is:
\begin{equation*}
    \Sigma_{ij} = {h^{(i)}}^\top (K^\top K)^{-1} h^{(j)},\qquad i,j=1,\dots,2k_1+k_2.
\end{equation*}
The new functional of interest matrix $\tilde{H}$ equals:
\begin{equation}
    \label{eq:reduced_matrix}
    \tilde{H} =\begin{pmatrix}
    I_{k_1}\otimes (1,\,-1) & \mathbf 0_{\,k_1\times k_2} \\[4pt]
    \mathbf 0_{\,k_2\times 2k_1} & I_{k_2}
    \end{pmatrix}
     = 
    \begin{pmatrix}
    1      & -1     & 0      & \cdots & 0      & 0      & \cdots & 0      \\[6pt]
    0      & 0      & 1      & -1     & 0      & 0      & \cdots & 0      \\[6pt]
    \vdots & \vdots & \ddots & \ddots & \ddots & \vdots &        & \vdots \\[6pt]
    0      & \cdots & 0      & 0      & 1      & -1     & 0      & \cdots \\[6pt]
    \hline
    0      & \cdots & 0      & 0      & 0      & 0      & 1      & 0      \\[6pt]
    \vdots &        &        &        &        &        & \ddots & \vdots \\[6pt]
    0      & \cdots & 0      & 0      & 0      & 0      & 0      & 1
    \end{pmatrix}
    \in \mathbb R^{k\times (2k_1 + k_2)}.  
\end{equation}
In words, $\tilde H$ is block diagonal: it has a $[1,-1]$ row for each of the $k_1$ sign-split functionals, followed by $k_2$ identity rows for the non-negative functionals. 
Note that the new dimension $2k_1 + k_2$ being in between $k$ and $2k$ means that the reduction is the most attractive when $k \ll p$. 
The maximum quantile over $x \geq 0$ is at $x=0$. The random variable is $\inf_{\tilde{H}z = 0, z \geq 0} \|z-\varepsilon\|^2_\Sigma$ with $\varepsilon \sim \mathcal N(0,\Sigma)$.
Since that sets the last $k_2$ variables to $0$ directly and then every pair of variables corresponding to the same functional collapse to one, this is effectively a $k_1$-dimensional problem. 
Let $Q_0$ be the $1-\alpha$ quantile of this random variable.
Note that this is a multiple functional problem, so as in the general case, the box given by the $k$ problems:
\begin{align*}
     \min_{x}/\max_{x} \quad & \tilde{H}_i^\top x\\
      \st \quad  & \Vert x-y \Vert_\Sigma^2 \leq Q_0\\ 
    & x \geq 0,
\end{align*}
encapsulates the convex region $\{\mu \in \mathbb R^{2k_1+k_2} : \lambda(\mu, y) \leq Q_0\}$. 
Furthermore, there is at most one point in the polyhedron $\tilde{H}x = \mu, x\geq 0$, and it can be shown that the maximum over any polyhedron of the form $\tilde{H}x=\mu, x\geq 0$ equals $Q_0$.

Analogous to the $k = 1$ case, a slightly weaker condition than $K$ being full column rank is the collection of expanded functionals $h^{(i)}$ all being in the row space of $K$, in which case, one can replace the inverses by pseudoinverses. 
In the collection of expanded functionals, we allow for any split $h = h^1-h^2$ such that $h^1 \geq 0, h^2 \geq 0$ as long as $h^1,h^2\in\row(K)$ for the reduction to work.

\subsubsection{Special case of box constraints}

We next develop an explicit algorithm for the special case of box constraints that arises frequently in practice.
These constraints are $x_l \leq x \leq x_u$, where some components of $x_l$ or $x_u$ may be infinite (formally $-\infty$ or $+\infty$), so that variables may be two-sided bounded, one-sided bounded, or unconstrained.
As a first step, we make an affine change of coordinates $z = Tx + \delta$ with invertible $T$ such that variables originally bounded from above and below are mapped to $[0,1]$, variables bounded only on one side are mapped to be non-negative, and unconstrained variables remain free in $\mathbb R$. 
After this change of coordinates, we obtain
\begin{equation} 
    \label{eq:box_change_coord}
    y = KT^{-1}(z-\delta) + \varepsilon 
    \quad\Rightarrow\quad 
    y' = K'z + \varepsilon,
\end{equation}
with $y' = y+KT^{-1}\delta$, $K' = KT^{-1}$, and new constraints given by:
\begin{equation*}
    0 \leq z_i \leq 1 \quad (i=1,\dots,s), \quad 
    z_i \geq 0 \quad (i=s+1,\dots,s+t), \quad 
    z_i \in \mathbb R \quad (i=s+t+1,\dots,p).
\end{equation*}

Write the index sets $A:=\{1,\dots,s\}$ (two-sided bounded), $B:=\{s+1,\dots,s+t\}$ (one-sided bounded, $z_i\ge 0$), and $U:=\{s+t+1,\dots,p\}$ (unconstrained).
Depending on the dimensionality of the problem, it might be feasible to calibrate regions of $\lambda^1$ or $\lambda^2_u$ by finding the maximum quantile numerically. 
We therefore start by extending the results for quantile maximization, such as \eqref{thm:maximizing_2u} and \eqref{thm:maximizing_one_term} to the case in which some variables are unconstrained. 
We then use these results to develop reductions for high-dimensional cases.
\begin{lemma} 
    \label{lemma:unconstrained_problems}
    Let $x^* = (x_c^*, x_u^*)$ with $x_c^* \in \mathbb{R}^l$ and $x_u^* \in \mathbb{R}^{p-l}$, and consider constraint set $\mathcal{X}$ of the form $\{x \in \mathbb{R}^p: Ax_c \leq b\}$. 
    Then, the distributions of the random variables $Z^2_{u, x}$ and $Z^1_x$ do not depend on the value of $x_u^*$. 
    In particular, by \Cref{thm:maximizing_2u}, the quantile global maximum is attained at a point of the form $\hat{x} = (\hat{x}_c,\hat{x}_u)$, where $\hat{x}_c$ is an extreme point of $Ax_c \leq b$ viewed as a polyhedron in $\mathbb{R}^l$, and $\hat{x}_u$ is arbitrary. 
\end{lemma}
\begin{proof}
    A direct computation after defining the test statistic $\lambda^1$ with the new constraint set yields:
    \begin{equation*}
        \mathcal{T}^1_{x}(\varepsilon) = \lambda^1(Hx, Kx + \varepsilon) = \inf_{H\xi = Hx^*, A\xi_c \leq b, \xi_u \in \mathbb{R}^{p-l}} \|K(x-x^*) - \varepsilon \|^2_2.
    \end{equation*}
    Setting $(z_c, z_u) := (x_c, x_u)-(\xi_c, \xi_u)$ as the new optimization variable, we get:
    \begin{equation*}
    \mathcal{T}^1_{x}(\varepsilon) = \inf_{Hz = 0, Az_c \geq Ax^*_c-b, z_u \in \mathbb{R}^{p-l}} \|Kz - \varepsilon \|^2_2,
    \end{equation*}
    which is independent of $x_u$. 
    The same computation with the extra term $\inf_{z} \|Kz - \varepsilon\|_2^2$, which is independent of both $x_c^*$ and $x_u^*$, proves the result for $\lambda^2_u$. 
    One then shows joint convexity of $\mathcal{T}^1_{x}(\varepsilon)$ and $\mathcal{T}^2_{u, x}(\varepsilon)$ and monotonicity across recessing directions $Ay_c \leq 0$, and the same results as in \Cref{thm:maximizing_2u} can be obtained over $x_c$. 
\end{proof}
In particular, the quantile maximizer over $z$, $\hat{z}$ is one of the $2^s$ points of the form:
\begin{equation*}
    \hat{z}_i \in \{0, 1\} \quad (i=1,\dots,s), \qquad 
    \hat{z}_i = 0 \quad (i=s+1,\dots,s+t), \qquad 
    \hat{z}_i \in \mathbb R \quad (i=s+t+1,\dots,p).
\end{equation*}

We propose the following reduction in the case it is unfeasible to compute the quantiles at $2^s$ points. 
Assume the model $y' = K'z + \varepsilon$ from \eqref{eq:box_change_coord}, and rename from now on $(K', y')$ to $(K, y)$ for simplicity. 
Assume $K$ has full column rank and $k=1$. 
Split the variables by constraint type as $x = x_A + x_B + x_U$, where $x_A=(x_1,\dots,x_s,0,\dots,0)$ with $0\le x_i\le 1$, $x_B=(0,\dots,0,x_{s+1},\dots,x_{s+t},0,\dots,0)$ with $x_i\ge 0$, and $x_U=(0,\dots,0,x_{s+t+1},\dots,x_p)$ unconstrained. 

Write $h=h^+-h^-$ with $h^+,h^-\ge 0$.
After observing $y=Kx^\star+\varepsilon$, define the reduced observation:
\begin{equation}
    \tilde y =
    \begin{bmatrix}
    h^{+,\top}(K^\top K)^{-1}K^\top y\\[2pt]
    h^{-,\top}(K^\top K)^{-1}K^\top y
    \end{bmatrix}
    =
    \begin{bmatrix}
    h^{+,\top}x_A + h^{+,\top}x_B + h^{+,\top}x_U\\[2pt]
    h^{-,\top}x_A + h^{-,\top}x_B + h^{-,\top}x_U
    \end{bmatrix}
    + \tilde\varepsilon.
\end{equation}
Introduce the six-dimensional vector $\tilde x := \big(h^{+,\top}x_A,h^{+,\top}x_B,h^{+,\top}x_U,h^{-,\top}x_A,h^{-,\top}x_B,h^{-,\top}x_U\big) \in \mathbb R^6$, so that
\begin{equation}
    \tilde y =
    \underbrace{\begin{bmatrix}
    1 & 1 & 1 & 0 & 0 & 0\\
    0 & 0 & 0 & 1 & 1 & 1
    \end{bmatrix}}_{=:~\tilde K\in\mathbb R^{2\times 6}}
    \tilde x + \tilde\varepsilon,
\end{equation}
with $\tilde\varepsilon \sim \mathcal N(0,\Sigma)$ where covariance matrix $\Sigma$ given by:
$
    \Sigma =
    \begin{bmatrix}
    (h^+)^\top (K^\top K)^{-1} h^+ & (h^+)^\top (K^\top K)^{-1} h^- \\[6pt]
    (h^-)^\top (K^\top K)^{-1} h^+ & (h^-)^\top (K^\top K)^{-1} h^-
    \end{bmatrix}.
$
The functional of interest is $\tilde h=(1,1,1,-1,-1,-1)$.
The new constraints on $\tilde x$ are:
\begin{equation*}
    0 \le \tilde x_1 \le M_+^{A},\qquad \tilde x_2 \ge 0,\qquad \tilde x_3 \in \mathbb R,\qquad
    0 \le \tilde x_4 \le M_-^{A},\qquad \tilde x_5 \ge 0,\qquad \tilde x_6 \in \mathbb R,
\end{equation*}
where $M_+^{A}:=\sum_{i\in A} h_i^+$ and $M_-^{A}:=\sum_{i\in A} h_i^-$. 
Note that $\tilde h$ lies in the row space of the reduced system $\tilde K$, hence the resulting interval is finite.

For calibration, \Cref{lemma:unconstrained_problems} gives us directly that only four points need to be checked for the maximum quantiles of $Z^1_x$ and $Z^2_{u, x}$ (they are the same in this case since $K$ is surjective): $$(0, 0, 0,0,0,0), (M^A_+, 0, 0,0,0,0), (0, 0, 0,M^A_-,0,0), (M^A_+, 0, 0,M^A_-,0,0),$$ where the third and sixth coordinates were set to $0$ arbitrarily. 
Let $D$ be the maximum $1-\alpha$ quantile of $\inf_{z_1 \geq -x_1, z_4 \geq -x_4, z_i \geq 0, i \in {2,3,4,6}}\|\tilde{K}z -\varepsilon\|^2_2$ with $\varepsilon \sim \mathcal{N}(0, \Sigma)$ over $x_1 \in \{0,M^A_+\}$, $x_4 \in \{0, M^A_-\}$. 
The final interval is then given by:
\begin{equation*}
    \begin{aligned}
         \min_{\tilde{x}}/\max_{\tilde{x}} \quad &\tilde{h}^\top \tilde{x} \\
          \st \quad  & \Vert \tilde{K}\tilde{x}-\tilde{y} \Vert_\Sigma^2 \leq D \\ 
        &   0 \le \tilde x_1 \le M_+^{A} \\ 
        & \tilde x_2 \ge 0 \\
        & \tilde x_3 \in \mathbb R \\
        & 0 \le \tilde x_4 \le M_-^{A} \\
        &\tilde x_5 \ge 0 \\
        &\tilde x_6 \in \mathbb R.
    \end{aligned}
\end{equation*}

The same construction extends directly to the case in which $K$ is not full column rank but $h^+$ and $h^-$ are in its row space, and also to the case of multiple functionals. 
After the affine change of variables described above, each functional $h_i$ is split into non-negative parts $h_i^+, h_i^-$ as in the previous section (if needed), and the resulting expanded system $(\tilde y,\tilde x,\tilde \varepsilon)$ has the same structure with block matrix $\tilde H$ as in \eqref{eq:reduced_matrix}. 
The same procedure of splitting $x$ into three applies, and the feasible region for $\tilde x$ is now determined by the box constraints: coordinates coming from two-sided bounded variables are restricted to $[0,1]$, those from one-sided bounds are constrained to be non-negative, and those from unconstrained variables remain free in $\mathbb R$. 

\subsection[Row and null splitting for arbitrary forward models]{Row and null splitting for arbitrary forward models $K$}

Methods such as the TFM reductions in \Cref{sec:reductions} require that positive/negative splits of each functional lie in $\row(K)$, which holds if $K$ has full column rank. 
To handle high-dimensional problems where this fails, we develop a \say{split} technique to decompose the problem into two components: a full column rank component and a null space part. 
This construction is completely general and makes no rank assumption on $K \in \mathbb{R}^{n\times p}$.

Let $r \leq \min(n,p)$ be the rank of $K$ and let $K^\dagger$ be its Moore-Penrose pseudoinverse. 
Set $M:=K^\dagger K$, the orthogonal projector onto the row space of $K$. 
For any $H\in\mathbb R^{k\times p}$, write
\begin{equation*}
    H_\parallel := HM,\qquad H_\perp:=H(I-M),\qquad\text{thus } Hx=H_\parallel x+H_\perp x.
\end{equation*}
We construct a $1-\alpha$ confidence region for $Hx$ as a Minkowski sum:
\begin{equation*}
    \mathcal R_\alpha(y)\ =\ \mathcal R_{\parallel,\,\alpha_1}(y)\ \oplus\ \mathcal R_{\perp,\,\alpha_2}(y),\qquad \alpha_1+\alpha_2=\alpha,
\end{equation*}
where $\oplus$ denotes the (set) Minkowski sum. 
The default split is $\alpha_1 = \alpha_2 = \alpha/2$. 

\emph{Row space part ($\parallel$).}
We have
\begin{equation*}
    HK^\dagger y \sim \mathcal N\!\big(H_\parallel x,\Sigma_\parallel\big),
    \qquad 
    \Sigma_\parallel := H K^\dagger (K^\dagger)^\top H^\top,
\end{equation*}
together with the constraint $x \geq 0$. Observe that $\Sigma_\parallel = (H K^\dagger)(H K^\dagger)^\top$ is positive semidefinite of rank at most $r=\mathrm{rank}(K)$. 
If $\mathrm{rank}(H K^\dagger)=k$, then $\Sigma_\parallel$ is positive definite and the natural confidence set that ignores the constraints
\begin{equation}
    \label{eq:ellipsoid_unc_parallel_part}
    \big\{\mu\in\mathbb R^k:\ \|\mu - HK^\dagger y\|_{\Sigma_\parallel}^2\le\ Q_{\chi^2_{k,\,1-\alpha_1}}\big\}.
\end{equation}
is a bounded ellipsoid. 
If $\mathrm{rank}(H K^\dagger)<k$, $\Sigma_\parallel$ is singular and the set above should be interpreted using the Moore-Penrose pseudoinverse $\Sigma_\parallel^\dagger$; the set is then unbounded along $\ker(\Sigma_\parallel)$, corresponding to redundant functionals that vanish on $\row(K)$.

Different options exist to build the $1-\alpha_1$ confidence region for the row space part. 
First and simplest, one can ignore the constraint $Ax \leq b$ and report directly \eqref{eq:ellipsoid_unc_parallel_part}. 
This unconstrained part also admits the formulation \eqref{eq:unconstrained_lambda1}, which in this case gives 
\begin{equation}
    \label{eq:ellipsoid_unc_parallel_part_2}
    \big\{\mu\in\mathbb R^k:\ \|\mu - HK^\dagger y\|^2_{\Sigma_\parallel}\ \le\ Q_{\chi^2_{n,\,1-\alpha_1}} - \min_{x'} \|Kx'-y\|^2_2 \big\}.
\end{equation}
To add the constraint, one can view this as a constrained problem of the form $\tilde{y} = \tilde{K}x + \tilde{\varepsilon}$
with $\tilde{K} = H_\parallel$, $\tilde{H} = \tilde{K}$, and $\tilde \varepsilon \sim \mathcal{N}(0, \Sigma_\parallel)$. 
In this particular case, where the functionals of interests exactly match the observation matrix, the test inversion formalism yields, after a calculation, the intersection of \eqref{eq:ellipsoid_unc_parallel_part} with $\{H_\parallel x : Ax\leq b\}$. This may lead to a possible improvement in the volume of the confidence region, depending on the geometry of the problem. If the row-space part satisfies the feasibility conditions for a TFM reduction, then this component may likewise be treated using a TFM-based approach. 

\emph{Null space part ($\perp$).}
Whether the null space part will be bounded depends on the strength of the constraints. 
One can build the strict bounds region:
\begin{equation*}
    \big\{ H_\perp x : x \in \mathbb{R}^p, Ax \leq b,  \|Kx-y\|^2_2\leq Q_{\chi^2_n, 1-\alpha_2}\} = \{\mu \in \mathbb{R}^k : \inf_{H_\perp x = \mu, Ax \leq b} \|Kx-y\|^2_2 \leq Q_{\chi^2_n, 1-\alpha_2}\big\},  
\end{equation*}
where the constant $Q_{\chi^2_n, 1-\alpha_2}$ can be improved, according to \Cref{cor:improved_constants_1}, to the $(1-\alpha_2)$ quantile of $\inf_{H_\perp z = 0, z\geq 0 }\|Kz - \varepsilon\|_2^2$ with $\varepsilon \sim \mathcal N(0, I)$.

\emph{Testing set membership.}
With $\alpha_1+\alpha_2=\alpha$, $\mathcal R_\alpha(y)\ :=\ \mathcal R_{\parallel,\,\alpha_1}(y)\ \oplus\ \mathcal R_{\perp,\,\alpha_2}(y)$ is a $1-\alpha$ confidence region for $Hx$.

An important question is whether it is computationally easy to check whether a particular point $\mu$ belongs to the Minkowski sum set $A \oplus B$. 
If $\mathcal{R}_\parallel$ and $\mathcal{R}_\perp$ are given by expressions of the form $\mathcal{R}_{\parallel, \perp} = \{\mu : g_{\parallel, \perp}(\mu, y) \leq 0\}$ for convex, easy-to-evaluate functions $g_\parallel$ and $g_\perp$ (as is the case for all the TFM reductions and our global methods), then the Minkowski sum is convex and testing membership for for a given $\mu$ is equivalent to solving the following convex feasibility program:
\begin{equation*}
    \text{find } a \quad \text{such that} \quad g_\parallel(a, y) \leq 0, \quad  g_\perp(\mu-a, y)\leq 0. 
\end{equation*}
If $g_\parallel(\mu) = \inf_{H_1x = \mu, x \geq 0} \|Kx-y\|^2_2 - d_1(y)$ and $g_\perp(\mu) = \inf_{H_2\tilde{x} = \mu, \tilde{x} \geq 0} \|\tilde{K}x-\tilde{y}\|^2_2 - d_2(\tilde{y})$, then to check whether $\mu$ is in the Minkowski sum of the two sets, we can solve the following joint feasibility convex quadratic program:
\begin{equation*}
    \begin{aligned}
    \text{find } & a,\, b,\, x,\, \tilde{x} \\
    \text{s.t. } & a + b = \mu \\
    & Hx = a, \quad x \ge 0, \quad \|Kx - y\|^2_2 \le d_1(y) \\
    & \tilde{H}\tilde{x} = b, \quad \tilde{x} \ge 0, \quad \|\tilde{K}\tilde{x} - \tilde{y}\|^2_2 \le d_2(\tilde{y}).
    \end{aligned}
\end{equation*}
This region can also be embedded into a hyper-rectangle by solving for the minimum and maximum of $h_i^\top x$ for $i = 1,\dots,k$ over this convex, easy-to-describe region.

\section{Numerical comparisons}
\label{sec:numerical}

To illustrate the methods in this work, we consider the following non-negatively constrained problem:
\begin{equation}
    \label{eq:toy_problem_setup}
    y = Kx + \varepsilon, 
    \quad \varepsilon \sim \mathcal N(0, I), 
    \quad  x^* \geq 0, 
    \quad K 
    = 
    \begin{pmatrix}
        2,1,1 \\
        0,1,1
    \end{pmatrix}, 
    \quad H 
    = 
    \begin{pmatrix}
        1,-1,0 \\
        0,1,-1
    \end{pmatrix}.
\end{equation}
Note that neither functional in $H$ lies in the row space of $K$. 
Hence, the non-negativity constraint is essential to obtain finite-area confidence regions in $\mathbb{R}^2$. 
By checking the condition in \Cref{thm:finite_or_infinite_regions}, we confirm that the only $d \geq 0$ with $Kd=0$ is the origin $d = 0$, which implies $Hd=0$ and thus guarantees finite regions under the added constraints. 
We construct $68\%$ and $95\%$ confidence regions from a single observation of $y$. We focus on regions built from $\lambda_1$ or $\lambda^2_u$ (which are identical in this example since $K$ is surjective), and we assess two improvements: 
(i) using the $\mu$-description instead of the $x$-description bounding box,
and 
(ii) using the optimal constants from \Cref{cor:improved_constants_2u,cor:improved_constants_1} instead of the naive $\chi^2$ thresholds.

\subsection{Details of competing methods}
\label{subsec:allthemethods}

Our first benchmark is given by the Simultaneous Strict Bounds (SSB) method \eqref{eq:starkSSBk}, which constructs the product of intervals:
\begin{equation}
    \label{eq:ssbk_toy_x}
    \prod_{i = 1}^2
    \quad
    \begin{aligned}
     \min_{x}/\max_{x} \quad & h_i^\top x  \\
      \st  \quad & \|Kx-y\|^2_2 \leq Q_{\chi^2_2,1-\alpha} \\ 
      & x \geq 0.
    \end{aligned}
\end{equation}
We refer to this method as SSB$_x$, since it is the $x$-description bounding box of the $\lambda^1$ inversion region.
Two improvements are immediate.
First, consider the $\mu$-description, i.e., the convex set:
\begin{equation}
    \label{eq:ssbk_toy_mu}
    \big\{ \mu \in \mathbb{R}^2: \inf_{Hx =\mu, x \geq 0} \|Kx-y\|^2_2  \leq Q_{\chi^2_2,1-\alpha}\big\}
    = 
    \big\{ Hx : x \in \mathbb{R}^3, x\geq 0, \|Kx-y\|^2_2  \leq  Q_{\chi^2_2,1-\alpha}\big\},
\end{equation}
which we call SSB$_\mu$ and which SSB$_\mu$ $\subseteq$ SSB$_x$ by construction.
Second, tighten the constant $Q_{\chi^2_2,1-\alpha}$ by replacing it (in either description) with the $1-\alpha$ quantile of $Z^1_0$ per \Cref{cor:improved_constants_1}, i.e., the quantile of $\inf_{Hz =0, z \geq 0} \|Kz-\varepsilon\|^2_2$ with $\varepsilon \sim \mathcal{N}(0,I)$. 
As noted in \Cref{sec:chi_bar_sq}, this is a chi-bar-squared distribution. 
A direct analysis shows that for this example, the distribution of $Z^1_0$ is $\frac{1}{2}\chi^2_1 + \frac 1 2 \chi^2_2$, so its quantiles can be computed by bisection on the cumulative density function. 
In higher dimensions, we estimate these quantiles by Monte Carlo sampling of $\varepsilon$. 
At $68\%$ and $95\%$, the corresponding quantiles are $1.644$ and $5.139$, compared to $2.279$ and $5.991$ for the $\chi^2_2$ distribution. 
We denote the methods with improved constants in \eqref{eq:ssbk_toy_x} and \eqref{eq:ssbk_toy_mu} as QuantileZero$_x$ and QuantileZero$_\mu$, respectively.

We also compare against a Bonferroni-corrected product of intervals by constructing marginal intervals for $h_1^\top x$ and $h_2^\top x$ using a Bonferroni-corrected quantiles at zero of $\inf_{h_i^\top z =0, z \geq 0} \|Kz-\varepsilon\|^2_2, i = 1,2$ to achieve joint coverage at levels $68\%$ and $95\%$.

Finally, we apply the split technique of \Cref{sec:reductions} by decomposing $H=H_\parallel+H_\perp$ using the projectors onto the row and null spaces of $K$, yielding
\begin{equation*}
    H_\parallel=\begin{pmatrix}1&-\tfrac12&-\tfrac12\\0&0&0\end{pmatrix},
    \qquad 
    H_\perp=\begin{pmatrix}0&-\tfrac12&\tfrac12\\0&1&-1\end{pmatrix}.
\end{equation*}
The row space contribution reduces to $HK^\dagger y\sim\mathcal N(H_\parallel x,\Sigma_\parallel)$ with $\Sigma_\parallel=\mathrm{diag}(5/4,0)$, so only the first component is noisy: $(HK^\dagger y)_1=\tfrac12 y_1-y_2\sim\mathcal N((H_\parallel x)_1,5/4)$, while the second component vanishes deterministically. 
This yields a simple $(1-\alpha_1)$ Gaussian interval for $(H_\parallel x)_1$.
The null space contribution enforces $x_2=x_3$, so feasible $x$ take the form $(x_1,t,t)$ with $x_1,t\ge0$. 
A $(1-\alpha_2)$ feasible set for $H_\perp x$ is obtained as $\{H_\perp x: x\geq 0, \|Kx-y\|_2^2 \leq R_{1-\alpha_2}\}$, with $R_{1-\alpha_2}$ chosen either as a $\chi^2_2$ quantile (which we call the naive version) or the quantile at zero of a $\lambda^1$ test statistic (which we call refined). 
Combining both parts with $\alpha_1+\alpha_2=\alpha$, the final $1-\alpha$ confidence region is the Minkowski sum $\mathcal R_\alpha(y)=\mathcal R_{\parallel,\alpha_1}(y)\oplus\mathcal R_{\perp,\alpha_2}(y)$, where the first factor is a one-dimensional segment and the second a line-constrained residual set. For this example, we take $\alpha_1 = \alpha_2 = \alpha/2$, though this split could be optimized.

\subsection[Visualization of confidence regions for fixed y]{Visualization of confidence regions for fixed $y$}
\label{subsec:visualize}

\begin{figure}[!t]
    \centering
    \begin{subfigure}{0.49\linewidth}
        \centering
        \includegraphics[width=\linewidth]{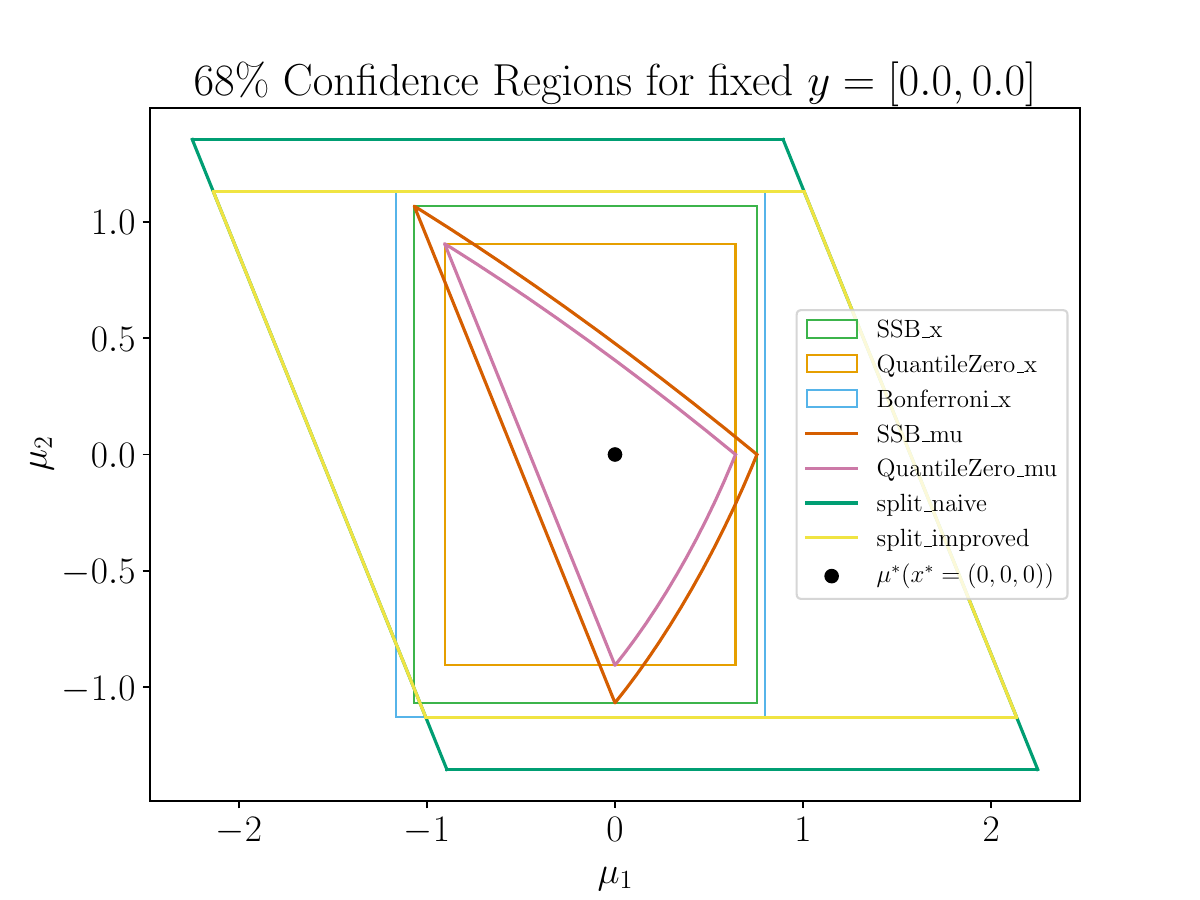}
    \end{subfigure}\hfill
    \begin{subfigure}{0.49\linewidth}
        \centering
        \includegraphics[width=\linewidth]{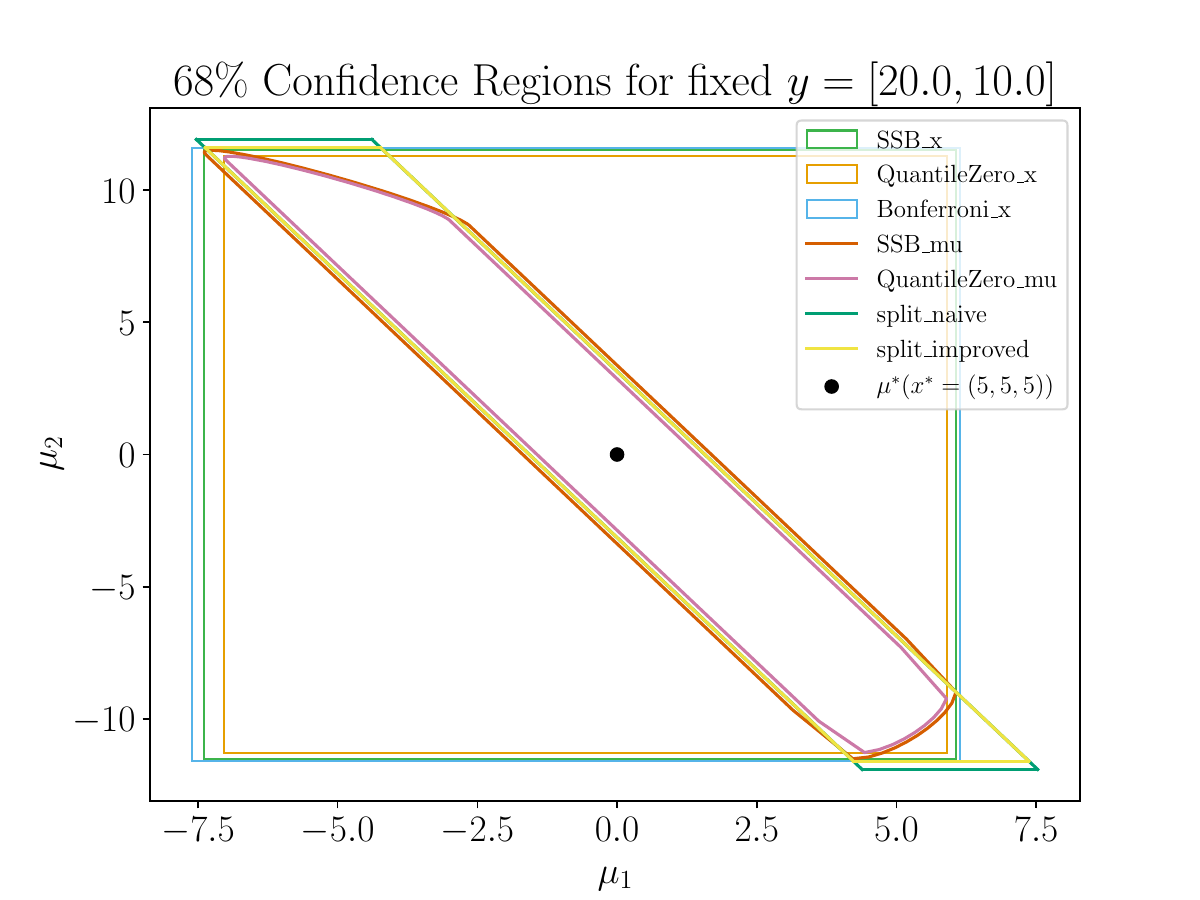}
    \end{subfigure}
    \caption{Confidence regions of at level $1-\alpha = 68\%$ for $Hx^*$ in the problem setup \eqref{eq:toy_problem_setup}, comparing the different methods in \Cref{subsec:allthemethods} when $y = (0,0)$ (left) and $y = (20,10)$ (right), which are the noiseless observations for $x^* = (0,0,0)$ and $x^* = (5,5,5)$, respectively.}
    \label{fig:confreg}
\end{figure}

We fix $y$ and compare the resulting confidence regions $\mathcal{R}_\alpha(y)$ for the different methods in \Cref{subsec:allthemethods}, deferring empirical coverage under resampling of $y$ to \Cref{subsec:areas_coverage}.
\Cref{fig:confreg} shows results for $y = (0,0)$ (left) and $y = (20,10)$ (right), the noiseless observations of $x^* = (0,0,0)$ and $x^* = (5,5,5)$. 

We observe how the $\mu-$descriptions produce non-rectangular convex sets whose bounding boxes are precisely the $x-$description bounding boxes, hence they occupy smaller area.
Furthermore, the tightening the constant from $\chi^2_2$ to the quantile at the origin shrinks both the bounding boxes and the $\mu-$descriptions. 
The Bonferroni method is additionally conservative due to both the correction and its restriction to axis-aligned products of intervals.
The split methods are Minkowski sums of two line segments in $\mathbb{R}^2$; they are fairly conservative when $y=(0,0)$ in area, but they match the $\mu$-descriptions more closely when $y=(20,10)$.

\subsection{Comparison of areas and coverages}
\label{subsec:areas_coverage}

\begin{figure}[!t]
    \centering
    \includegraphics[width=0.99\linewidth]{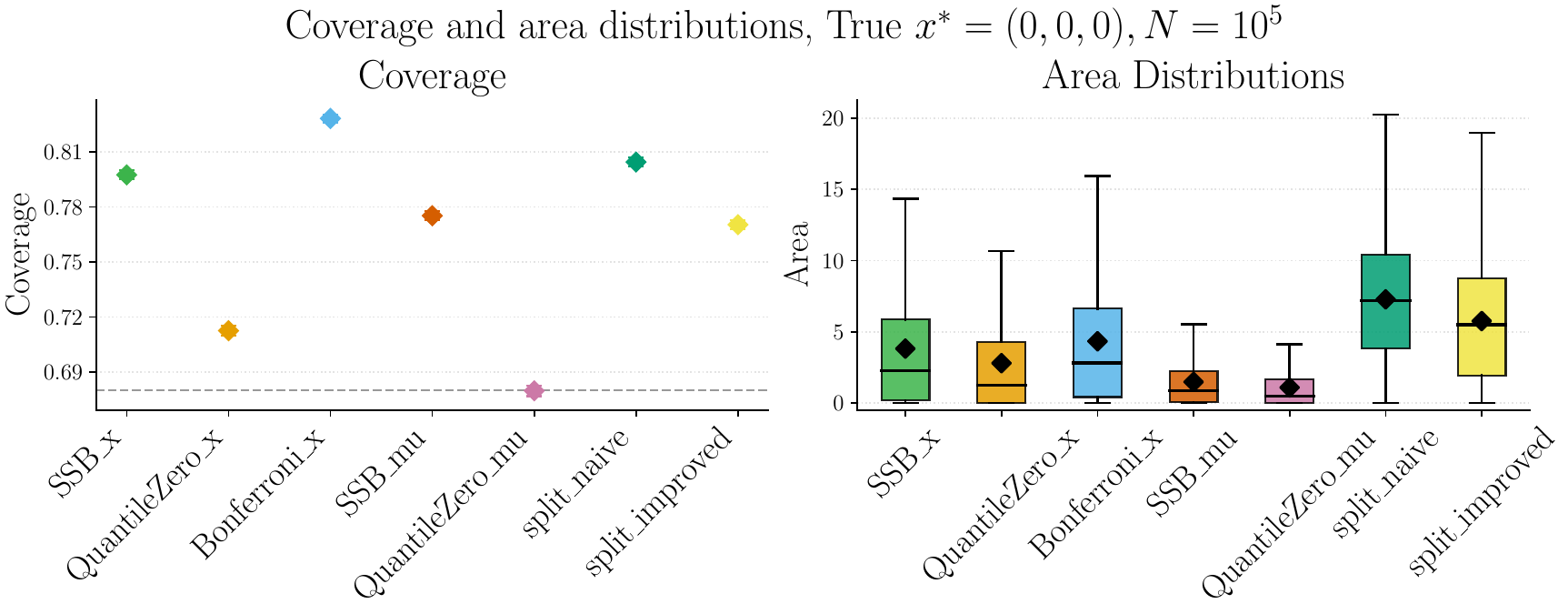} 
    \includegraphics[width=0.99\linewidth]{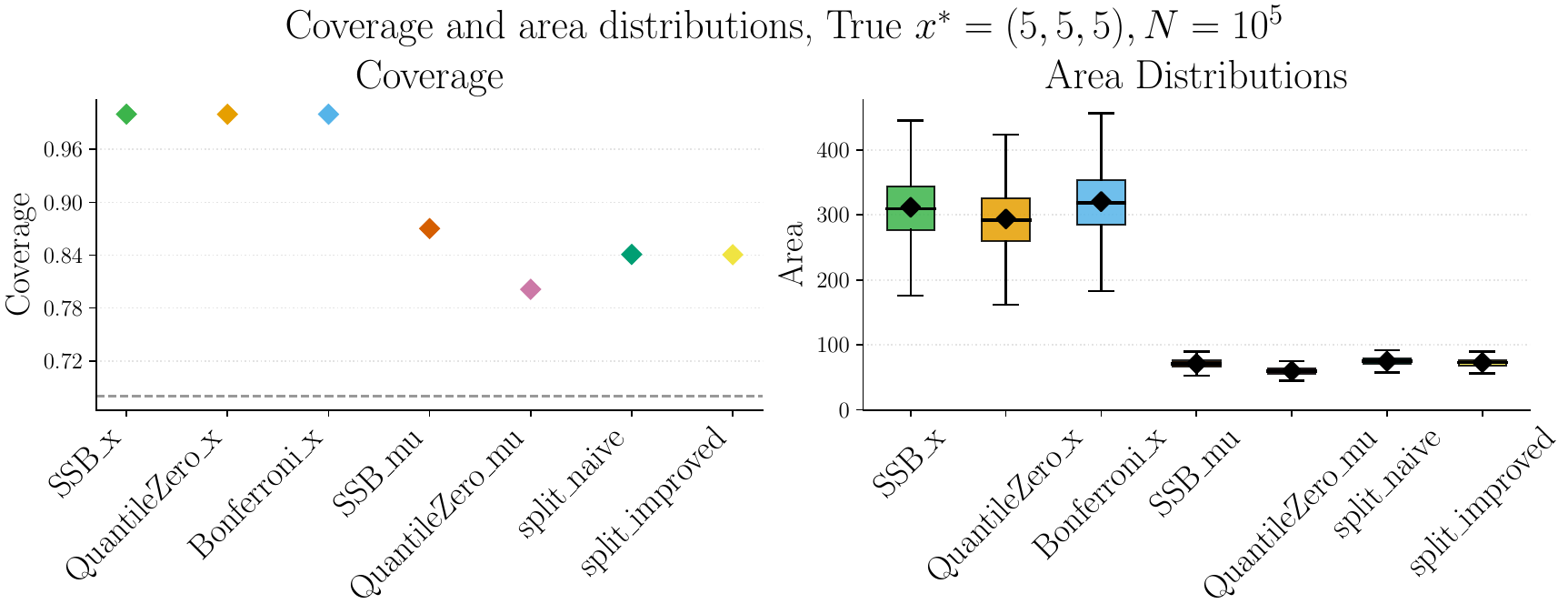}
    \caption{Coverage probabilities and area distributions of the constructed confidence sets for two scenarios: (top) $x^\star = (0,0,0)$, and (bottom) $x^\star = (5,5,5)$. Each panel shows the empirical coverage (left) and the distribution of the areas (right) for different methods, with $N = 10^5$ samples of $y \sim \mathcal{N}(Kx^*, I)$. In each plot, the diamond shows the empirical average.}
    \label{fig:100k}
\end{figure}

We evaluate empirical coverage and area distributions of different methods that construct $\mathcal{R}_\alpha(y)$ for $1-\alpha = 0.68$, by sampling $y \sim \mathcal{N}(Kx^*, I)$ a total of $N = 10^5$ times and constructing all regions for $x^* = (0,0,0)$ and $x^* = (5,5,5)$.
For each region and $y$, we compute its area using polar quadrature numerical integration, and check if they contain $Hx^*$. 
\Cref{fig:100k} reports empirical coverage (left) and area distributions (right) for $x^* = (0,0,0)$ (top) and $x^* = (5,5,5)$ (bottom). 

We observe that results align well with our theoretical predictions and the fixed-$y$ comparison in \Cref{subsec:visualize}. 
Because the $\lambda^1$ quantile is maximized at $x^* = (0,0,0)$, the $\mu$-description with the correct constant (quantile at the origin) attains near-exact coverage at this point and also achieves the smallest average area. 
For $x^* > 0$, including our example of $(5,5,5)$, this method overcovers, since the test statistic is not pivotal and the quantile outside of the origin is smaller. 
When $x^*$ is away from the origin, areas increase substantially for all methods, but the $\mu$-descriptions retain a sizable advantage in average area over their $x$-description counterparts.
We expect this improvement in terms of area to also hold and be typically magnified for higher-dimensional problems as the number of functionals $k$ grows.

\section{Discussion}
\label{sec:conclusion}

In this paper, we introduced a unified test–inversion framework for confidence regions in linear inverse problems with linear or conic constraints.
Our analysis ties three statistics, $\lambda^2_c$, $\lambda^2_u$, and $\lambda^1$, to concrete optimization templates and calibration rules.
We established convexity and recession–monotonicity for $\lambda^2_u$ and $\lambda^1$, identified where their worst–case quantiles occur, and gave tight stochastic dominance relations.
Over closed convex cones, we proved that the maximizing quantiles for $\lambda^2_u$ and $\lambda^1$ occur at $x=0$ and admit $\bar\chi^2$ laws with linked intrinsic–volume weights.
We developed practical high–dimensional methods, including refined TFM reductions, multi–functional generalizations, box–constraint reductions, and a row/null splitting strategy.
We supplied optimal global thresholds for nonnegative problems for $\lambda^2_u$ and $\lambda^1$, improving classical $\chi^2$ bounds and Burrus–style intervals.

The convex-analytic parts of our framework extend beyond Gaussian noise. 
In particular, for general log-concave likelihoods, by taking test statistics to be based on the negative log likelihood $-2 \log p(y|x)$ as in \cite{BatlleDisproof}, the statistic $\lambda^1$ and its induced quantile remain convex, since its proof of convexity relies only on log-concavity and Pr\'{e}kopa’s theorem. 
By contrast, $\lambda^2_u$ is no longer convex in general, since the Gaussian-specific structure of least squares projections and the Pythagorean decomposition are lost. 
Thus, in the log-concave setting, one is restricted to $\lambda^1$ for tractable convex confidence regions. 
The structural results concerning cones, linear constraints, monotonicity along recession directions, and boundedness, however, do not depend on Gaussianity and therefore continue to hold for arbitrary log-concave noise distributions. 
What does not extend are the exact chi-square and chi-bar-squared distributional laws, whose intrinsic-volume decompositions rely on Gaussian orthogonal invariance. 
Nevertheless, the distributions at $x=0$ that maximize the quantiles for closed convex cones can still be sampled by solving convex programs.

Although this paper concentrates on the multiple-functional case, there remain many interesting future research directions for the single-functional case.
We focused on test statistics amenable to closed-form analysis of their distributional properties, but it is possible to directly attack the calibration of the original $\lambda_c^2$ statistic by solving the related chance-constrained problem (CCP) specified in \cite{BatlleDisproof}.
CCP's are generally difficult to solve, but there are a variety of computational approaches one can take.
Scenario approaches such as those explored in \cite{Nemirovski2006, campi_garatti_2011} provide one computational approach.
``Scenarios'' in this context refer to realizations of the underlying random process, which is known in the theoretical setup of this paper.
The realized scenarios are then used to approximate the true CCP and theoretical results provide certainty that the result of the scenario optimization is feasible to the original.
Since scenarios add constraints to the approximate optimization, there is a trade-off between the certainty in the solution feasibility and the optimality.
\cite{campi_garatti_2011} provides an interesting approach to characterize this relationship.
A second possible computational approach is to use a convex approximation of the CCP.
Since the CCP for this problem is essentially that of finding a largest quantile, it can be framed as a ``Value-at-Risk'' (VaR) optimization problem.
\cite{cvar_paper} considers the relaxation of VaR to ``Conditional Value-at-Risk'' (CVaR) which is convex and thus easier to optimize.
\cite{convex_approx_ccp} generalizes the notion of convex approximations to CCP and provides several computational options.

Since it is often difficult to verify or be sure that the above approaches produce a valid quantile, there is still ample motivation to discover new analytical characterizations.
The original approach to generating functional confidence sets was based on the likelihood ratio test statistic \citep{BatlleDisproof}.
However, as we have explored in this paper, it is possible to use others which come with analytical upside.
Unlike this paper which has focused on modifications of the likelihood ratio test statistic, it may be possible to consider test statistics for which the resulting CCP is more easily solvable.
For example, \cite{prekopa2003} shows that when $\varepsilon$ follows a log-concave distribution and the test statistic is quasi-concave in the parameter that solving the CCP is equivalent to solving a convex problem.
Since our Gaussian assumption on $\varepsilon$ satisfies the log-concavity criterion, only mathematical creativity stands in the way of finding a quasi-concave test statistic with good statistical properties.
Although $\lambda^1$ defined in this paper satisfies this property, we hypothesize that this test statistics is one of many possibilities.

In addition to generating altogether new test statistics with friendly mathematical properties, there is room to explore particular configurations of the test statistics presented in this paper.
For example, although the Burrus conjecture was refuted in \cite{BatlleDisproof}, we conjecture that there are \emph{some} scenarios in which the $\chi^2_1$ quantiles are valid due to the multiplicity of empirical circumstances in which the resulting intervals were valid \cite{patil2022, stanley_unfolding, stanley2025confidenceintervalsfunctionalsconstrained}.
It is possible that there exists a class of $(K, h)$ such that $\chi^2_1$ achieves finite-sample validity.
Additionally, there may be asymptotic senses in which these quantiles or quantiles from a different distribution are valid.
Situations in which one or both the number of observations ($n$) and the parameter dimension ($p$) are large could be interesting situations to investigate.
Or, since parameter settings violating the Burrus conjecture were typically found along the constraint boundary, it may be possible to prove validity additionally assuming that the true parameter is sufficiently far from the constraint boundary.

As noted in \Cref{sec:test_statistics}, even partial results about maximizing the quantile function induced by the natural log-likelihood ratio test $\lambda^2_c$ remain elusive. 
We conjecture that the correctly calibrated $\lambda^2_c$ tests could yield regions better than those obtained from $\lambda^2_u$ or $\lambda^1$, but this problem remains open. 
A deeper analysis and comparison between $\lambda^2_u$ and $\lambda^1$, in particular about their $\bar\chi^2$ distributions at the origin would also be of interest, leading to a recommendation of one over the other depending on the geometry of particular problems.

Finally, although this work has been motivated by inverse problems in the physical sciences, we believe these results may be useful to the more general statistical literature on constrained and shape-restricted methods.
For example, \cite{hengartner_stark_1995} explores finite-sample confidence envelopes for densities that are known to be monotonic or have $k$ modes relative to a positive weight function.
The confidence envelopes are defined via the strict bounds approach from \cite{Stark1992} and thus rely on simultaneous confidence sets of the distribution.
It is possible that our approach could reduce the conservatism inherent in confidence bands constructed in this way.
Furthermore, the constrained inference literature centered on the chi-bar-squared distribution could also potentially benefit from our theoretical insights, as we have treated scenarios outside of the Type A and B problems to which the approaches in \cite{Silvapulle2001} are confined.

\section*{Acknowledgments}

PB, JRL, and HO acknowledge support from the Air Force Office of Scientific Research under MURI awards number FA9550-20-1-0358 (Machine Learning and Physics-Based Modeling and Simulation), FOA-AFRL-AFOSR-2023-0004 (Mathematics of Digital Twins), the Department of Energy under award number DE-SC0023163 (SEA-CROGS: Scalable, Efficient, and Accelerated Causal Reasoning Operators, Graphs and Spikes for Earth and Embedded Systems), the National Science Foundations under award number 2425909 (Discovering the Law of Stress Transfer and Earthquake Dynamics in a Fault Network using a Computational Graph Discovery Approach) and the DoD Vannevar Bush Faculty Fellowship Program.

\bibliographystyle{plain}
\bibliography{paper/bibliography.bib}

\begin{thebibliography}{10}

\bibitem{angelopoulos2023gentle}
Anastasios~{N.} Angelopoulos and Stephen Bates.
\newblock A gentle introduction to conformal prediction and distribution-free uncertainty quantification.
\newblock {\em Foundations and Trends in Machine Learning}, 16(4):494--591, 2023.

\bibitem{BatlleDisproof}
Pau Batlle, Michael Stanley, Pratik Patil, Mikael Kuusela, and Houman Owhadi.
\newblock Optimization-based frequentist confidence intervals for functionals in constrained inverse problems: Resolving the burrus conjecture.
\newblock {\em arXiv preprint arXiv:2310.02461}, 2023.

\bibitem{Burrus1964}
Walter~Ross Burrus.
\newblock {\em Utilization of a Priori Information in the Statistical Interpretation of Measured Distribution}.
\newblock PhD thesis, The Ohio State University, 1964.

\bibitem{campi_garatti_2011}
{M.}~{C.} Campi and {S.} Garatti.
\newblock A sampling-and-discarding approach to chance-constrained optimization: Feasibility and optimality.
\newblock {\em Journal of Optimization Theory and Applications}, 148(2):257--280, 2011.

\bibitem{cornish1938moments}
Edmund~{A.} Cornish and Ronald~{A.} Fisher.
\newblock Moments and cumulants in the specification of distributions.
\newblock {\em Revue de l'Institut international de Statistique}, pages 307--320, 1938.

\bibitem{fisher1960percentile}
Ronald~{A.} Fisher and Edmund~{A.} Cornish.
\newblock The percentile points of distributions having known cumulants.
\newblock {\em Technometrics}, 2(2):209--225, 1960.

\bibitem{hengartner_stark_1995}
Nicolas~W. Hengartner and Philip~B. Stark.
\newblock Finite-sample confidence envelopes for shape-restricted densities.
\newblock {\em The Annals of Statistics}, 23(2):525--550, 1995.

\bibitem{Kibzun1995-oz}
Andrey~{I.} Kibzun and Yuri~{S.} Kan.
\newblock {\em Stochastic programming problems with probability and quantile functions}.
\newblock Wiley Interscience Series in Systems and Optimization. John Wiley \& Sons, 1995.

\bibitem{kuusela2017strict}
Mikael Kuusela and Philip~{B.} Stark.
\newblock Shape-constrained uncertainty quantification in unfolding steeply falling elementary particle spectra.
\newblock {\em Annals of Applied Statistics}, 2017.

\bibitem{lei2018distribution}
Jing Lei, Max G’Sell, Alessandro Rinaldo, Ryan~{J.} Tibshirani, and Larry Wasserman.
\newblock Distribution-free predictive inference for regression.
\newblock {\em Journal of the American Statistical Association}, 113(523):1094--1111, 2018.

\bibitem{lei2014distribution}
Jing Lei and Larry Wasserman.
\newblock Distribution-free prediction bands for non-parametric regression.
\newblock {\em Journal of the Royal Statistical Society Series B: Statistical Methodology}, 76(1):71--96, 2014.

\bibitem{li_constrained_bootstrap}
Jessie Li.
\newblock The proximal bootstrap for constrained estimators.
\newblock {\em Journal of Statistical Planning and Inference}, 236:106245, 2025.

\bibitem{miller_simultaneous}
Rupert~{G.} Miller.
\newblock {\em Simultaneous Statistical Inference}.
\newblock Springer, 1981.

\bibitem{Nemirovski2006}
Arkadi Nemirovski and Alexander Shapiro.
\newblock {\em Scenario Approximations of Chance Constraints}, pages 3--47.
\newblock Springer London, London, 2006.

\bibitem{convex_approx_ccp}
Arkadi Nemirovski and Alexander Shapiro.
\newblock Convex approximations of chance constrained programs.
\newblock {\em SIAM Journal on Optimization}, 17(4):969--996, 2007.

\bibitem{OLeary1986}
Dianne~{P.} O’Leary and Bert~{W.} Rust.
\newblock Confidence intervals for inequality-constrained least squares problems, with applications to ill-posed problems.
\newblock {\em SIAM Journal on Scientific and Statistical Computing}, 7(2):473–489, 1986.

\bibitem{patil2022}
Pratik Patil, Mikael Kuusela, and Jonathan Hobbs.
\newblock Objective frequentist uncertainty quantification for atmospheric \(\mathrm{CO}\_2\) retrievals.
\newblock {\em SIAM/ASA Journal on Uncertainty Quantification}, 10(3):827--859, 2022.

\bibitem{Prekopa1971}
A.~Pr{\'e}kopa.
\newblock Logarithmic concave measures with application to stochastic programming.
\newblock {\em Acta Sci. Math. (Szeged)}, 32:301--316, 1971.

\bibitem{Prekopa1973}
A.~Pr{\'e}kopa.
\newblock On logarithmic concave measures and functions.
\newblock {\em Acta Sci. Math. (Szeged)}, 34:335--343, 1973.

\bibitem{prekopa2003}
Andr{\'a}s Pr{\'e}kopa.
\newblock Probabilistic programming.
\newblock {\em Handbooks in Operations Research and Management Science}, 10:267--351, 2003.

\bibitem{robertson1988}
Tim Robertson, {R.}~{T.} Wright, and {R.}~{L.} Dykstra.
\newblock {\em Order Restricted Statistical Inference}.
\newblock John Wiley and Sons, 1988.

\bibitem{cvar_paper}
{R.}~Tyrrell Rockafellar and Stanislav Uryasev.
\newblock Optimization of conditional value-at-risk.
\newblock {\em Journal of Risk}, 2:21--42, 2000.

\bibitem{Rust1972-cf}
{B.}~{W.} Rust and {W.}~{R.} Burrus.
\newblock {\em Mathematical Programming and the Numerical Solution of Linear Equations}.
\newblock Elsevier Science, 1972.

\bibitem{rust1994confidence}
Bert~{W.} Rust and Dianne~{P.} O'leary.
\newblock Confidence intervals for discrete approximations to ill-posed problems.
\newblock {\em Journal of Computational and Graphical Statistics}, 3(1):67--96, 1994.

\bibitem{seberlee}
George {A.}~{F.} Seber and Alan~{J.} Lee.
\newblock {\em Linear Regression Analysis}.
\newblock Wiley, 2003.

\bibitem{Silvapulle2001}
Mervyn~{J.} Silvapulle and Pranab~{K.} Sen.
\newblock {\em Constrained Statistical Inference: Inequality, Order, and Shape Restrictions}.
\newblock Wiley, 2001.

\bibitem{stanley_thesis}
Michael Stanley.
\newblock {\em {Uncertainty Quantification for Ill-Posed Inverse Problems in the Physical Sciences: Confidence Intervals via Optimization and Sampling}}.
\newblock PhD thesis, Carnegie Mellon University, 2024.

\bibitem{stanley2025confidenceintervalsfunctionalsconstrained}
Michael Stanley, Pau Batlle, Pratik Patil, Houman Owhadi, and Mikael Kuusela.
\newblock Confidence intervals for functionals in constrained inverse problems via data-adaptive sampling-based calibration.
\newblock {\em arXiv preprint arXiv:2502.02674}, 2025.

\bibitem{stanley_unfolding}
Michael Stanley, Pratik Patil, and Mikael Kuusela.
\newblock Uncertainty quantification for wide-bin unfolding: one-at-a-time strict bounds and prior-optimized confidence intervals.
\newblock {\em Journal of Instrumentation}, 17(10):P10013, 2022.

\bibitem{Stark1992}
Philip~{B.} Stark.
\newblock Inference in infinite-dimensional inverse problems: Discretization and duality.
\newblock {\em Journal of Geophysical Research: Solid Earth}, 97(B10):14055--14082, 1992.

\bibitem{stark1987velocity}
Philip~{B.} Stark and Robert~{L.} Parker.
\newblock Rigorous velocity bounds from soft $\tau(p)$ and $x(p)$ data.
\newblock {\em Geophysical Journal International}, 89(3):987--996, 1987.

\bibitem{tenorio2007confidence}
{L.} Tenorio, {A.} Fleck, and {K.} Moses.
\newblock Confidence intervals for linear discrete inverse problems with a non-negativity constraint.
\newblock {\em Inverse Problems}, 23(2):669, 2007.

\bibitem{vovk2005algorithmic}
Vladimir Vovk, Alexander Gammerman, and Glenn Shafer.
\newblock {\em Algorithmic Learning in a Random World}.
\newblock Springer, 2005.

\bibitem{universal_inference}
Larry Wasserman, Aaditya Ramdas, and Sivaraman Balakrishnan.
\newblock Universal inference.
\newblock {\em Proceedings of the National Academy of Sciences}, 117(29):16880--16890, 2020.

\end{thebibliography}

\end{document}